\documentclass[12pt, leqno]{amsart}
\usepackage{
amsmath,
amssymb,
amsthm,
mathrsfs,
amsfonts,
ytableau,
enumitem,
comment,
upgreek,
soul,
yhmath,
mathtools,
kotex,
stackengine,
caption,
tabularx,
stackrel,
amscd,
bbm,
diagbox,tabu,stackengine}
\usepackage[breaklinks]{hyperref}
\usepackage[mathscr]{eucal}
\usepackage{tikz,tikz-cd}
\usetikzlibrary{decorations.pathreplacing,decorations.pathmorphing}
\usetikzlibrary{arrows,arrows.meta}
\usetikzlibrary{shapes,shapes.geometric,shapes.symbols}
\usetikzlibrary{matrix,calc}
\usepackage[poly,all]{xy}
\usepackage[colorinlistoftodos, textwidth = 2.3cm]{todonotes}
\usepackage{graphicx}
\usepackage{xcolor}
\usepackage{hyperref}
\hypersetup{
    colorlinks=true,
    linkcolor=blue,
    citecolor=magenta,
    urlcolor=cyan}
\usepackage[capitalise, nameinlink]{cleveref}
\ytableausetup{mathmode, boxsize=1.1em,aligntableaux=center}

\crefname{figure}{{\sc Figure}}{{\sc Figure}}
\crefname{section}{Section}{Sections}
\crefname{example}{Example}{Examples}
\crefname{theorem}{Theorem}{Theorems}
\crefname{lemma}{Lemma}{Lemmas}
\crefname{proposition}{Proposition}{Propositions}
\crefname{thm}{Theorem}{Theorems}
\crefname{lem}{Lemma}{Lemmas}
\crefname{prop}{Proposition}{Propositions}
\crefname{figure}{Figure}{Figures}
\crefname{fig}{Figure}{Figures}
\crefname{remark}{Remark}{Remarks}
\crefname{rem}{Remark}{Remarks}
\crefname{cor}{Corollary}{Corollaries}
\crefname{corollary}{Corollary}{Corollaries}
\crefname{conjecture}{Conjecture}{Conjectures}
\crefname{conj}{Conjecture}{Conjectures}
\crefname{ex}{Example}{Examples}

\newtheorem{theorem}{Theorem}[section]
\newtheorem{proposition}[theorem]{Proposition}
\newtheorem{lemma}[theorem]{Lemma}
\newtheorem{corollary}[theorem]{Corollary}

\newtheorem*{claim*}{Claim}

\theoremstyle{definition}
\newtheorem{algorithm}[theorem]{Algorithm}
\newtheorem{example}[theorem]{Example}
\newtheorem{definition}[theorem]{Definition}
\newtheorem{remark}[theorem]{Remark}

\numberwithin{equation}{section} \numberwithin{figure}{section}
\numberwithin{table}{section}

\def \hp {0.4}
\def \vp {0.5}
\def \lw {0.5mm}
\def \ccc {3.5mm}

\def \Z {\mathbb Z}
\def \C {\mathbb C}

\newcommand{\nc}{\newcommand}
\nc{\SG}{\mathfrak{S}}
\nc{\Des}{\mathrm{Des}}
\nc{\DesLR}[2]{\mathrm{Des}_{\rm #1}(#2)}
\nc{\sfRSP}{\mathsf{RSP}}
\nc{\sfB}{\mathsf{B}}
\nc{\sfRW}{\mathsf{RW}}
\nc{\tsfB}{\widetilde{\mathsf{B}}}
\nc{\omin}{\overline{\mathrm{min}}}
\nc{\omax}{\overline{\mathrm{max}}}
\nc{\inv}[1]{\mathrm{Inv}(#1)}
\nc{\invL}[1]{\mathrm{Inv}_L(#1)}
\nc{\coinv}[1]{\mathrm{Coinv}(#1)}
\nc{\coinvL}[1]{\mathrm{Coinv}_L(#1)}
\nc{\set}{\mathrm{set}}
\nc{\comp}{\mathrm{comp}}
\nc{\peak}{\mathrm{peak}}
\nc{\setLE}[1]{{L}(#1)}
\nc{\SGR}[1]{\ensuremath{\Sigma_R(#1)}}
\nc{\MSGR}[1]{\ensuremath{\Sigma'_R(#1)}}
\nc{\SGL}[1]{\ensuremath{\Sigma_L(#1)}}
\nc{\poset}[1]{\mathsf{P}_{#1}}
\nc{\Rposet}[1]{\mathsf{RP}_{#1}}
\nc{\DHTMP}[1]{M_{#1}}
\nc{\opi}{\overline{\pi}}
\nc{\ch}{\mathrm{ch}}
\nc{\tch}{\widetilde{\mathrm{ch}}}
\nc{\calR}{\mathcal{R}}
\nc{\calF}{\mathcal{F}}
\nc{\calC}{\mathcal{C}}
\nc{\calE}{\mathcal{E}}
\nc{\sfT}{\mathsf{T}}
\nc{\sfS}{\mathsf{S}}
\nc{\sfc}{\mathsf{c}}
\nc{\sfr}{\mathsf{r}}
\nc{\sfM}{\mathsf{M}}
\nc{\inc}[1]{\mathsf{inc}(#1)}
\nc{\sourceTD}[1]{\ensuremath{T^\rightarrow_{#1}}}
\nc{\sinkTD}[1]{T^\uparrow_{#1}}
\nc{\sfFD}{{\sf F}_D}
\nc{\colindexR}[1]{\mathbf{c}(#1)}
\nc{\colindexL}[1]{\mathbf{c}_{\rm L}(#1)}
\nc{\rowindexH}[1]{\mathbf{r}(#1)}
\nc{\tHindex}[1]{\mathsf{h}(#1)}
\nc{\removablenodes}{\mathsf{RN}(\alpha)}
\nc{\sign}[1]{\mathrm{sign}(#1)}
\nc{\setNLP}[1]{{NL}(#1)}
\nc{\PI}{\mathtt{P}_n(I)}
\nc{\rmc}{\mathrm{c}}
\nc{\rmr}{\mathrm{r}}
\nc{\rmt}{\mathrm{t}}
\nc{\mapf}{\boldsymbol{\Xi}}
\nc{\readingSYCTR}[1]{\mathrm{w}_{#1}}
\nc{\readingBTLR}{\mathrm{w}_{\rm BL}}
\nc{\readingLRTB}{\mathrm{w}_{\rm LT}}
\nc{\readingTBLR}{\mathrm{w}_{\rm TL}}
\nc{\readingLRBT}{\mathrm{w}_{\rm LB}}
\nc{\readingBTRL}{\mathrm{w}_{\rm BR}}
\nc{\readingRLBT}{\mathrm{w}_{\rm RB}}
\nc{\readingRLBTBT}{\mathrm{w}_{\rm RB1}}
\nc{\readingdr}{\mathrm{w}_{\rm dr}}
\nc{\sfST}[1]{\mathsf{ST}(#1)}
\nc{\SF}{\mathsf{SF}}
\nc{\ST}{\mathsf{ST}}
\nc{\CST}{\mathsf{RCST}}
\nc{\shYQS}[1]{\overline{S}_{#1}}
\nc{\PYQS}[1]{\widetilde{S}_{#1}}
\nc{\ShSYCT}{\ensuremath{\mathsf{ShSYCT}}}
\nc{\SPYCT}{\ensuremath{\mathsf{SPYCT}}}
\nc{\SPIT}{\ensuremath{\mathsf{SPIT}}}
\nc{\RShSCT}{\ensuremath{\mathsf{RShSCT}}}
\nc{\SIT}{\ensuremath{\mathsf{SIT}}}
\nc{\SRT}{\ensuremath{\mathsf{SRT}}}
\nc{\SET}{\ensuremath{\mathsf{SET}}}
\nc{\SRCT}{\ensuremath{\mathsf{SRCT}}}
\nc{\SYCT}{\ensuremath{\mathsf{SYCT}}}
\nc{\ShSIT}{\ensuremath{\mathsf{ShSIT}}}
\nc{\ShSET}{\ensuremath{\mathsf{ShSET}}}
\nc{\Peak}{\ensuremath{\mathrm{Peak}}}
\nc{\calB}{\mathcal{B}}
\nc{\PCp}{\mathsf{Pcomp}}
\nc{\Cp}{\mathsf{Comp}}
\nc{\strel}{\mathsf{st}}
\nc{\tcd}{\mathsf{cd}}
\nc{\trd}{\mathsf{rd}}
\nc{\scd}{\mathsf{sd}}
\nc{\alphamax}{\ensuremath{\alpha_{\rm max}}}
\nc{\tal}{\tilde{\alpha}}
\nc{\BDes}{\mathbf{Des}}
\nc{\Udes}{\mathrm{Des}_{\uparrow}}
\nc{\Ldes}{\mathrm{Des}_{\leftarrow}}
\nc{\QSQ}[1]{{\widetilde{Q}}_{#1}}
\nc{\bQSQ}[1]{{\overline{Q}}_{#1}}
\nc{\bESQ}[1]{{\overline{X}}_{#1}}
\nc{\calP}{\mathcal{P}}
\nc{\calT}{\mathcal{T}}
\nc{\modR}{\textbf{\bf mod-}H_n(0)}
\nc{\modHCl}{\text{\bf mod-}HCl_n(0)}
\nc{\HClmod}{HCl_n(0)\text{-\bf mod}}
\nc{\HClbmod}{HCl_\bullet(0)\text{-\bf mod}}
\nc{\modA}{\text{\bf mod-}A}
\nc{\modRb}{\text{\bf mod-}H_\bullet(0)}
\nc{\modHClb}{\text{\bf mod-}HCl_\bullet(0)}
\nc{\Rmod}{H_n(0)\text{\bf -mod}}
\nc{\Rbmod}{H_\bullet(0)\text{\bf -mod}}
\nc{\calG}{\mathcal{G}}
\nc{\tcalG}{\widetilde{\mathcal{G}}}
\nc{\Qsym}{\mathrm{QSym}}
\nc{\Nsym}{\mathrm{NSym}}
\nc{\bfF}{\mathbf{F}}
\nc{\bfI}{\mathbf{I}}
\nc{\bfT}{\mathbf{T}}
\nc{\bfrw}{\mathbf{read}}
\nc{\bfb}{\mathbf{b}}
\nc{\bfr}{\mathbf{r}}
\nc{\bfc}{\mathbf{c}}
\nc{\bfx}{\mathbf{x}}
\nc{\bfX}{\mathbf{X}}
\nc{\AlgoOne}{\Phi_{\alpha;\rho}}
\nc{\AlgoTwo}{\Psi_{\sigma;\alpha}}
\nc{\tbfF}{\widetilde{\mathbf{F}}}
\nc{\bfM}{\mathbf{M}}
\nc{\bfQ}{\mathbf{Q}}
\nc{\bfs}{\mathbf{s}}
\nc{\tM}{\widetilde{M}}
\nc{\tK}{\widetilde{K}}
\nc{\tQ}{{\ensuremath{\widetilde{Q}}}}
\nc{\tbQ}{{\ensuremath{\widetilde{\bfQ}}}}
\nc{\bfP}{\mathbf{P}}
\nc{\bfN}{\mathbf{N}}
\nc{\valley}{\mathrm{valley}}
\nc{\rmR}{\mathrm{R}}
\nc{\Ind}{\mathrm{Ind}}
\nc{\readingword}[1]{\mathrm{rw}_{#1}}
\nc{\id}{\mathrm{id}}
\nc{\Hom}{\mathrm{Hom}}
\nc{\ocalF}{\overline{\calF}}
\nc{\osfB}{\overline{\mathsf{B}}}
\nc{\tosfB}{\widetilde{\mathsf{B}}}
\nc{\PO}{{\ensuremath{(P,\omega)}}}
\nc{\rwRC}[1]{\ensuremath{\mathrm{w}_{\rm #1}}}
\nc{\RWrM}[2]{\ensuremath{\mathsf{RW}_{#1}(#2)}}
\nc{\sinkAlgo}{\mathtt{s}(T)}
\nc{\sourceE}[1]{E_{#1}}
\nc{\sourcetauE}[1]{\ensuremath{T}^\rightarrow_{#1}}
\nc{\Tsupalpha}{T^{\rm sup}_\alpha}
\nc{\Tsourcealpha}{T^{\rm source}_\alpha}
\nc{\Trowalpha}{T^{\rightarrow}_\alpha}
\nc{\Tsinkalpha}{T^{\leftarrow}_\alpha}
\nc{\Tcolalpha}{T^{\uparrow}_\alpha}
\nc{\Tspalpha}{T^{\rm sp}_\alpha}
\nc{\rmread}{\mathsf{read}}
\nc{\rbread}{\mathrm{read}}
\nc{\sourcehtauE}[1]{\hat{\tau}^\rightarrow_{#1}}
\nc{\sourceSIT}[1]{\calT_{#1}}
\nc{\sourceSPIT}[1]{\calT_{#1}}
\nc{\sinkSIT}[1]{\calT'_{#1}}
\nc{\sinkSPIT}[1]{\calT''_{#1}}
\nc{\sourceSET}[1]{\sfT_{#1}}
\nc{\sinkSET}[1]{\sfT'_{#1}}
\nc{\sinktauE}[1]{\ensuremath{\tau^\la_{#1}}}
\nc{\sinkSYCT}[1]{\ensuremath{\hat{\tau}'_{#1}}}
\nc{\sourceSYCT}[1]{\ensuremath{\hat{\tau}_{#1}}}
\nc{\mDIF}[1]{\mathcal{V}_{#1}}
\nc{\mProj}[1]{\mathcal{P}_{#1}}
\nc{\mESF}[1]{X_{#1}}
\nc{\mQSQ}[1]{\widetilde{\mathbf{Q}}_{#1}}
\nc{\mHQSQ}[1]{\mathbf{Q}_{#1}}
\nc{\mSQ}[1]{\widetilde{\mathbf{Q}}_{#1}}
\nc{\mPYQS}[1]{\widetilde{\mathbf{S}}_{#1}}
\nc{\Tsup}[1]{T^{\rm sup}_{#1}}
\nc{\umin}{\underline{\min}}
\nc{\rmInt}{\mathrm{Int}}
\nc{\rmst}{\mathrm{st}}
\nc{\classQS}[1]{\mathcal{E}(#1)}
\nc{\calEsa}{\mathcal{E}^\sigma(\alpha)}
\nc{\classRQS}[1]{\mathbf{E}(#1)}
\nc{\classYQS}[1]{\widehat{\mathcal{E}}(#1)}
\nc{\YQS}[1]{\widehat{\mathcal{S}}_{#1}}
\nc{\mYQS}[1]{\widehat{\mathbf{S}}_{#1}}
\nc{\QS}[1]{\mathcal{S}_{#1}}
\nc{\mQS}[1]{\mathbf{S}_{#1}}
\nc{\bfSsa}{\mathbf{S}^\sigma_\alpha}
\nc{\bfSsaE}{\mathbf{S}^\sigma_{\alpha,E}}
\nc{\bal}{{\boldsymbol{\upalpha}}}
\nc{\sfP}{\mathsf{P}}
\nc{\sfRP}{\mathsf{RP}}
\nc{\RWC}{\ensuremath{{\mathsf {RW}}_{\rm c}}}
\nc{\RWR}{\ensuremath{{\mathsf {RW}}_{\rm r}}}
\nc{\ra}{\rightarrow}
\nc{\readSearles}{\mathrm{rw}}
\nc{\readSRCT}{\underline{\mathsf{read}}}
\nc{\readSYCT}{\overline{\mathsf{read}}}
\nc{\readSIT}{\mathsf{read}}
\nc{\readSRT}{\mathsf{w}}
\nc{\readSET}{\mathsf{read}}
\nc{\mbalpha}{\mathbf{Y}_\alpha}
\nc{\rev}{\mathrm{r}}
\nc{\precdot}{\prec\mathrel{\mkern-3mu}\mathrel{\cdot}}
\nc{\balpha}{Y_\alpha}
\nc{\SPCT}{\mathrm{SPCT}}
\nc{\RT}{\mathrm{RT}}
\nc{\RCT}{\mathrm{RCT}}
\nc{\SRET}{\mathrm{SRET}}
\nc{\SPYCTsa}[2]{\mathrm{SPYCT}^{#1}(#2)}
\nc{\DIF}[1]{\mathfrak{S}^*_{#1}}
\nc{\RDIF}[1]{\mathcal{R}\mathfrak{S}^*_{#1}}
\nc{\mRDIF}[1]{\mathcal{R}\mathcal{V}_{#1}}
\nc{\ESF}[1]{\mathcal{E}_{#1}}
\nc{\mRESF}[1]{\mathcal{R}X_{#1}}
\nc{\RESF}[1]{\mathcal{R}\mathcal{E}_{#1}}
\nc{\posetP}[1]{\mathcal{P}_{#1}}
\nc{\posetSEE}[1]{P_{\mQS{#1}}}
\nc{\PQS}[2]{\mathcal{S}^{#1}_{#2}}
\nc{\mPQS}[2]{\mathbf{S}^{#1}_{#2}}
\nc{\YRQS}[1]{\mathcal{R}\mathcal{S}_{#1}}
\nc{\mYRQS}[1]{\mathbf{R}\mathbf{S}_{#1}}
\nc{\RQS}[1]{\mathcal{R}_{#1}}
\nc{\mRQS}[1]{\mathcal{R}\widehat{\mathbf{S}}_{#1}}
\nc{\tH}{\mathtt{H}}
\nc{\sfposet}{\mathsf{poset}}
\nc{\sfSP}{\mathsf{SP}}
\nc{\sh}{\mathrm{sh}}
\newcommand{\smallnearrow}{{\scalebox{0.5}{$\nearrow$}}}
\newcommand{\smallswarrow}{{\scalebox{0.5}{$\swarrow$}}}
\nc{\bfY}{\mathbf{Y}}
\nc{\ble}{\boldsymbol{\leq}}
\nc{\scrD}{\mathscr{D}}
\nc{\scrE}{\mathscr{E}}
\nc{\scrC}{\mathscr{C}}
\nc{\scrG}{\mathscr{G}}
\nc{\HnZmodR}{\ensuremath{\mathbf{Q}^{r}}}
\nc{\ourMP}[1]{\mathbf{M}_{#1}}
\nc{\rad}{\mathrm{rad}}
\nc{\bdI}{\boldsymbol{I}}
\nc{\soc}{\mathrm{soc}}
\nc{\coker}{\mathrm{Coker}}

\makeatletter
\def\@tocline#1#2#3#4#5#6#7{\relax
  \ifnum #1>\c@tocdepth 
  \else
    \par \addpenalty\@secpenalty\addvspace{#2}%
    \begingroup \hyphenpenalty\@M
    \@ifempty{#4}{%
      \@tempdima\csname r@tocindent\number#1\endcsname\relax
    }{%
      \@tempdima#4\relax
    }%
    \parindent\z@ \leftskip#3\relax \advance\leftskip\@tempdima\relax
    \rightskip\@pnumwidth plus4em \parfillskip-\@pnumwidth
    #5\leavevmode\hskip-\@tempdima
      \ifcase #1
       \or\or \hskip 1em \or \hskip 2em \else \hskip 3em \fi%
      #6\nobreak\relax
    \dotfill\hbox to\@pnumwidth{\@tocpagenum{#7}}\par
    \nobreak
    \endgroup
  \fi}
\makeatother

\makeatletter
\newcommand{\oset}[3][0ex]{%
  \mathrel{\mathop{#3}\limits^{
    \vbox to#1{\kern-4\ex@
    \hbox{$\scriptstyle#2$}\vss}}}}
\makeatother

\nc\preceqdot{\mathrel{\ooalign{$\preceq$\cr\hidewidth\raise0.225ex\hbox{$\cdot\mkern0.5mu$}\cr}}}

\nc{\Deq}{\oset{D}{\simeq}}
\nc{\Meq}{\oset{M}{\simeq}}
\nc{\Ceq}{\oset{K}{\simeq}}
\nc{\Leq}{\oset{L}{\simeq}}

\newlength\cellsize 

\savebox2{%
\begin{picture}(12,12)
\put(0,0){\line(1,0){12}}
\put(0,0){\line(0,1){12}}
\put(12,0){\line(0,1){12}}
\put(0,12){\line(1,0){12}}
\end{picture}}

\nc\cellify[1]{\def\thearg{#1}\def\nothing{}%
\ifx\thearg\nothing\vrule width0pt height\cellsize depth0pt%
  \else\hbox to 0pt{\usebox2\hss}\fi%
  \vbox to 12\unitlength{\vss\hbox to 12\unitlength{\hss$#1$\hss}\vss}}
\nc\tableau[1]{\vtop{\let\\=\cr
\setlength\baselineskip{-12000pt}
\setlength\lineskiplimit{12000pt}
\setlength\lineskip{0pt}
\halign{&\cellify{##}\cr#1\crcr}}}
\nc\ctab[2]{\ensuremath{{#1} = \begin{array}{c}
\tableau{#2}
\end{array}
}}
\nc\ctabT[1]{\ensuremath{\begin{array}{c}
\tableau{#1}
\end{array}
}}

\makeatletter
\newcommand*\bigcdot{\mathpalette\bigcdot@{.5}}
\newcommand*\bigcdot@[2]{\mathbin{\vcenter{\hbox{\scalebox{#2}{$\m@th#1\bullet$}}}}}
\makeatother
\title[Equivalence classes of lower and upper descent weak Bruhat intervals]
{Equivalence classes of lower and upper descent weak Bruhat intervals}
\author[S.-I. Choi]{Seung-Il Choi}
\address{Center for Quantum structures in Modules and Spaces, Seoul National University, Seoul 08826, Republic of Korea}
\email{ignatioschoi@snu.ac.kr}

\author[S.-Y. Nam]{Sun-Young Nam}
\address{Department of Mathematics, Sogang University, Seoul 04107, Republic of Korea}
\email{synam.math@gmail.com}

\author[Y.-T. Oh]{Young-Tak Oh}
\address{Department of Mathematics / Institute for Mathematical and Data Sciences, Sogang University, Seoul 04107, Republic of Korea}
\email{ytoh@sogang.ac.kr}
\date{\today}
\begin{document}
\begin{abstract}
Let \(\rmInt(n)\) denote the set of nonempty left weak Bruhat intervals in the symmetric group \(\SG_n\). 
We investigate the equivalence relation \(\Deq\) on \(\rmInt(n)\), where \(I \Deq J\) if and only if there exists a descent-preserving poset isomorphism between \(I\) and \(J\). 
For each equivalence class \(C\) of \((\rmInt(n), \Deq)\), a partial order \(\preceq\) is defined by \([\sigma, \rho]_L \preceq [\sigma', \rho']_L\) if and only if \(\sigma \preceq_R \sigma'\). 
Kim--Lee--Oh (2024) showed that the poset \( (C, \preceq) \) is isomorphic to a right weak Bruhat interval.

In this paper, we focus on lower and upper descent weak Bruhat intervals, specifically those of the form \([w_0(S), \sigma]_L\) or \([\sigma, w_1(S)]_L\), where \(w_0(S)\) is the longest element in the parabolic subgroup \(\SG_S\) of \(\SG_n\), generated by \(\{s_i \mid i \in S\}\) for a subset \(S \subseteq [n-1]\), and \(w_1(S)\) is the longest element among the minimal-length representatives of left \(\SG_{[n-1] \setminus S}\)-cosets in \(\SG_n\).
We begin by providing a poset-theoretic characterization of the equivalence relation \(\Deq\). 
Using this characterization, the minimal and maximal elements within an equivalence class \(C\) are identified when \(C\) is a lower or upper descent interval. 
Under an additional condition, a detailed description of the structure of \((C, \preceq)\) is provided. 
Furthermore, for the equivalence class containing \([w_0(S), \sigma]_L\), an injective hull of \(\sfB([w_0(S), \sigma]_L)\) is given, and for the equivalence class containing \([\sigma, w_1(S)]_L\), a projective cover of \(\sfB([\sigma, w_1(S)]_L)\) is given. 
Here, \(\sfB(I)\) denotes the weak Bruhat interval module of the \(0\)-Hecke algebra associated with $I \in \rmInt(n)$.
The results obtained are applied to investigate lower descent intervals arising from quotient modules and submodules of projective indecomposable modules of the \(0\)-Hecke algebra.
\end{abstract}
\maketitle

\tableofcontents

\section{Introduction}
Weak Bruhat intervals play a crucial role not only in the combinatorics of Coxeter groups but also in the representation theory of both generic and degenerate Hecke algebras associated with these groups.

For each positive integer \(n\), let \(\SG_n\) denote the symmetric group on the set \([n] := \{1, 2, \ldots, n\}\). 
In this paper, we focus on weak Bruhat intervals in \(\SG_n\). 
We define \(\rmInt(n)\) to be the set of all nonempty weak Bruhat intervals in \(\SG_n\).
Unless explicitly stated otherwise, ``weak Bruhat interval" will always refer to a left weak Bruhat interval throughout this paper. 
A central problem concerning \(\rmInt(n)\) is the classification of its elements according to suitable equivalence relations.

We begin by introducing two significant and well-known equivalence relations.
Let \( \poset{n} \) denote the set of posets with ground set \( [n] \).
Each poset \( P \in \mathsf{P}_n \) can be naturally regarded as the labeled poset \( (P, \omega) \), where the labeling \( \omega: P \to [n] \) is given by \( \omega(i) = i \) for all \( i \in P \). 
This labeling assigns to each element its own index in the ground set \( [n] \), so that the label uniquely encodes the identity of the element. 
In this way, the poset \( P \) is treated as a labeled poset whose labeling \( \omega \) respects the underlying set structure.
Consequently, to each poset \( P\in \sfP_n \), 
one can associate the following generating function for its $P$-partitions: 
\[
K_{P} := \sum_{f: \text{$P$-partition}}
x_1^{|f^{-1}(1)|} x_2^{|f^{-1}(2)|} \cdots.
\]
The first equivalence relation arises in the context of $P$-partition generating functions.
For $P_1,P_2 \in \poset{n}$, define $P_1\Ceq P_2$ if $K_{P_1}=K_{P_2}$.
The classification of posets in $\poset{n}$ with respect to this equivalence relation remains a long-standing open problem (for example, see \cite{72Stan, 07RSW, 09MPW, 14MW, 20LW, 23AvDjMc}). 
Additionally, Bj\"{o}rner--Wachs demonstrated in \cite[Theorem 6.8]{91BW} that the intervals in $\rmInt(n)$ correspond precisely to regular posets in $\poset{n}$ under the mapping $P \mapsto \SGL{P}$, where $\SGL{P}$ denotes the set $\{\sigma \in \SG_n \mid \text{$\sigma(i) \leq \sigma(j)$ for all $i, j\in [n]$ with $i \preceq_P j$}\}$.
Through this correspondence, we can induce the equivalence relation $\Ceq$ on $\rmInt(n)$, defined by $I \Ceq J$ if $P_{I}\Ceq P_{J}$, where $P_I$ and $P_J$ denote the unique regular posets in $\poset{n}$ such that $\SGL{P_I}=I$ and $\SGL{P_J}=J$, respectively.
The classification of intervals in $\rmInt(n)$ under this equivalence relation also remains an open problem (see \cite[Section 5]{22JKLO}). 

The second equivalence relation arises in the context of $H_n(0)$-modules associated with the posets in \(\poset{n}\) 
(for the definition of the $0$-Hecke algebra $H_n(0)$, see~\cref{Subsec: modules of 0Hecke alg}). 
Let \(P \in \sfP_n\). 
In \cite[Definition 3.18]{02DHT}, Duchamp--Hivert--Thibon defined a right \(H_n(0)\)-module \(M_P\) associated with \(P\). Building on this concept, Kim--Lee--Oh introduced a left \(H_n(0)\)-module $\bfM_P$ associated with \(P\) in \cite[Definition 2.8]{23KLO} through a slight modification of the original module.
For $P_1,P_2 \in \poset{n}$, we define $P_1\Meq P_2$ if the modules $\bfM_{P_1}$ and $\bfM_{P_2}$ are isomorphic as left $H_n(0)$-modules. 
The relation $\Meq $ is particularly important since it refines $\Ceq$. 
However, the classification of posets in $\poset{n}$ with respect to $\Meq $ remains an open problem.
As noted earlier, according to \cite[Theorem 6.8]{91BW}, we can define the equivalence relation $\Meq$ on $\rmInt(n)$ 
by stating that $I \Meq J$ if 
$P_{I}\Meq P_{J}$.  
On the other hand, Jung--Kim--Lee--Oh \cite{22JKLO} associated to each left weak Bruhat interval \(I\) an $H_n(0)$-module \(\sfB(I)\), referred to as the {\em weak Bruhat interval module} associated with \(I\). By definition, \(\sfB(I)\) is equal to \(\bfM_{P_I}\), and hence \(I \Meq J\) if and only if \(\sfB(I) \cong \sfB(J)\) as $H_n(0)$-modules
(see \cref{Subsec: modules of 0Hecke alg}).
Although the set under consideration is restricted from $\poset{n}$ to 
$\rmInt(n)$, 
the classification of intervals in $\rmInt(n)$ under $\Meq$
remains an open problem as well (see \cite{22JKLO,23KLO}). 

We introduce the equivalence relation $\Deq$ on $\rmInt(n)$, which plays a central role in this study. This relation is defined as $I \Deq J$ if there exists a (left) descent-preserving poset isomorphism between $I$ and $J$. Recently, Kim--Lee--Oh showed that $\Deq$ refines $\Meq$ and coincides with $\Meq$ on the subset of weak Bruhat intervals corresponding to {\it regular Schur-labeled posets} on $[n]$. 
Furthermore, they successfully classified weak Bruhat intervals within this subset with respect to $\Deq$. 
They also conjectured that $\Deq$ and $\Meq$ are, in fact, identical (\cite[Theorem 5.5, Theorem 4.7, and Conjecture 7.2]{23KLO}). 
Recently, Yang--Yu \cite{24YY} proved that this conjecture holds for all weak Bruhat interval modules in arbitrary finite Coxeter types including type $A$.

Given an equivalence class $C$ of $(\rmInt(n), \Deq)$, define a partial order $\preceq$ on $C$ by $[\sigma,\rho]_L \preceq [\sigma',\rho']_L$ if and only if $\sigma \preceq_R \sigma'$. 
It was shown in \cite[Theorem 4.6]{23KLO} that $(C, \preceq)$ forms an interval. Specifically, with $\min C := [\sigma_0, \rho_0]_L$ and $\max C := [\sigma_1, \rho_1]_L$, the poset $(C, \preceq)$ is isomorphic to the right weak Bruhat interval $[\sigma_0, \sigma_1]_R$.
The main purpose of the present paper is to investigate the equivalence classes of weak Bruhat intervals that are in distinguished form but not necessarily regular Schur-labeled. Specifically, these intervals take the form 
\[
[w_0(S), \rho]_L \quad \text{or} \quad [\sigma, w_1(S)]_L,
\]
where \( w_0(S) \) is the longest element in the parabolic subgroup \( \SG_{S} \) of \( \SG_n \), generated by \( \{s_i \mid i \in S\} \) for each subset \( S \) of \( [n-1] \), and \( w_1(S) \) is the longest element in the set of minimal length representatives for left \( \SG_{[n-1]\setminus S} \)-cosets. 
Bj\"{o}rner--Wachs demonstrated in \cite[Theorem 6.2]{88BW} that for \( S \subseteq S' \subseteq [n-1] \), the set \( \{w \in \SG_n \mid S \subseteq \DesLR{R}{w} \subseteq S'\} \) is given by the weak Bruhat interval \([w_0(S), w_1(S')]_L \). 
Therefore, we refer to an interval of the form \( [w_0(S), \rho]_L \) as a {\em lower descent interval} and an interval of the form \( [\sigma, w_1(S)]_L \) as an {\em upper descent interval}. 
From the context of the representation theory of $0$-Hecke algebras, lower and upper descent intervals arise from suitable quotient modules and submodules from projective modules with basis \([w_0(S), w_1(S')]_L \) (see \cref{Lower and Upper Descent Weak Bruhat Intervals from Tableau-Cyclic}).  

In \cref{poset-theoretical characterization}, we provide a poset-theoretical characterization of the equivalence relation $\Deq$ on $\rmInt(n)$. 
More concretely, we show that for $I, J \in \rmInt(n)$, $I \Deq J$ if and only if one of $P_I$ and $P_J$ can be obtained from the other by repeatedly applying label changes to pairs that are comparable but not in the covering relation (\cref{characterization of the equivalence relation}).
This characterization facilitates a substantially enhanced understanding of the equivalence relation within a combinatorial framework.

In \cref{Sec: The equivalence classes of weak Bruhat Intervals}, we investigate the equivalence class \( C \) of a lower or upper descent interval. 

In \cref{Sec41: The canonical diagram posets}, given a diagram \( D \) with \( n \) boxes, we introduce two posets, \( P_{F_D^\downarrow} \) and \( P_{F_D^\rightarrow} \), referred to as the \emph{canonical diagram posets} associated with \( D \). 
In \cite[Section 5.1]{24CKO}, it is shown that every upper descent interval appears as \( P_{F_D^\rightarrow} \) for some diagram \( D \). 
Here, we show that every lower descent interval appears as \( P_{F_D^\downarrow} \) for some diagram \( D \) (\cref{Lem: diagram alpha rho}).

In \cref{The minimal and maximal elements of the equivalence class}, we characterize \(\min C\) and \(\max C\) by identifying posets \(P, Q \in \poset{n}\) such that \(\min C = \SGL{P}\) and \(\max C = \SGL{Q}\). 
In particular, if \(C\) is the equivalence class of \([w_0(S), \rho]_L\) (respectively \([\sigma, w_1(S)]_L\)), then \(\min C = [w_0(S), \rho]_L\) (respectively \(\max C = [\sigma, w_1(S)]_L\)) (\cref{Thm: Descriptions of minC_D and maxC_D}).
For general lower or upper descent intervals, explicitly describing the poset structure of \( (C, \preceq) \) is challenging. 
However, under a suitable condition, we establish that \( (C, \preceq) \cong (\ST(D), \ble) \) as posets, where \( D \) is a diagram related to the given interval and \( (\ST(D), \ble) \) denotes the poset of standard tableaux on \( D \).
In this case, the equivalence class \( C \) of a lower descent interval has an upper descent interval as \( \max C \), and the equivalence class \( C \) of an upper descent interval has a lower descent interval as \( \min C \) (\cref{Thm: in case of Fsw = Fdown}). 
Using this property, we determine the injective hull of \( \sfB([w_0(S), \rho]_L) \) and the projective cover of \( \sfB([\sigma, w_1(S)]_L) \) (\cref{how to find injective hull}).

In \cref{Lower and Upper Descent Weak Bruhat Intervals from Tableau-Cyclic}, we examine lower and upper descent intervals arising from a projective indecomposable \( H_n(0) \)-module. 
In the context of the representation theory of \(0\)-Hecke algebras, lower and upper descent intervals arise from appropriately selected quotient modules and submodules of the projective modules with the basis \([w_0(S), w_1(S')]_R\), where \(S \subseteq S' \subseteq [n-1]\). 

In \cref{Lower descent intervals from quotient modules}, we consider the $H_n(0)$-modules arising from the sequence of surjective \( H_n(0) \)-module homomorphisms
\begin{equation*}
\begin{tikzcd}
\bfP_\alpha \arrow[r, tail, twoheadrightarrow] &
\mDIF{\alpha} \arrow[r, tail, twoheadrightarrow] & \mESF{\alpha} \arrow[r, tail, twoheadrightarrow] & \mYQS{\alpha,\calC} \arrow[r, tail, twoheadrightarrow] & \bfF_\alpha
\end{tikzcd}
\end{equation*}
given in \cite[Corollary 4.6]{22CKNO1}.
It was shown in \cite{22JKLO} that all the modules in this series are, up to isomorphism, weak Bruhat interval modules and
that the corresponding intervals are lower descent intervals. 
We show that all of them satisfy the condition in \cref{Thm: in case of Fsw = Fdown} and provide explicit descriptions of the equivalence classes \( C \), along with the poset structures of \( (C, \preceq) \) 
(\cref{Prop: The case of mDIF}, \cref{Prop: Case where Y = X} and \cref{Prop: Case where Y = hS}).
Combining these with \cref{how to find injective hull}, 
we determine injective hulls of $\mDIF{\alpha}$, $\mESF{\alpha}$, and $\mYQS{\alpha,\calC}$ in a uniform manner
(\cref{Coro: injective hull of V}, \cref{Coro: injective hull of X}, and \cref{Coro: injective hull of hS}).
It should be noted that our injective hull of \( \mESF{\alpha} \) is a projective indecomposable module, and that an injective hull of \( \mDIF{\alpha} \) was already constructed in \cite[Theorem 4.11]{22CKNO2} in a different manner. 

In \cref{Upper descent intervals from submodules I}, by applying the anti-involution twist of \( \uptheta \circ \upchi \) to the above sequence, we obtain the sequence of injective $H_n(0)$-module homomorphisms 
\[
\begin{tikzcd}
\bfF_{\alpha^\rmc} \arrow{r}{} &
\mRQS{\alpha,\calC} \arrow{r}{} & \mRESF{\alpha} \arrow{r}{} & \mRDIF{\alpha} \arrow{r}{} & \bfP_{\alpha^\rmt}
\end{tikzcd}
\]
(see \cref{involutive images of quotient modules}).
We provide explicit descriptions of the equivalence classes, the poset structures, and the projective covers for these modules (\cref{twisted modules} and \cref{covers of twisted modules}).

In \cref{Upper descent intervals from submodules II}, we discuss the submodules $\HnZmodR_{\alpha}$ of projective 
\( H_n(0) \)-modules that arise from the representation theory of \(0\)-Hecke--Clifford algebras. 
Here, $\alpha$ ranges over the set of peak compositions of $n$. 
These modules are spanned by standard peak immaculate tableaux and are related to the quasisymmetric Schur \( Q \)-functions.
Given a peak composition \(\alpha\) of \(n\), we provide an explicit description of the equivalence class, the poset structure, as well as a projective cover and injective hull for \(\HnZmodR_{\alpha}\) 
(\cref{Thm: Class of Qr} and \cref{Thm: Proj. and Inj. of Qr}). As an importance consequence, we derive that \(\HnZmodR_{\alpha}\) is indecomposable 
(\cref{Coro: indecomposability of Qr}).

\section{Preliminaries}

Throughout this paper, let \( n \) be a positive integer. 
Define \( [n] \) as \( \{1, 2, \ldots, n\} \), and set \( [0] := \emptyset \).

\subsection{Compositions}
A \emph{composition} $\alpha$ of a nonnegative integer $n$, denoted by $\alpha \models n$, is a finite ordered list of positive integers $(\alpha_1, \alpha_2, \ldots, \alpha_l)$ satisfying $\sum_{i=1}^l \alpha_i = n$.
We call $l$ the \emph{length} of $\alpha$ and denote it by $\ell(\alpha)$.
For convenience, we define the empty composition $\emptyset$ to be the unique composition of size and length $0$.
If $\alpha_1 \ge \alpha_2 \ge \cdots \ge \alpha_l$, then we say that $\alpha$ is a \emph{partition} of $n$.

Given $\alpha = (\alpha_1, \alpha_2, \ldots,\alpha_l) \models n$ and $I = \{i_1 < i_2 < \cdots < i_p\} \subset [n-1]$,
let
\begin{align*}
&\set(\alpha) := \{\alpha_1,\alpha_1+\alpha_2,\ldots, \alpha_1 + \alpha_2 + \cdots + \alpha_{l-1}\}, \\
&\comp(I) := (i_1,i_2 - i_1,\ldots,n-i_p).
\end{align*}
The set of compositions of $n$ is in bijection with the set of subsets of $[n-1]$ under the correspondence $\alpha \mapsto \set(\alpha)$ (or $I \mapsto \comp(I)$). 
Define 
\begin{itemize}
\item $\alpha^\rmr$ by the reverse composition $(\alpha_l, \ldots, \alpha_1)$, 

\item $\alpha^\rmc$ by the complement composition satisfying $\set(\alpha^\rmc) = [n-1] \setminus \set(\alpha)$, 

\item $\alpha^\rmt$ by the transpose composition $ (\alpha^\rmr)^\rmc$, and 

\item $\tal$ by the partition obtained by sorting the parts of $\alpha$ in the weakly decreasing order.
\end{itemize}
And, for $I \subset [n-1]$, define 
$I^\rmr:=\set(\comp(I)^\rmr)$ and $I^\rmt:=\set(\comp(I)^\rmt)$.

Let $\alpha = (\alpha_1,\ldots, \alpha_l) \models n$.
We define the \emph{composition diagram} $\tcd(\alpha)$ of $\alpha$ as a left-justified array of $n$ boxes where the $i$th row from the bottom has $\alpha_i$ boxes for $1 \le i \le l$.
We also define the \emph{ribbon diagram} $\trd(\alpha)$ of $\alpha$ as the connected skew diagram without $2 \times 2$ boxes, such that the $i$th column from the left has $\alpha_i$ boxes.
For instance, if $\alpha = (1,3,2)$, then
\[
\tcd(\alpha) = 
\begin{array}{c}
\begin{ytableau}
~ & ~ \\
~ & ~ & ~ \\
~ 
\end{ytableau}
\end{array}
\quad \text{and} \quad
\trd(\alpha) = 
\begin{array}{c}
\begin{ytableau}
~ & ~ \\ 
\none & ~ \\ 
\none & ~ & ~ \\ 
\none & \none & ~ \\
\end{ytableau}
\end{array}.
\]

\subsection{The weak Bruhat orders on the symmetric group}\label{Sec: weak Bruaht order}

The symmetric group $\SG_n$ is generated by simple transpositions $s_i := (i,i+1)$ with $1 \le i \le n-1$. 
An expression of $\sigma \in \SG_n$ that uses the minimal number of simple transpositions is called a \emph{reduced expression for $\sigma$}.
This minimal number is denoted by $\ell(\sigma)$ and called the \emph{length} of $\sigma$. 

The {\em left descent set} and {\em right descent set} of a permutation $\sigma$ are defined by 
\begin{align*}
\DesLR{L}{\sigma} &:= \{i \in [n-1] \mid \ell(s_i \sigma) < \ell(\sigma)\}  \quad \text{and}\\ 
\DesLR{R}{\sigma} &:= \{i \in [n-1] \mid \ell(\sigma s_i) < \ell(\sigma)\},
\end{align*}
respectively.
The \emph{left weak Bruhat order} $\preceq_L$ and 
\emph{right weak Bruhat order} $\preceq_R$ on $\SG_n$
are defined to be the partial order on $\SG_n$ whose covering relation $\preceq_L^c$ and $\preceq_R^c$ are given as follows: 
\begin{align*}
& \sigma \preceq_L^c s_i \sigma \ \text{ if and only if } \ i \notin \DesLR{L}{\sigma} \text{ and }\\
& \sigma \preceq_R^c \sigma s_i \ \text{ if and only if }\  i \notin \DesLR{R}{\sigma},
\end{align*}
respectively.
Given $\sigma, \rho \in \SG_n$,
if $\sigma \preceq_L \rho$, then 
the \emph{$($left$)$ weak Bruhat interval} from $\sigma$ to $\rho$ is defined by
\[
[\sigma,\rho]_L := \{\gamma \in \SG_n \mid \sigma \preceq_L \gamma \preceq_L \rho \},
\]
and if $\sigma \preceq_R \rho$, then 
the \emph{right weak Bruhat interval}  from $\sigma$ to $\rho$ is defined by
\[
[\sigma,\rho]_R := \{\gamma \in \SG_n \mid \sigma \preceq_R \gamma \preceq_R \rho \}.
\]

For $S \subseteq [n-1]$, let $\SG_{S}$ be the parabolic subgroup of $\SG_n$ generated by $\{s_i \mid i\in S\}$ and $w_0(S)$ the longest element in $\SG_{S}$. When $S=[n-1]$, we simply write $w_0$ for $w_0(S)$.
An element $w \in \SG_n$ can be written uniquely as $w=zu$, where $z\in \SG^S$ and $u \in \SG_S$, with the property that  $\ell(w)=\ell(z)+\ell(u)$.
Here $\SG^S:=\{z\in \SG_n \mid \DesLR{R}{z} \subseteq S^{\rm c}\}$ is the set of minimal length representatives for left $\SG_S$-cosets, where $S^{\rm c}=[n-1]\setminus S$. 

\begin{theorem} {\rm (\cite[Theorem 6.2]{88BW})} \label{Bjorner and Wachs}
Given $S \subseteq T \subseteq [n-1]$,
the set $\{w\in \SG_n \mid S \subseteq \DesLR{R}{w} \subseteq T\}$ is exactly the weak Bruhat interval $[w_0(S), w_1(T)]_L$,
where $w_0(S)$ is the longest element in $\SG_{S}$ and $w_1(T)$ is the longest element in $\SG^{{T}^{\rm c}}$.
\end{theorem}

Since $w_1(T)=w_0 w_0({T}^{\rm c})$, \cref{Bjorner and Wachs} can be rewritten as 
\begin{equation}\label{variation of descent class}
\{w\in \SG_n \mid S \subseteq \DesLR{L}{w} \subseteq T\}= [w_0(S), w_0({T}^{\rm c})w_0]_R.
\end{equation}

\subsection{Modules of the \texorpdfstring{$0$}{0}-Hecke algebras from intervals and posets}
\label{Subsec: modules of 0Hecke alg}

The $0$-Hecke algebra $H_n(0)$ is the associative $\C$-algebra with $1$ generated by the elements $\opi_1,\opi_2,\ldots,\opi_{n-1}$ subject to the following relations:
\begin{align}
\begin{aligned}\label{Rel: 0-Hecke algebra}
\opi_i^2 &= -\opi_i \quad \text{for $1\le i \le n-1$},\\
\opi_i \opi_{i+1} \opi_i &= \opi_{i+1} \opi_i \opi_{i+1}  \quad \text{for $1\le i \le n-2$},\\
\opi_i \opi_j &= \opi_j \opi_i \quad \text{if $|i-j| \ge 2$}.
\end{aligned}
\end{align}
Another set of generators consists of $\pi_i:= \opi_i + 1$ for $i=1,2,\ldots,n-1$ with the same relations as above except that $\pi_i^2 = \pi_i$.

For any reduced expression $s_{i_1} s_{i_2} \cdots s_{i_p}$ for $\sigma \in \SG_n$, let $\opi_{\sigma} := \opi_{i_1} \opi_{i_2} \cdots \opi_{i_p}$ and $\pi_{\sigma} := \pi_{i_1} \pi_{i_2 } \cdots \pi_{i_p}$.
It is well known that these elements are independent of the choices of reduced expressions, and both $\{\opi_\sigma \mid \sigma \in \SG_n\}$ and $\{\pi_\sigma \mid \sigma \in \SG_n\}$ are $\mathbb C$-bases for $H_n(0)$.

According to \cite{79Norton}, there are $2^{n-1}$ pairwise inequivalent irreducible $H_n(0)$-modules and $2^{n-1}$ pairwise inequivalent projective indecomposable $H_n(0)$-modules, which are naturally indexed by compositions of $n$.
For a composition $\alpha$ of $n$, let $\bfF_{\alpha}$ denote the $1$-dimensional $\C$-vector space corresponding to the composition $\alpha$ of $n$, spanned by a vector $v_{\alpha}$. 
For each $1\le i \le n-1$, define an action of the generator $\pi_i$ of $H_n(0)$ as follows:
\[
\pi_i(v_\alpha) = \begin{cases}
0 & i \in \set(\alpha),\\
v_\alpha & i \notin \set(\alpha).
\end{cases}
\]
This module is the irreducible $1$-dimensional $H_n(0)$-module corresponding to $\alpha$.
And, the projective indecomposable $H_n(0)$-module corresponding to $\alpha$
is given by the submodule 
$\bfP_\alpha = H_n(0) \pi_{w_0(\set(\alpha)^\rmc)} \opi_{w_0(\set(\alpha))}$ of the regular representation of $H_n(0)$.

One can construct modules of the $0$-Hecke algebra using various combinatorial objects. 
This paper focuses on modules arising from weak Bruhat intervals and posets in $\poset{n}$. 
First, we review the weak Bruhat interval $H_n(0)$-modules introduced by Jung--Kim--Lee--Oh (\cite{22JKLO}).

\begin{definition}{\rm (\cite[Definition 1]{22JKLO})} \label{left WBI}
Let $I \in \rmInt(n)$. 
\begin{enumerate}[label = {\rm (\alph*)}]
\item The \emph{weak Bruhat interval module associated with $I$}, denoted by $\sfB(I)$, is the left $H_n(0)$-module with  $\C I$ as the underlying space and with the $H_n(0)$-action defined by
\begin{equation*}\label{Hecke algebra action: left and pi}
\pi_i \cdot \gamma := 
\begin{cases}
\gamma & \text{if $i \in \DesLR{L}{\gamma}$}, \\
0 & \text{if $i \notin \DesLR{L}{\gamma}$ and $s_i\gamma \notin I$,} \\
s_i \gamma & \text{if $i \notin \DesLR{L}{\gamma}$ and $s_i\gamma \in I$}
\end{cases} 
\end{equation*}
for $1 \leq i \leq n-1$ and $\gamma \in I$.

\item The \emph{negative weak Bruhat interval module associated with $I$}, denote by $\osfB(I)$, is the left $H_n(0)$-module with  $\C I$ as the underlying space and with the $H_n(0)$-action defined by
\begin{equation*}\label{Hecke algebra action: left and pi bar}
\opi_i \star \gamma := 
\begin{cases}
-\gamma & \text{if $i \in \DesLR{L}{\gamma}$}, \\
0 & \text{if $i \notin \DesLR{L}{\gamma}$ and $s_i\gamma \notin I$,} \\
s_i \gamma & \text{if $i \notin \DesLR{L}{\gamma}$ and $s_i\gamma \in I$}
\end{cases} 
\end{equation*}
for $1 \leq i \leq n-1$ and $\gamma \in I$.
\end{enumerate}
\end{definition}
Indeed $\osfB(I)$ is the $\uptheta$-twist of $\sfB(I)$. For the definition of the $\uptheta$-twist, see \cref{Upper descent intervals from submodules I}.

Next, we review the left $H_n(0)$-modules arising from posets.
Recall that $\poset{n}$ is the set of posets with ground set $[n]$.
Given any poset $P \in \poset{n}$, 
we define 
\begin{equation}\label{definition of left linear extensions}
\SGL{P} := \{\sigma \in \SG_n \mid \text{$\sigma(i) \leq \sigma(j)$ for all $i, j\in [n]$ with $i \preceq_P j$}\}.
\end{equation}

\begin{definition} {\rm (\cite[Definition 2.8]{23KLO})}
Let $P \in \poset{n}$. Define the $H_n(0)$-module $\ourMP{P}$ to be the underlying space $\C\SGL{P}$ and with the $H_n(0)$-action:
\begin{align*}
\pi_i \cdot \gamma:= 
\begin{cases}
\gamma & \text{if $i \in \DesLR{L}{\gamma}$},\\
0 & \text{if $i \notin \DesLR{L}{\gamma}$ and $s_i \gamma \notin \SGL{P}$},\\
\gamma s_i & \text{if $i \notin \DesLR{L}{\gamma}$  and $s_i \gamma \in \SGL{P}$}
\end{cases}
\end{align*}
for $i \in [n-1]$ and $\gamma \in \SGL{P}$.
\end{definition}

It should be noted that in the case where $\SGL{P}$ is a left weak Bruhat interval, $\ourMP{P}$ is identical to $\sfB(\SGL{P})$.

\section{An equivalence relation on \texorpdfstring{$\mathrm{Int}(n)$}{Intn} and its  poset-theoretic characterization}
\label{poset-theoretical characterization}

Let \( \rmInt(n) \) denote the set of nonempty weak Bruhat intervals in \( \SG_n \).  
In this section, we present a poset-theoretic characterization of an equivalence relation on \( \rmInt(n) \) introduced by Kim--Lee--Oh in \cite{23KLO}.  
We begin by recalling the definition of this equivalence relation.

For \( I_1, I_2 \in \rmInt(n) \), a poset isomorphism $f: (I_1, \preceq_L) \to (I_2, \preceq_L)$ is said to be \emph{descent-preserving} if  
\[
\DesLR{L}{\gamma} = \DesLR{L}{f(\gamma)} \quad \text{for all } \gamma \in I_1.
\]
Define an equivalence relation \( \Deq \) on \( \rmInt(n) \) such that \( I_1 \Deq I_2 \) if there exists a descent-preserving poset isomorphism between \( (I_1, \preceq_L) \) and \( (I_2, \preceq_L) \).
This equivalence relation plays an important role in refining the classification of the \( H_n(0) \)-modules \( \sfB(I) \) associated with \( I \in \rmInt(n) \).  
Specifically, if \( I_1 \Deq I_2 \), then $\sfB(I_1) \cong \sfB(I_2)$ as \( H_n(0) \)-modules.

Kim--Lee--Oh conjectured that the converse also holds for all intervals in \( \rmInt(n) \) \cite[Conjecture 7.2]{23KLO}, and verified the conjecture in the case where the intervals arise from regular Schur-labeled posets, that is, for intervals in the subset
\[
\left\{ \Sigma_L(P) \mid \text{$P$ is a regular Schur-labeled poset on } [n] \right\} \subset \rmInt(n),
\]
as established in \cite[Sections 4 and 5]{23KLO}.

More recently, Yang--Yu \cite{24YY} proved that this conjecture holds in full generality, for all weak Bruhat interval modules in arbitrary finite Coxeter types, including type \( A \).
We begin this section by presenting a detailed characterization of the equivalence relation \( \Deq \).

\begin{lemma}
Let \( I_1 , I_2 \in \rmInt(n) \). 
Then the following are equivalent.
\begin{enumerate}[label = {\rm (\alph*)}]
\item \( I_1 \Deq I_2 \).

\item \( \sfB(I_1) \cong \sfB(I_2) \) as \( H_n(0) \)-modules.

\item \( \osfB(I_1) \cong \osfB(I_2) \) as \( H_n(0) \)-modules.
\end{enumerate}
\end{lemma}
\begin{proof}
The equivalence (a) \(\Leftrightarrow\) (b) was established in \cite[Theorem 4.11]{24YY}. Furthermore, the equivalence (b) \(\Leftrightarrow\) (c) follows from \cite[Theorem 4 (2)]{22JKLO}.
\end{proof}

Every weak Bruhat interval corresponds to the set of linear extensions of a regular poset. 
We now recall the definition of regular posets.
\begin{definition}{\rm (\cite[p. 110]{91BW})}\label{def: regular posets}
A poset $P \in \poset{n}$ is said to be \emph{regular} if the following holds:
for all $x,y,z \in [n]$ with
$x \preceq_P z$, 
if $x < y < z$ or $z<y<x$, then $x \preceq_P y$ or $y \preceq_P z$.
\end{definition}
We denote by $\Rposet{n}$ the set of all regular posets in $\poset{n}$.
The following theorem shows how regular posets can be characterized in terms of weak Bruhat intervals.

\begin{theorem}{\rm (\cite[Theorem 6.8]{91BW})}
\label{thm: left interval and regular}
Let $U$ be a nonempty subset of $\SG_n$.
Then, the following conditions are equivalent:
\begin{enumerate}[label = {\rm (\arabic*)}]
\item $U$ is a weak Bruhat interval.

\item $U = \Sigma_L(P)$ for some $P \in \Rposet{n}$.
\end{enumerate}  
\end{theorem}

Consider the map
\[
\upeta: \poset{n} \to \mathscr{P}(\SG_n), \quad P \mapsto \Sigma_L(P),
\]
where $\mathscr{P}(\SG_n)$ is the power set of $\SG_n$.
One can see that $\upeta$ is injective.
Combining this with \cref{thm: left interval and regular}, we obtain a one-to-one correspondence 
\[
\upeta|_{\Rposet{n}}: \Rposet{n} \to \rmInt(n), \quad P \mapsto \Sigma_L(P).
\]
For $I \in \rmInt(n)$, we denote by $P_I$ the regular poset such that $\SGL{P_I} = I$.
Throughout this paper, we will identify any weak Bruhat interval $I\in \SG_n$ with the regular poset $P_I \in \Rposet{n}$.
Based on this identification, we provide a poset-theoretical characterization of the equivalence relation $\Deq$ on $\rmInt(n)$.

\begin{definition}{\rm (\cite[Definition 1.5]{23AvDjMc})}.
An edge-decorated poset is a poset $P$ such that each edge in its Hasse diagram is assigned to be either weak or strict. 
\end{definition}

Recall that each poset \( P \in \mathsf{P}_n \) can be naturally identified with the labeled poset \( (P, \omega) \), where the labeling \( \omega: P \to [n] \) is defined by \( \omega(i) = i \) for all \( i \in P \). 
Consequently, \( P \) inherits the structure of an edge-decorated poset.
Following the convention for drawing Hasse diagrams of posets, we draw a bold edge (referred to as a {\em strict} edge) between \(x\) and \(y\) when \(x \preceq_P y\) and \(\omega(x) > \omega(y)\), and a plain edge (referred to as a {\em weak} edge) when \(x \preceq_P y\) and \(\omega(x) < \omega(y)\) in the Hasse diagram of \(P\).
Here $\preceq_{P}$ is used to denote the partial order of $P$
and $\le$ the usual order on $[n]$.

We say that two labeled posets \((P, \omega)\) and \((Q, \tau)\) are isomorphic, denoted by \((P, \omega) \Leq (Q, \tau)\), if there exists a poset isomorphism from \(P\) to \(Q\) that maps strict edges and weak edges in \(P\) to strict edges and weak edges in \(Q\), respectively.
Then it follows from the definition of $(P,\omega)$-partition generating function that 
\[
(P, \omega) \Leq (Q,\tau) \Rightarrow (P,\omega) \Ceq (Q,\tau)
\]
(see \cite[Lemma 3.6]{14MW}).

The symmetric group \(\SG_n\) acts on \(\poset{n}\) by composing the labeling with permutations. Specifically, for \(\sigma \in \SG_n\) and \((P, \omega) \in \poset{n}\), the action is defined as  
\[
\sigma \cdot (P, \omega) = (P, \sigma \circ \omega).
\]
\begin{lemma}\label{Lem: regular poset and indirectly comparable pair}
Let $P \in \Rposet{n}$ and let $1 \leq i \leq n-1$.
\begin{enumerate}[label = {\rm (\alph*)}]
\item 
If $(i,i+1)$ is a comparable pair that is not in the covering relation in $P$, then $s_i \cdot P$ is equal to $P$ as an edge-decorated poset and remains regular.

\item 
If $(i,i+1)$ is in the covering relation in $P$, then $s_i \cdot P$ is not equal to $P$ as an edge-decorated poset, but it remains regular.

\item 
If \((i, i+1)\) is an incomparable pair in \(P\), then \(s_i \cdot P\) equals \(P\) as an edge-decorated poset. Moreover, either \(\SGL{P} = \SGL{s_i \cdot P}\) or \(s_i \cdot P\) is not regular.
\end{enumerate}
\end{lemma}
\begin{proof}
The assertions are straightforward when \( n = 1 \) or \( n = 2 \).  
Thus, we assume \( n \geq 3 \).

(a) Since \(i\) and \(i+1\) are comparable but not in the covering relation in \(P\), it is evident that \(s_i \cdot P\) is identical to \(P\) as edge-decorated posets.

To show that \(s_i \cdot P\) remains regular, choose any triple \((x, y, z)\) in \(s_i \cdot P\) with \(x \preceq_{s_i \cdot P} z\). 
If none of \(x, y, z\) are equal to \(i\) or \(i+1\), the triple \((x, y, z)\) remains unchanged between \(P\) and \(s_i \cdot P\), so the regularity is preserved. 
If only one of \(x, y, z\) equals \(i\) or \(i+1\), the regularity of \(P\) ensures that \(s_i \cdot P\) remains regular as well.
Otherwise, the following cases are possible:
\begin{itemize}
\item \(x = i < y = i+1 < z\),
\item \(x < y = i < z = i+1\),
\item \(z = i < y = i+1 < x\), or
\item \(z < y = i < x = i+1\).
\end{itemize}
In all these cases, the assumption that \((i, i+1)\) is a comparable pair in \(P\) guarantees that either \(x \preceq_{s_i \cdot P} y\) or \(y \preceq_{s_i \cdot P} z\). 
Therefore, \(s_i \cdot P\) satisfies the regularity conditions.

(b) If \(i\preceq_{P} i+1\) is in the covering relation in \(P\), then \(i+1 \preceq_{s_i \cdot P} i\) is in the covering relation in \(s_i \cdot P\), making \(s_i \cdot P\) distinct from \(P\) as an edge-decorated poset. 

The regularity of \(s_i \cdot P\) follows directly from the argument in (a), as the reversal of a covering relation does not disrupt regularity.

(c) Let \((i, i+1)\) be an incomparable pair in \(P\). 
Since the relative orders in the Hasse diagrams of the connected components containing \(i\) and \(i+1\) in \(P\) remain unchanged, it follows that \(s_i \cdot P\) is identical to \(P\) as an edge-decorated poset.

We now prove the second assertion. First, consider the case where every element \(x \in P\) other than \(i, i+1\) is either incomparable to both \(i\) and \(i+1\), or comparable to both \(i\) and \(i+1\). In this case, the set of linear extensions remains unchanged, and thus \(\SGL{P} = \SGL{s_i \cdot P}\).
Next, suppose there exists an element \(x \neq i, i+1\) in \(P\) that is comparable to either \(i\) or \(i+1\), but not both. We focus on the case where \(x \prec_P i\) or \(i \prec_P x\), noting that the argument for the case where \(x \prec_P i+1\) or \(i+1 \prec_P x\) follows analogously.
We have the following four subcases:
\begin{itemize}
\item
\(x \prec_P i\) and $x < i$: \ 
It holds that \(x \prec_{s_i \cdot P} i+1\) in \(s_i \cdot P\) and 
$x<i<i+1$. However, \(x \npreceq_{s_i \cdot P} i\) and \(i \npreceq_{s_i \cdot P} i+1\) in \(s_i \cdot P\). Thus \(s_i \cdot P\) is not regular.

\item 
\(x \prec_P i\) and $i+1<x$: \ 
It holds that $i<i+1<x$. 
However, \(x \npreceq_P i+1\) and \(i+1 \npreceq_P i\) in \(\cdot P\). 
This contradicts the assumption that $P$ is regular. 

\item \(i \prec_P x\) and $x < i$: \ 
It holds that \(i+1 \prec_{s_i \cdot P} x\) in \(s_i \cdot P\) and $x<i<i+1$. 
However, \(i+1 \npreceq_{s_i \cdot P} i\) and \(i \npreceq_{s_i \cdot P} x\) in \(s_i \cdot P\). Thus \(s_i \cdot P\) is not regular.

\item \(i \prec_P x\) and $i+1<x$: \ 
It holds that $i<i+1<x$. 
However, \(i \npreceq_P i+1\) and \(i+1 \npreceq_P x\) in \(\cdot P\). 
This contradicts the assumption that $P$ is regular. 
\end{itemize}
\end{proof}

\begin{remark}
From the perspective of equivalence relations, we note that in \cref{Lem: regular poset and indirectly comparable pair} (a), $s_i \cdot P \Meq P$, and thus $s_i \cdot P \Ceq P$; in (b), $s_i \cdot P \overset{K}{\not \simeq} P$, and thus $s_i \cdot P \overset{M}{\not \simeq} P$; and in (c), $s_i \cdot P \Ceq P$.
\end{remark}

\begin{example}
Consider the regular poset 
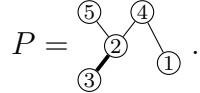
\[
P=\begin{tikzpicture}[baseline=3.5mm]
\def \hp {0.35}
\def \vp {0.45}
\def \ccc {1mm}
\node[shape=circle,draw,minimum size=\ccc*3, inner sep=0pt] at (0, 0) (A5) {\tiny $3$};
\node[shape=circle,draw,minimum size=\ccc*3, inner sep=0pt] at (\hp*1, \vp) (A1) {\tiny $2$};
\node[shape=circle,draw,minimum size=\ccc*3, inner sep=0pt] at (0*\hp, \vp*2) (A3) {\tiny $5$};
\node[shape=circle,draw,minimum size=\ccc*3, inner sep=0pt] at (2*\hp, \vp*2) (A4) {\tiny $4$};
\node[shape=circle,draw,minimum size=\ccc*3, inner sep=0pt] at (3*\hp, \vp*0.5) (A2) {\tiny $1$};

\draw[line width = \lw] (A5) -- (A1);
\draw (A1) -- (A3);
\draw (A1) -- (A4) -- (A2);
\end{tikzpicture} \ .
\]
Note that $(3,4)$ is a comparable non-covering pair,  $(2,3)$ is a covering pair, and $(1,2),(4,5)$ are incomparable pairs in $P$.
Now, consider the posets obtained by applying simple transpositions to $P$:
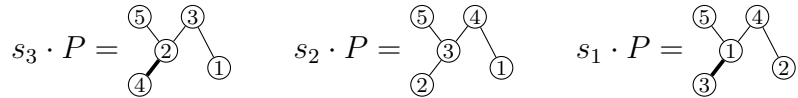
\[
s_3 \cdot P = \begin{tikzpicture}[baseline=3.5mm]
\def \hp {0.35}
\def \vp {0.45}
\def \ccc {1mm}
\node[shape=circle,draw,minimum size=\ccc*3, inner sep=0pt] at (0, 0) (A5) {\tiny $4$};
\node[shape=circle,draw,minimum size=\ccc*3, inner sep=0pt] at (\hp*1, \vp) (A1) {\tiny $2$};
\node[shape=circle,draw,minimum size=\ccc*3, inner sep=0pt] at (0*\hp, \vp*2) (A3) {\tiny $5$};
\node[shape=circle,draw,minimum size=\ccc*3, inner sep=0pt] at (2*\hp, \vp*2) (A4) {\tiny $3$};
\node[shape=circle,draw,minimum size=\ccc*3, inner sep=0pt] at (3*\hp, \vp*0.5) (A2) {\tiny $1$};

\draw[line width = \lw] (A5) -- (A1);
\draw (A1) -- (A3);
\draw (A1) -- (A4) -- (A2);
\end{tikzpicture} 
\qquad 
s_2 \cdot P = \begin{tikzpicture}[baseline=3.5mm]
\def \hp {0.35}
\def \vp {0.45}
\def \ccc {1mm}
\node[shape=circle,draw,minimum size=\ccc*3, inner sep=0pt] at (0, 0) (A5) {\tiny $2$};
\node[shape=circle,draw,minimum size=\ccc*3, inner sep=0pt] at (\hp*1, \vp) (A1) {\tiny $3$};
\node[shape=circle,draw,minimum size=\ccc*3, inner sep=0pt] at (0*\hp, \vp*2) (A3) {\tiny $5$};
\node[shape=circle,draw,minimum size=\ccc*3, inner sep=0pt] at (2*\hp, \vp*2) (A4) {\tiny $4$};
\node[shape=circle,draw,minimum size=\ccc*3, inner sep=0pt] at (3*\hp, \vp*0.5) (A2) {\tiny $1$};

\draw (A5) -- (A1) -- (A3);
\draw (A1) -- (A4) -- (A2);
\end{tikzpicture}
\qquad
s_1 \cdot P = \begin{tikzpicture}[baseline=3.5mm]
\def \hp {0.35}
\def \vp {0.45}
\def \ccc {1mm}
\node[shape=circle,draw,minimum size=\ccc*3, inner sep=0pt] at (0, 0) (A5) {\tiny $3$};
\node[shape=circle,draw,minimum size=\ccc*3, inner sep=0pt] at (\hp*1, \vp) (A1) {\tiny $1$};
\node[shape=circle,draw,minimum size=\ccc*3, inner sep=0pt] at (0*\hp, \vp*2) (A3) {\tiny $5$};
\node[shape=circle,draw,minimum size=\ccc*3, inner sep=0pt] at (2*\hp, \vp*2) (A4) {\tiny $4$};
\node[shape=circle,draw,minimum size=\ccc*3, inner sep=0pt] at (3*\hp, \vp*0.5) (A2) {\tiny $2$};

\draw[line width = \lw] (A5) -- (A1);
\draw (A1) -- (A3);
\draw (A1) -- (A4) -- (A2);
\end{tikzpicture}
\]
We observe that the posets \( P \), \( s_3 \cdot P \), and \( s_1 \cdot P \) are identical as edge-decorated posets. 
However, while \( s_3 \cdot P \) is regular, \( s_1 \cdot P \) is not. 
Additionally, \( K_P \neq K_{s_2 \cdot P} \), since \( 12345 \notin \SGL{K_P} \).
\end{example}

We now present a poset-theoretic characterization of \(\Deq\). For the proof, we introduce the notation  
\[
\calR_{\rm st}(P) := \{(x, y) \in [n]^2 \mid x \preceq_{P} y \text{ and } x \neq y\} \quad \text{for  \(P \in \poset{n}\)}.
\]

\begin{theorem} 
\label{characterization of the equivalence relation}
Let $I,J \in \rmInt(n)$. Then 
$I \Deq J$ if and only if $P_J= s_{i_r}\cdots s_{i_2}s_{i_1} \cdot P_I$ for some nonnegative integer $r$,
where $(i_k,i_k+1)$ is a comparable pair that is not in the covering relation in $s_{i_{k-1}}\cdots s_{i_2}s_{i_1} \cdot P_I$ for $1\le k \le r$.
\end{theorem}
\begin{proof}
The assertion is straightforward when \( n = 1 \) or \( n = 2 \).  
Thus, we assume \( n \geq 3 \).
In the trivial case where $I = J$, the assertion is clear, as for any weak Bruhat interval $I$ in $\SG_n$, there exists a unique poset $P_I \in \poset{n}$ such that $\SGL{P_I} = I$. 
Now, assume $I \Deq J$ with $I \neq J$. 
By \cite[Theorem 4.6]{23KLO}, for \(1 \leq i \leq n-1\), the assertion is established by verifying the one-step equivalence: \(I \Deq I s_i\) if and only if \((i, i+1)\) is a comparable pair that is not in the covering relation in \(P_I\); in this case \(P_{I s_i} = s_i \cdot P_I\).

To begin with, we show that if \((i, i+1)\) is a comparable pair in \(P_I\) that does not belong to the covering relation, then \(P_{I s_i} = s_i \cdot P_I\) always holds.
Suppose that $\SGL{P_I} = [\sigma,\rho]_L$.
Since $P_I$ is regular and $(i,i+1)$ is a comparable pair that is not in the covering relation in $P_I$, it follows from \cref{Lem: regular poset and indirectly comparable pair} that $s_i \cdot P_I$ is regular.
So $\SGL{s_i \cdot P_I}$ is a weak Bruhat interval, we put it by $[\sigma',\rho']_L$.
We claim that $\sigma' = \sigma s_i$ and $\rho' = \rho s_i$.

For a given regular poset \( P \), let \( \SGL{P} = [\delta, \eta]_L \). 
It is straightforward to verify that \( \delta \) and \( \eta \) can be determined from \( P \) as follows: for \( 1 \leq k \leq n \),  
\begin{equation}\label{how to determine min and max}
\begin{aligned}  
\delta(k) &= \left| \{x \mid x \preceq_{P} k\} \cup \{x \mid \text{$x$ is incomparable to $k$ in $P$ and \( x < k \)}\} \right|, \\  
\eta(k) &= \left| \{x \mid x \preceq_{P} k\} \cup \{x \mid \text{$x$ is incomparable to $k$  in \( P \) and \( x > k \)}\} \right|.  
\end{aligned}
\end{equation}
Applying this property to the posets $P_I$ and $s_i \cdot P_I$, we observe that $\sigma'(i) = \sigma(i+1)$, $\rho'(i) = \rho(i+1)$, and
all other entries remain unchanged. 
Therefore, we have 
\begin{equation}\label{Eq: si P_I equivalent P_I s_i}
\SGL{s_i \cdot P_I} =I s_i, \text{ equivalently }s_i \cdot P_I=P_{Is_i}.
\end{equation}

Now, we prove the ``if" direction.
We have two cases.
\begin{itemize}
\item \(i \prec_{P_I} i+1\): \ 
From \eqref{definition of left linear extensions} it holds that $\gamma(i) < \gamma(i+1)$ for each $\gamma \in \SGL{P}$.
Since $(i,i+1)$ is not in the covering relation in $P_I$, $\gamma(i)+1 < \gamma(i+1)$.
Then one can easily see that 
\[
\DesLR{L}{\gamma} = \DesLR{L}{\gamma s_i} \quad \text{for each $\gamma \in \SGL{P_I}$}.
\]
Since $\SGL{P_I} = I$ and $\SGL{s_i \cdot P_I} = I s_i$, this implies a descent-preserving bijection $I \to I s_i, \ \gamma \mapsto \gamma s_i$.
Hence, \(I \Deq I s_i\). 

\item \(i+1 \prec_{P_I} i\): \ 
This case follows from a symmetric argument.
\end{itemize}

Next, we prove the ``only if" direction.
Suppose $I \Deq I s_i$.
Then there exists a descent-preserving bijection \(f:I \to I s_i,\ \gamma \mapsto \gamma s_i\). 
Let \(\SGL{P_I} = [\sigma, \rho]_L\). 
We have two cases.
\begin{itemize}
\item \(\rho(i) < \rho(i+1)\): \ 
Since \(\sigma \preceq_L \rho\), it follows that \(\sigma(i) < \sigma(i+1)\). 
By \cite[Theorem 6.8]{91BW}, the strict relations in \(P_I\) are given by
\[
\calR_{\rm st}(P_I) = \{(x, y) \mid \sigma(x) < \sigma(y) \text{ and } \rho(x) < \rho(y)\}.
\]
This implies \(i \prec_{P_I} i+1\). 
If $i+1$ covers $i$, then it follows from \eqref{definition of left linear extensions} that \(\sigma(i+1) = \sigma(i) + 1\).
It contradicts the descent-preserving map $f$. 
Hence, \((i, i+1)\) is a comparable pair but not in the covering relation in $P_I$.

\item \(\rho(i) > \rho(i+1)\): \ 
This case follows from a symmetric argument.
\end{itemize}

This completes the proof.
\end{proof}

\section{Equivalence classes of lower and upper descent weak Bruhat intervals}
\label{Sec: The equivalence classes of weak Bruhat Intervals}

In this section, we investigate the equivalence classes of lower and upper descent intervals, specifically weak Bruhat intervals of the forms 
\begin{align*}
[w_0(S),\rho]_L \quad \text{or} \quad [\sigma,w_1(S)]_L, \quad (S \subseteq [n-1]).
\end{align*}
We first demonstrate that these intervals can be represented as \(\Sigma_L(P)\) for some special posets \(P \in \poset{n}\), which we refer to as \emph{canonical diagram posets}. 
Next, given an equivalence class \(C\), we provide the posets \(P\) and \(Q\) in \(\poset{n}\) such that 
\[
\min C = \SGL{P} \quad \text{and} \quad \max C = \SGL{Q}.
\]

\subsection{The canonical diagram posets}
\label{Sec41: The canonical diagram posets}

We recall the notion of canonical diagram posets as described in \cite[Section 5.1]{24CKO}. Then, we examine their images under the typical involutions defined on \(\poset{n}\).

In this subsection, we consider an \( n \)-element subset of \( \mathbb{N}^2 \), which we will treat as a diagram composed of \( n \) boxes located in the first quadrant. 
We identify each point \( (i,j) \) with the empty rectangle whose vertices are at 
\( (i - 1, j - 1) \), \( (i , j - 1) \), \( (i - 1, j ) \), and \( (i , j) \) (see \cref{Ex: S = 25 and its diagram}).
Let \( \mathfrak{D}_n \) be the set of \( n \)-element subsets of \( \mathbb{N}^2 \) such that the corresponding diagram has no empty rows or columns within the smallest rectangle that can completely enclose the diagram. 
In particular, we consider the composition diagrams \(\tcd(\alpha)\) and the ribbon diagrams \(\trd(\alpha)\) of size \(n\) as elements of \( \mathfrak{D}_n \) by positioning the lower leftmost box at \((1,1)\).
Unless explicitly stated otherwise, we will assume throughout this section that $D \in \mathfrak{D}_n$.
\vspace*{3mm}

\noindent
\textbf{Convention.}
In this paper, we regard a filling $F$ of $D$ with positive integers as a map 
$$F: D \to \Z_{>0}, \quad (i,j) \mapsto F(i,j),$$
where $F(i,j)$ denotes the entry in the box $(i,j)$ of $F$.
Given a filling \( F \) of \( D \), we primarily utilize four reading words, denoted as \( \readingTBLR(F) \), \( \readingLRTB(F) \), \( \readingLRBT(F) \), and \( \readingBTLR(F) \), which are defined as follows: 
\begin{itemize}
\item $\readingTBLR(F)$ is the word obtained by reading the entries of $F$ from top to bottom in each column, starting with the leftmost column.

\item $\readingLRTB(F)$ is the word obtained by reading the entries of $F$ from left to right across each row, starting with the topmost row.

\item $\readingLRBT(F)$ is the word obtained by reading the entries of $F$ from left to right across each row, starting with the bottommost row.

\item $\readingBTLR(F)$ is the word obtained by reading the entries of $F$ from bottom to top in each column, starting with the leftmost column.
\end{itemize}
A filling \( F \) on \( D \) with positive integers is called \emph{standard} if it contains the entries \( 1, 2, \ldots, n \), each appearing exactly once. Every word containing distinct entries from \( 1 \) to \( n \) is regarded as a permutation in \( \SG_n \).

\begin{definition} \label{def: diagram posets}
Let \( F \) be a standard filling of \( D \). 
\begin{enumerate}[label = {\rm (\alph*)}]
\item 
Let \( P_F \) denote the poset in \( \poset{n} \) with partial order \( \preceq_{P_F} \), defined by
\[
i \preceq_{P_F} j \quad \text{if and only if} \quad x_i \leq x_j \text{ and } y_i \leq y_j,
\]
where \( (x_i, y_i) \) represents the position of \( i \) in the filling \( F \) for each \( 1 \leq i \leq n \).

\item
Define \( F_D^\downarrow \) (respectively \( F_D^\rightarrow \)) as the filling of \( D \) obtained by placing the integers \( 1, 2, \ldots, n \) sequentially, without repetition, down each column (respectively across each row), from top to bottom (respectively left to right), starting from the leftmost column (respectively the uppermost row). We refer to \( P_{F_D^\downarrow} \) and \( P_{F_D^\rightarrow} \) as the \emph{canonical diagram posets} associated with \( D \).
\footnote{Note that in \cite{24CKO}, only the poset \( P_{F_D^\rightarrow} \) is referred to as the \emph{canonical diagram poset associated with \( D \)}.}
\end{enumerate}
\end{definition}

\begin{remark}\label{Rem: reading permutations from F}
In \cref{def: diagram posets}(a),  
suppose \( P_F \) is regular.  
Let \( \SGL{P_F} = [\sigma, \rho]_L \). 
Using \eqref{how to determine min and max}, one can obtain \( \sigma \) and \( \rho \) directly from \( F \). 
To be precise, for \( 1 \leq k \leq n \),  
\begin{align*}
\sigma(k) = & |\{1 \leq x \leq n \mid \text{$x$ is lower-left of $k$ in $F$}\} \\ 
& \cup \{1 \leq x < k \mid \text{$x$ is strictly upper-left or strictly lower-right of $k$ in $F$} \}|, \quad \text{and} \\
\rho(k) = & |\{1 \leq x \leq n \mid \text{$x$ is lower-left of $k$ in $F$}\} \\
& \cup \{k < x \leq n \mid \text{$x$ is strictly upper-left or strictly lower-right of $k$ in $F$} \}|.
\end{align*}
\end{remark}

Since \( P_{F_D^\downarrow} \) and \( P_{F_D^\rightarrow} \) are regular posets, the sets \( \SGL{P_{F_D^\downarrow}} \) and \( \SGL{P_{F_D^\rightarrow}} \) form weak Bruhat intervals. 
Here, we investigate the properties of these intervals.
To begin with, we introduce the necessary lemma and notations.

(i) Given a diagram $D \in \mathcal{D}_n$, let $D^\rmt$ be the diagram obtained by transposing the coordinates, that is, $D^\rmt = \{(j, i) \mid (i, j) \in D\}$. 
Similarly, for a standard filling $F$ on $D$, let $F^\rmt$ be the corresponding standard filling on $D^\rmt$, defined by setting
\[
F^\rmt(i, j) := F(j, i) \quad \text{for }  (j, i) \in D. 
\]

\begin{lemma}\label{Lem: PF = PF'}
Let $F$ be a standard filling on $D$.
Then we have $P_F = P_{F^\rmt}$.  
\end{lemma}
\begin{proof}
To prove the assertion, we start by noting that $(F^\rmt)^\rmt = F$. 
Thus, it suffices to show that for all \( 1 \leq i, j \leq n \), we have that
\[
i \preceq_{P_{F}} j \implies i \preceq_{P_{F^\rmt}} j.
\]

For $1 \leq i \leq n$, let $(x_i,y_i)$ denote the position of $i$ in the filling $F$.
Suppose that $i \preceq_{P_F} j$.
By the definition of $P_F$, this implies that
\begin{equation}\label{Eq: i preceq_PF j}
x_i \leq x_j \quad \text{and} \quad y_i \leq y_j.
\end{equation}
In the filling $F^\rmt$, the positions of $i$ and $j$ are transposed, with $i$ now at $(y_i, x_i)$ and $j$ at $(y_j, x_j)$. 
The condition \eqref{Eq: i preceq_PF j} can be rewritten as $y_i \leq y_j$ and $x_i \leq x_j$.
Therefore, by the definition of the partial order for $P_{F^\rmt}$, it follows that $i \preceq_{P_{F^\rmt}} j$.
\end{proof}

(ii) A filling \( T \) on \( D \) with entries \( 1, 2, \ldots, n \) is called a \emph{standard tableau} on \( D \) if its entries are distinct and \( T(i,j) \le T(k,l) \) whenever \( i \le k \) and \( j \le l \). 
Let $\ST(D)$ denote the set of standard tableaux on $D$.
Define 
\begin{equation} \label{def of TD and TD prime}
T'_D \text{ (respectively $T_D$)} 
\end{equation}
to be the standard tableau on $D$ obtained by sequentially placing the integers $1,2,\ldots,n$ without repetition, along the columns (respectively rows), from bottom to top (respectively left to right), starting with the leftmost column (respectively the bottommost row).
 
(iii) For \( P \in \poset{n} \), define \( \overline{P} \) as the poset in \( \poset{n} \) where the order relation is given by
\[
u \preceq_{\overline{P}} v \quad \Longleftrightarrow \quad n + 1 - u \preceq_{P} n + 1 - v.
\]
The map ${}^- : \poset{n} \rightarrow \poset{n}, \ P \mapsto \overline{P}$, is clearly an involution.
Furthermore, given a standard filling \( F \) on \( D \), we denote by \( \overline{F} \) the standard filling obtained from $F$ by replacing each entry \( i \) with \( n - i + 1 \) for all \( 1 \leq i \leq n \).

(iv) 
Let \( k \) and \( l \) denote the number of rows and columns of \( D \), respectively. 
For \( 1 \le i \le k \) and \( 1 \le j \le l \), let \( r_i \) and \( c_j \) be the number of boxes in the \( i \)th row (from the top) and the \( j \)th column (from the left) of \( D \), respectively. 
Then, define 
\[
\sfr(D) := (r_1, r_2, \ldots, r_k) \quad \text{and} \quad \sfc(D) := (c_1, c_2, \ldots, c_l).
\]

With the prerequisites (i)--(iv), we determine the minimal and maximal elements of the intervals \( \SGL{P_{F_D^\downarrow}} \) and \( \SGL{P_{F_D^\rightarrow}} \), which correspond to specific reading words derived from the standard tableaux. 

\begin{proposition}
\label{Prop: two kinds of P_D and intervals}
Let $D \in \mathfrak{D}_n$. 

\begin{enumerate}[label = {\rm (\alph*)}]
\item 
$\SGL{P_{F_D^\rightarrow}} = [\readingLRTB(T'_D),\readingLRTB(T_D)]_L$ 
and $\SGL{P_{F_D^\downarrow}} = [\readingTBLR(T'_D),\readingTBLR(T_D)]_L$.

\item
$\readingLRTB(T_D)=w_1(\set(\sfr(D)))$, thus $\SGL{P_{F_D^\rightarrow}}$ is an upper descent interval. 

\item
$\readingTBLR(T'_D) = w_0(\set(\sfc(D))^\rmc)$,
thus $\SGL{P_{F_D^\downarrow}}$ is a lower descent interval. 
\end{enumerate}
\end{proposition}
\begin{proof}
(a) 
The first equality was established in \cite[Theorem 5.3]{24CKO}, so we will only prove the second.
By \cref{Lem: PF = PF'}, we have that 
$P_{F_{D}^\downarrow}=P_{(F_{D}^\downarrow)^\rmt}$.
Applying \( (F_{D}^\downarrow)^\rmt=(\overline{F_{D^\rmt}^\rightarrow})\)
and $P_{\overline{F_{D^\rmt}^\rightarrow}}=\overline{P}_{F_{D^\rmt}^\rightarrow}$,
to this equality yields that $P_{F_{D}^\downarrow}=\overline{P}_{F_{D^\rmt}^\rightarrow}$.
On the other hand, from the first equality, it follows that
\[
\SGL{\overline{P}_{F_{D^\rmt}^\rightarrow}} = \SGL{P_{F_{D^\rmt}^\rightarrow}}w_0 = [\readingLRTB(T_{D^\rmt})w_0, \readingLRTB(T'_{D^\rmt})w_0]_L.
\]
Observe that $\readingLRTB(T_{D^\rmt})w_0 = \readingTBLR(T'_{D})$ and 
$\readingLRTB(T'_{D^\rmt})w_0 = \readingTBLR(T_{D})$.
Therefore, 
\[
\SGL{P_{F_{D}^\downarrow}} = [\readingTBLR(T'_{D}),\readingTBLR(T_{D})]_L.
\]

(b) 
Suppose that $k$ is the number of rows in $D$.
For $1 \leq i \leq k-1$, let $\overline{\sfr}_i := r_{i+1} + r_{i+2} + \cdots + r_{k}$. 
By the definitions of \(\readingLRTB\) and \(T_D\), we can express \(\readingLRTB(T_D)\) as follows:
\[
\readingLRTB(T_D) = \underbrace{\overline{\sfr}_1+1 \ \overline{\sfr}_1+2 \ \cdots \ n}_{\text{$1$st row}} \ \underbrace{\overline{\sfr}_2+1 \ \overline{\sfr}_2+2 \ \cdots \ \overline{\sfr}_{1}}_{\text{$2$nd row}} \ \cdots \ \underbrace{1 \ 2 \ \cdots \ r_k}_{\text{$k$th row}}
\]
Here, the notation ``$i$th row" refers to the $i$th row of $T_D$ when counted from the top.
By definition, we have that $\readingLRTB(T_D) = w_1(\set(\sfr(D)))$.
Thus, $\SGL{P_{F_D^\rightarrow}}$ is an upper descent interval. 

(c) This assertion can be proven in a manner similar to (b).
\end{proof}

Next, we show that every lower descent interval is of the form \( \SGL{P_{F_D^\downarrow}} \), and every upper descent interval is of the form \( \SGL{P_{F_D^\rightarrow}} \). 
We begin by reviewing the result in \cite[Section 5.1]{24CKO}. 
For $P\in \poset{n}$, let 
$\SGR{P}:=\{\sigma^{-1} \mid \sigma \in \SGL{P}\}$.

\begin{lemma}{\rm (\cite[Theorem 5.6]{24CKO})}\label{contruction of a diagram for upper lower interval}
Given \(\alpha \models n\) and \(\rho \in \SG_n\) with \(w_0(\set(\alpha)) \preceq_L \rho\), the diagram \(E_{\alpha;\rho} \in \mathfrak{D}_n\) \footnote{It should be remarked that in \cite[Theorem 5.6]{24CKO}, the notation \(D_{\alpha;\rho}\) is employed in place of \(E_{\alpha;\rho}\).}, constructed via \cite[Algorithm 5.4]{24CKO}, satisfies  
\[
\SGR{P_{F^\rightarrow_{E_{\alpha;\rho}}}} = \mapf([w_0(\alpha), \rho]_L),
\]
where \(\mapf: \SG_n \to \SG_n\) is the bijection defined by \(\mapf(\gamma) = w_0 \gamma^{-1}\).
\end{lemma}

In \cref{contruction of a diagram for upper lower interval}, we let
\begin{equation}\label{Eq: def D_sigma_S}
D_{\sigma;S} := E_{\comp(S)^\rmt;\sigma w_0}.
\end{equation}

\begin{example}{\rm (cf. \cite[Example 5.5]{24CKO})}
Let $\sigma = 267935148$ and $S = \{4,6\}$.
Then $\comp(S)^\rmt=(1,1,2,2,1,1,1)\models 9$ and $\sigma w_0 =8 4 1 5 3 9 7 6 2$.
For simplicity, let $\alpha:=\comp(S)^\rmt$ and $\rho:=\sigma w_0$. Then $D_{\sigma;S} = E_{\alpha;\rho}$.
Since $w_0(\set(\alpha)) = 3 2 1 5 4 9 8 7 6$ and $w_0(\alpha) \preceq_L \rho$, we can proceed with the construction of the diagram $E_{\alpha;\rho}$ using \cite[Theorem 5.6]{24CKO}:
\[
E_{\alpha;\rho} = 
\begin{tikzpicture}[baseline=8mm]
\foreach \c in {0,2,4}{
    \filldraw[color=black!15] (\hp*\c,\vp*0) rectangle (\hp*\c+\hp,\vp*1);
    \draw (\hp*\c,\vp*0) rectangle (\hp*\c+\hp,\vp*1);
}
\foreach \c in {1,2}{
    \filldraw [color=black!15] (\hp*\c,\vp*1) rectangle (\hp*\c+\hp,\vp*2);
    \draw (\hp*\c,\vp*1) rectangle (\hp*\c+\hp,\vp*2);
}
\foreach \c in {0,2,3,4}{
    \filldraw[color=black!15] (\hp*\c,\vp*2) rectangle (\hp*\c+\hp,\vp*3);
    \draw (\hp*\c,\vp*2) rectangle (\hp*\c+\hp,\vp*3);
}
\draw[blue,->] (0,0) -- (\hp*6,0) node[anchor=west]{\textcolor{blue}{\footnotesize $x$}};
\draw[blue,->] (0,0) -- (0,\vp*4) node[anchor=south]{\textcolor{blue}{\footnotesize $y$}};
\end{tikzpicture}
\]
\end{example}

In what follows, we present an algorithm that, given \( S \subseteq [n-1] \) and \( \rho \in \SG_n \) satisfying \( w_0(S) \preceq_L \rho \), produces a diagram \( D \in \mathfrak{D}_n \). 
This diagram constitutes a slight modification of the algorithm described in \cite[Algorithm 5.4]{24CKO}.

\begin{algorithm}
\label{Algo: Construction of  D_alpha_rho}
Let $S \subseteq [n-1]$ and $w_0(S) \preceq_L \rho$.

\begin{enumerate}[label = {\it Step \arabic*.}]
\item
Let $S^\rmc = \{z_1,z_2,\ldots,z_p\}$, and $z_0 := 0$ and $z_{p+1}:=n$.
For $1 \leq i \leq p+1$, 
let 
\[
X_i(S;\rho) := \{\rho(r) \mid z_{i-1}+1 \leq r \leq z_i\}.
\]

\item 
Let $e:= |\DesLR{L}{\rho}|$ and let $\DesLR{L}{\rho} = \{q_1 < q_2 < \cdots < q_e\}$, $q_0:=0$, and $q_{e+1}:=n$.
Then, for $j = 1,2,\ldots,e+1$, set
\[
Y_j(\rho) := \{c \mid q_{j-1}+1 \leq c \leq q_j\}.
\]

\item 
Let 
\[
D_{S;\rho} := \{(i,j) \mid X_i(S;\rho) \cap Y_j(\rho) \neq \emptyset \}.
\]
Return $D_{S;\rho}$.
\end{enumerate}  
\end{algorithm}

\begin{example}\label{Ex: S = 25 and its diagram}
Let $S = \{2,5\} \subseteq 6$, and consider $\rho = 231564$.
We apply \cref{Algo: Construction of D_alpha_rho} to obtain the diagram $D_{S;\rho}$.
The observation that $w_0(S) = 132465$ and $w_0(S) \preceq_L \rho$ gives us to proceed with the construction of $D_{S;\rho}$.
Since $S^\rmc = \{1,3,4\}$, {\it Step 1} gives
\[
X_1(S;\rho) = \{2\}, \quad X_2(S;\rho) = \{1,3\} \quad X_3(S;\rho) = \{5\} \quad \text{and} \quad X_4(S;\rho) = \{4,6\}.
\]
With $\DesLR{L}{\rho} =\{1,4\}$, {\it Step 2} gives
\begin{align*}
Y_{1}(\rho) = \{1\}, \quad Y_{2}(\rho) = \{2,3,4\}, \quad \text{and} \quad Y_{3}(\rho) = \{5,6\}.
\end{align*}
Finally, {\it Step 3} gives 
\[
D_{S;\rho} = \{(1,2),(2,1),(2,2),(3,3),(4,2),(4,3)\}.
\]
The diagrams $D_{S;\rho}$ is represented graphically as follows:
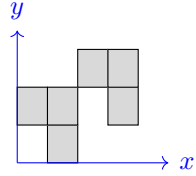
\[
\begin{tikzpicture}
\foreach \c in {1}{
    \filldraw[color=black!15] (\hp*\c,\vp*0) rectangle (\hp*\c+\hp,\vp*1);
    \draw (\hp*\c,\vp*0) rectangle (\hp*\c+\hp,\vp*1);
}
\foreach \c in {0,1,3}{
    \filldraw [color=black!15] (\hp*\c,\vp*1) rectangle (\hp*\c+\hp,\vp*2);
    \draw (\hp*\c,\vp*1) rectangle (\hp*\c+\hp,\vp*2);
}
\foreach \c in {2,3}{
    \filldraw[color=black!15] (\hp*\c,\vp*2) rectangle (\hp*\c+\hp,\vp*3);
    \draw (\hp*\c,\vp*2) rectangle (\hp*\c+\hp,\vp*3);
}
\draw[blue,->] (0,0) -- (\hp*5,0) node[anchor=west]{\textcolor{blue}{\footnotesize $x$}};
\draw[blue,->] (0,0) -- (0,\vp*3.5) node[anchor=south]{\textcolor{blue}{\footnotesize $y$}};
\end{tikzpicture}
\]
\end{example}

\begin{proposition}\label{Lem: diagram alpha rho}
Let $S \subseteq [n-1]$.
\begin{enumerate}[label = {\rm (\alph*)}]  
\item Let \( \sigma \in \SG_n \) satisfy \( w_1(S) \preceq_L \sigma \). The diagram \( D_{\sigma;S} \), constructed using \cite[Algorithm 5.4]{24CKO}, satisfies  
\[
\SGL{P_{F_{D_{\sigma;S}}^\rightarrow}} = [\sigma, w_1(S)]_L.
\]

\item Let \( \rho \in \SG_n \) satisfy \( w_0(S) \preceq_L \rho \). The diagram \( D_{S;\rho} \), constructed using \cref{Algo: Construction of  D_alpha_rho}, satisfies  
\[
\SGL{P_{F_{D_{S;\rho}}^\downarrow}} = [w_0(S), \rho]_L.
\]
\end{enumerate}
\end{proposition}
\begin{proof}
(a) By \cref{contruction of a diagram for upper lower interval}, it follows that  
$\SGL{P_{F^\rightarrow_{E_{\alpha;\rho}}}} = [\rho w_0, w_0(\set(\alpha))w_0]_L$.
Let $\rho w_0=\sigma$ and $w_0(\set(\alpha))w_0=w_1(S)$.
Then, $\rho=\sigma w_0$ and 
\[
w_0(\set(\alpha))=w_1(S)w_0=w_0(S^\rmt).
\]
Now, the assertion follows from the definition that \( D_{\sigma;S} = E_{\comp(S)^\rmt;\sigma w_0} \).

(b) The assertion can be proven using the same approach as in the proof of \cite[Theorem 5.4]{24CKO}, and thus we omit the proof.
\end{proof}

By combining \cref{Lem: diagram alpha rho} with \cref{Prop: two kinds of P_D and intervals}(b) and (c), we obtain the equalities  
\begin{align*}  
&\left\{ \SGL{P_{F_D^\rightarrow}} \mid D \in \mathfrak{D}_n \right\} = \left\{ [\sigma, w_1(S)]_L \mid S \subseteq [n-1], \, \sigma \preceq_L w_1(S) \right\}, \\  
&\left\{ \SGL{P_{F_D^\downarrow}} \mid D \in \mathfrak{D}_n \right\} = \left\{ [w_0(S), \rho]_L \mid S \subseteq [n-1], \, w_0(S) \preceq_L \rho \right\}.
\end{align*}

\subsection{The minimal and maximal elements of the equivalence class} 
\label{The minimal and maximal elements of the equivalence class}
We begin by recalling an essential structural theorem concerning an arbitrary equivalence class in \((\rmInt(n), \Deq)\). 
Let \( C \) be an equivalence class in \((\rmInt(n), \Deq)\). 
Define \(\xi_C := \rho \sigma^{-1}\) for any \( [\sigma, \rho]_L \in C \). 
By \cite[Proposition 4.1]{23KLO}, the permutation \(\xi_C\) is independent of the choice of representative \( [\sigma, \rho]_L \) in \( C \). 
Define a partial order \( \preceq \) on \( C \) by setting
\[
[\sigma, \rho]_L \preceq [\sigma', \rho']_L \quad \text{if and only if} \quad \sigma \preceq_R \sigma'.
\]
As shown in \cite[Theorem 4.6]{23KLO}, the set \(\{\sigma \mid [\sigma, \rho]_L \in C\}\) forms a right weak Bruhat interval. 
Let \(\sigma_0\) and \(\sigma_1\) denote the minimal and maximal elements of this interval, respectively. By definition, 
\begin{equation*}\label{poset str of C}
(C, \preceq) \cong [\sigma_0, \sigma_1]_R \quad \text{as posets}.
\end{equation*}
The minimal and maximal elements of  \((C,\preceq)\) are denoted by \(\min C\) and \(\max C\), respectively, so that \(\min C = [\sigma_0, \xi_C \sigma_0]_L\) and \(\max C = [\sigma_1, \xi_C \sigma_1]_L\).

In general, explicitly describing arbitrary equivalence classes remains an open problem. 
However, this issue was successfully addressed in \cite{23KLO} for the equivalence classes of \( \Sigma_L(P) \), where \( P \) is a regular Schur-labeled poset on \( [n] \). 
Here, when \( C \) is the equivalence class of a lower or upper descent interval, we provide the minimal and maximal elements of \( C \).

We construct a filling \( F_D^{\smallnearrow} \), derived from \( F_D^{\downarrow} \). 
Consider the sequence of fillings 
\begin{equation} \label{definition of $F_i$'s}
Z_0 = F_D^{\downarrow}, \ Z_1, \ Z_2,\ldots, \ Z_n,
\end{equation}
where each \( Z_i \) is obtained from \( Z_{i-1} \) through the following process.
\vspace*{2mm}

\noindent 
\textbf{Case 1:} If there exists an entry \( x > i\) in \( Z_{i-1} \) satisfying the conditions:
\begin{enumerate}[label = {\rm (\roman*)}]
\item \( x \) is positioned strictly above and to the right of \( i \), and 
\item for each \( j = i, i + 1, \ldots, x - 1 \), the entry \( j \) is positioned strictly below and weakly to the left of \( x \), but not in the covering relation in \( P_{Z_{i-1}} \),
\end{enumerate}
then select the uppermost entry among such \( x \)'s, denoted by \( \bfx \). 
Construct \( Z_i \) from \( Z_{i-1} \) by incrementing the entries \( i, i + 1, \ldots, \bfx - 1 \) by \( 1 \), and then swapping the original entry \( \bfx \) with \( i \).
\vspace*{2mm}

\noindent 
\textbf{Case 2:} If no such \( x \) exists, set \( Z_i := Z_{i-1} \).
\smallskip 

It is clear that this process terminates in a finite number of steps. 
We define \( F_D^{\smallnearrow} \) to be the final filling \(Z_n\).

\begin{example}\label{Ex: D and diagram} 
Consider the case where 
\[
\ctab{D}{
\empty & \empty & ~ & ~ \\ 
~ & ~ & \empty & ~ \\ 
\empty & ~ &}
\quad \text{and} \quad 
\ctab{F_D^\downarrow}{
\empty & \empty & 4 & 5 \\ 
1 & 2 & \empty & 6 \\ 
\empty & 3 &}\,\,.
\]
Then the sequence of fillings will be as follows:
\[
Z_0 = Z_1 = Z_2 = F_D^\downarrow, \ 
\ctab{Z_3}{
\empty & \empty & 3 & 5 \\ 
1 & 2 & \empty & 6 \\ 
\empty & 4 &}, \ 
\ctab{Z_4}{
\empty & \empty & 3 & 4 \\ 
1 & 2 & \empty & 6 \\ 
\empty & 5 &}, \ 
\ctab{Z_5}{
\empty & \empty & 3 & 4 \\ 
1 & 2 & \empty & 5 \\ 
\empty & 6 &}
\]
As a step-by-step analysis, 
for \(i = 1\) and \(i = 2\), there is no \(x\) in \(Z_{i-1}\) satisfying conditions (i) and (ii), so \(Z_1 = Z_2 = F_D^\downarrow\).
For \(i = 3\), the entry \(x = 4\) satisfies the conditions, resulting in \(Z_3\).
For \(i = 4\), the entry \(x = 5\) satisfies the conditions, producing \(Z_4\).
For \(i = 5\), the entry \(x = 6\) satisfies the conditions, yielding \(Z_5\).
Thus, the final filling is \(F_D^\smallnearrow = Z_5\). 
\end{example}

The filling $F_D^\smallnearrow$ is evidently a standard filling on $D$.
The following lemma pertains to the structure of \( F_D^{\smallnearrow} \). 

\begin{lemma}\label{Lem: FDdownarrow is well-defined}
Let $D \in \mathfrak{D}_n$.

\begin{enumerate}[label = {\rm (\alph*)}]
\item 
The posets \( P_{F_D^{\downarrow}} \) and \( P_{F_D^{\smallnearrow}} \) are equal as edge-decorated posets.

\item
If \( i \preceq_{P_{F_D^{\smallnearrow}}} i + 1 \), then \((i,i+1)\) is in the covering relation in \( P_{F_D^{\smallnearrow}} \).
\end{enumerate}
\end{lemma}
\begin{proof}
(a) 
Let \(\{Z_j\}_{0 \leq j \leq n}\) be the sequence generating $F_D^\smallnearrow$.
We claim that for each \(1 \leq i \leq n\), the posets \(P_{Z_{i-1}}\) and \(P_{Z_i}\) maintain the same edge-decorated structure (for the definition of $Z_i$, see \eqref{definition of $F_i$'s}).
If {\bf Case 1} is applied to \(Z_{i-1}\), then the entries \(i, i+1, \ldots, \bfx-1\) in \(Z_{i-1}\) are incremented by \(1\) based on conditions (i) and (ii). 
This increment does not alter their relative order. 
Additionally, the original entry \(\bfx\) is replaced by \(i\). 
Since there is an entry greater than \(\bfx\) positioned between \(\bfx\) and \(\bfx-1\) in \(P_{Z_{i-1}}\), as ensured by conditions (i) and (ii), the relative order among these entries also remains unchanged. 
Thus, the claim holds in this case.
If {\bf Case 2} is applied, then $Z_{i} = Z_{i-1}$, so that the relative order and edge decorations clearly remain unchanged. 
Hence, the claim also holds.

For all \(1 \leq i \leq n\), the posets \(P_{Z_{i-1}}\) and \(P_{Z_i}\) retain the same edge-decorated structure. 
As a result, the final poset \(P_{F_D^\smallnearrow}\) and the initial poset \(P_{F_D^\downarrow}\) are identical as edge-decorated posets.

(b) 
Assume, for the sake of contradiction, that there exists a pair \( (i, i+1) \) in \( P_{F_D^\smallnearrow} \) such that \( i \preceq_{P_{F_D^\smallnearrow}} i+1 \), but \( (i, i+1) \) is not in the covering relation in \( P_{F_D^\smallnearrow} \). 
This implies that \( i+1 \) is strictly above and to the right of \( i \) in \( F_D^\smallnearrow \), but there exists an entry \( k \neq i,i+1 \) within the smallest rectangle-shaped subdiagram containing the boxes for \( i \) and \( i+1 \).

Consider the sequence \(\{Z_j\}_{0 \leq j \leq n}\) generating $F_D^\smallnearrow$.
By the construction of $Z_{j}$ from $Z_{j-1}$, we know that the entries  $1,2,\ldots,i+1$ retain their position in $Z_{i+2}, Z_{i+3},\ldots, Z_n$.
Consequently, we have that 
\[
Z_i^{-1}(i) = (F_D^\smallnearrow)^{-1}(i) \quad \text{and} \quad Z_i^{-1}(i+1) = (F_D^\smallnearrow)^{-1}(i+1).
\]
In this situation, the entry $i+1$ satisfies conditions (i) and (ii).
As a result, the entry \(i\) is placed in \( (Z_i)^{-1}(i+1)\), equivalently \(i+1\) is placed in \( (Z_i)^{-1}(i)\). 
By the above argument, we have that 
\[
Z_i^{-1}(i) = (F_D^\smallnearrow)^{-1}(i+1) \quad \text{and} \quad Z_i^{-1}(i+1) = (F_D^\smallnearrow)^{-1}(i),
\]
that is, $i+1 \preceq_{P_{F_D^\smallnearrow}} i$.
This contradicts the assumption, and therefore, no such pair \( (i, i+1) \) can exist in \( P_{F_D^\smallnearrow} \) such that \(i \preceq_{P_{F_D^\smallnearrow}} i+1\) without \( (i, i+1) \) being in the covering relation.
\end{proof}

Given a diagram \( D \in \mathfrak{D}_n \), let \( D^* \) be the diagram obtained by reflecting \( D \) across the line \( y = -x \) and shifting it appropriately to ensure it remains in \( \mathcal{D}_n \). Explicitly, if \( D \) has \( r \) rows and \( c \) columns, then  
\[
D^* = \{(r-j+1, c-i+1) \mid (i, j) \in D\}.
\]
For a standard filling \( F \) on \( D \), define \( F^* \) as the standard filling on \( D^* \), given by  
\[
F^*(i, j) = F(r-j+1, c-i+1),
\] 
where \( r \) and \( c \) are the number of rows and columns of \( D \), respectively.
Using this notation, we define the filling \( F_{D_{\sigma;S}}^{\smallswarrow} \) as  
\begin{equation}\label{Def: FD^swarrow}
F_{D_{\sigma;S}}^{\smallswarrow} := \left(F_{D_{S^\rmc;w_0\sigma}}^{\smallnearrow}\right)^*.
\end{equation}

With these definitions in place, we are ready to state the main result of this subsection.

\begin{theorem}\label{Thm: Descriptions of minC_D and maxC_D}
Let \( S \) be a subset of \( [n-1] \). Suppose that \( \rho, \sigma \in \SG_n \) with \( w_0(S) \preceq_L \rho \) and \( \sigma \preceq_L w_1(S) \). Denote by \( C_{S;\rho} \) and \( C_{\sigma;S} \) the classes of the intervals \( [w_0(S), \rho]_L \) and \( [\sigma, w_1(S)]_L \), respectively. Then, the following results hold:
\begin{enumerate}[label = {\rm (\alph*)}, itemsep = 0.2em]
\item \(\min C_{S;\rho} = [w_0(S), \rho]_L\) and \(\max C_{S;\rho}  = \SGL{P_{F_{D_{S;\rho}}^{\smallnearrow}}}\).

\item \(\min C_{\sigma;S} = \SGL{P_{F_{D_{\sigma;S}}^{\smallswarrow}}}\) and \(\max C_{\sigma;S} = [\sigma, w_1(S)]_L\).
\end{enumerate}
\end{theorem}
\begin{proof}
(a) According to \cref{Lem: diagram alpha rho},  $\SGL{P_{F_{D_{S;\rho}}^\downarrow}} = [w_0(S), \rho]_L$. 
To establish the first equality, it suffices to show that $\min C_{S;\rho}=\SGL{P_{F_{D_{S;\rho}}^\downarrow}}$. 
By \cref{characterization of the equivalence relation} along with \cite[Theorem 4.6]{23KLO}, it is equivalent to verifying that \(P_{F_{D_{S;\rho}}^\downarrow}\) contains no pair $(i, i+1)$ satisfying the following conditions: 
\begin{itemize}
\item $(i, i+1)$ is a comparable pair that is not in the covering relation, and 
\item $i+1$ is less than $i$ in the partial order of  $P_{F_{D_{S;\rho}}^\downarrow}$. 
\end{itemize}
This property follows directly from \cref{def: diagram posets}.  

Similarly, we can derive the second equality using \cref{Lem: FDdownarrow is well-defined}.

(b)
Since \( w_0[\sigma, w_1(S)]_L = [w_0(S^\rmc), w_0 \sigma]_L \), it follows that
\[
I \in C_{\sigma;S} \iff w_0 I \in C_{S^\rmc; w_0 \sigma}
\]
and consequently,
\[
I \Deq I s_i \text{ for } I \in C_{\sigma;S}  \iff w_0 I \Deq w_0 I s_i \text{ for } w_0 I \in C_{S^\rmc;w_0\sigma}.
\]
Here, we use the notation 
$\zeta U:=\{\zeta \sigma \mid \sigma \in U\}$ 
for $U \subseteq \SG_n$ and $\zeta\in \SG_n$. 
This equivalence establishes that the map 
$$C_{\sigma;S} \to C_{S^\rmc;w_0\sigma}, \ I \mapsto w_0 I$$ 
is a bijection.
Moreover, for $F \in \mathfrak{D}_n$, 
it holds that $w_0 \SGL{P_F} = \SGL{P_{F^*}}$.
By applying these two properties to $F_{D_{\sigma;S}}^\smallswarrow$, 
the desired result follows from (a) and \eqref{Def: FD^swarrow}.
\end{proof}

\begin{remark}  \label{distinguished representatives}
For any \( [\sigma, \rho]_L \in C_{S;\rho} \), we have \( \DesLR{L}{\sigma} = \DesLR{L}{w_0(S)} = S \). This implies that \( \sigma \) cannot be of the form \( w_0(J) \) for any \( J \subset [n-1] \) other than \( S \). Therefore, the class \( C_{S;\rho} \) has a unique lower descent interval \( [w_0(S), \rho]_L \). Similarly, for the same reason, \( C_{\sigma;S} \) has a unique upper descent interval \( [\sigma, w_1(S)]_L \).
\end{remark}

\begin{example} 
(a) Let $n=6$, $S = \{2, 5\} \subseteq [5]$, and $\rho = 231564$. Then $w_0(S) \preceq_L \rho$. 
From \cref{{Thm: Descriptions of minC_D and maxC_D}}(a) it follows that 
$\min C_{S;\rho} = [w_0(S), \rho]_L$ and $\max C_{S;\rho} = \SGL{P_{F_{D_{S;\rho}}^{\smallnearrow}}}$. 
It can be easily seen that 
\[
F_{D_{S;\rho}}^\smallnearrow=
\begin{array}{c}
\tableau{
\empty & \empty & 3 & 4 \\ 
1 & 2 & \empty & 5 \\ 
\empty & 6 &}
\end{array} 
\]
Now, by \cref{Rem: reading permutations from F}, we obtain that  
$\min C_{S;\rho} = [132465, 231564]_L$ and $\max C_{S;\rho} = [134652, 235641]_L$.

Let us investigate \( C_{S; \rho} \) in more detail.  
One can observe that \(P_{Z_j}\) \( (1 \leq j \leq 4) \) are all the posets obtained from \( F_{D_{S; \rho}}^\downarrow \) by applying label changes that satisfy the condition in \cref{characterization of the equivalence relation}, where 
\[
\ctab{Z_1:=F_{D_{S;\rho}}^\downarrow}{\empty & \empty & 4 & 5 \\ 
1 & 2 & \empty & 6 \\ 
\empty & 3 &}\,,
\, \,\,  
\ctab{Z_2}{\empty & \empty & 3 & 5 \\ 
1 & 2 & \empty & 6 \\ 
\empty & 4 &}\,, \,\,\,
\ctab{Z_3}{
\empty & \empty & 3 & 4 \\ 
1 & 2 & \empty & 6 \\ 
\empty & 5 &}\,, \text{ and } 
Z_4: = F_{D_{S;\rho}}^\smallnearrow \,\,.
\]
Let $I_j:=\SGL{P_{Z_j}}$ for $1\le j \le 4$. 
By \cref{characterization of the equivalence relation}, 
\( C_{S; \rho} = \{ I_1, I_2, I_3, I_4 \} \), and its poset structure is given as follows:
\[
\begin{tikzpicture}
\node at (\hp*0,\vp*0) (A11) {\( I_1 \)};
\node at (\hp*0,\vp*-2) (A21) {\( I_2 \)};
\node at (\hp*0,\vp*-4) (A31) {\( I_3 \)};
\node at (\hp*0,\vp*-6) (A41) {\( I_4 \)};
\draw[->] (A11) -- (A21) node[right,midway,pos=0.5] {\tiny \( \cdot s_3 \)};
\draw[->] (A21) -- (A31) node[right,midway,pos=0.5] {\tiny \( \cdot s_4 \)};
\draw[->] (A31) -- (A41) node[right,midway,pos=0.5] {\tiny \( \cdot s_5 \)};
\def \hpp {-35mm}
\node at (\hp*0+\hpp,\vp*0) (B11) {\( P_{Z_1} \)};
\node at (\hp*0+\hpp,\vp*-2) (B21) {\( P_{Z_2} \)};
\node at (\hp*0+\hpp,\vp*-4) (B31) {\( P_{Z_3} \)};
\node at (\hp*0+\hpp,\vp*-6) (B41) {\( P_{Z_4} \)};
\draw[->,line width=0.3mm] (B11) -- (B21) node[right,midway,pos=0.5] {\tiny \( \cdot s_3 \)};
\draw[->,line width=0.3mm] (B21) -- (B31) node[right,midway,pos=0.5] {\tiny \( \cdot s_4 \)};
\draw[->,line width=0.3mm] (B31) -- (B41) node[right,midway,pos=0.5] {\tiny \( \cdot s_5 \)};
\draw [->,decorate,decoration={snake,amplitude=.4mm,segment length=3mm,post length=1mm}] (\hpp*0.7,\vp*-3) -- (\hp*-2,\vp*-3) node[above,pos=0.5] {\small $\SGL{\cdot}$};
\end{tikzpicture}
\]
In this figure, the down arrow 
$\begin{tikzpicture}[baseline=0mm]
\draw[->,line width=0.3mm] (\hp*0,\vp*0.6) -- (\hp*0,\vp*-0.2) node[right,midway,pos=0.5] {\tiny \( \cdot s_i \)};
\end{tikzpicture}$ on the left denotes the labeling change $i \leftrightarrow i+1$, and  the down arrow 
$\begin{tikzpicture}[baseline=0mm]
\draw[->] (\hp*0,\vp*0.6) -- (\hp*0,\vp*-0.2) node[right,midway,pos=0.5] {\tiny \( \cdot s_i \)};
\end{tikzpicture}$ on the right denotes the right multiplication by $s_i$.

(b) Let \( n = 6 \), \( S = \{1, 3, 4\} \subseteq [5] \), and \( \sigma = 546213 \). 
Since \( \sigma \preceq_L w_1(S) \), it follows from \cite[Algorithm 5.4]{24CKO} and \eqref{Eq: def D_sigma_S} that we have
\[
\ctab{D_{\sigma;S}}{
\empty & ~ \\ 
\empty & ~ & ~ \\ 
~ \\ 
~ & ~}\, .
\]
So, \cref{Thm: Descriptions of minC_D and maxC_D}(b) yields that $\min C_{\sigma;S} = \SGL{P_{F_{D_{\sigma;S}}^{\smallswarrow}}}$ and $\max C_{\sigma;S} = [\sigma, w_1(S)]_L$.
The fillings \( F_{D_{\sigma;S}}^{\smallswarrow} \) and \( F_{D_{\sigma;S}}^{\rightarrow} \) are given by
\[
F_{D_{\sigma;S}}^{\smallswarrow} =
\begin{array}{c}
\tableau{
\empty & 1 \\ 
\empty & 2 & 6 \\ 
3 \\ 
4 & 5}
\end{array}
\quad \text{and} \quad
F_{D_{\sigma;S}}^{\rightarrow} =
\begin{array}{c}
\tableau{
\empty & 1 \\ 
\empty & 2 & 3 \\ 
4 \\ 
5 & 6}
\end{array}.
\]
From \cref{Rem: reading permutations from F}, we compute
\[
\min C_{\sigma;S} = [542136, 643125]_L \quad \text{and} \quad \max C_{\sigma;S} = [546213, 645312]_L.
\]

By comparing this result with (a), we observe the structure and relationships between the interval classes \( C_{\sigma;S} \) and \( C_{S^\rmc;w_0\sigma} \) through their respective fillings and the poset representations of their elements.
\end{example}

For a complete understanding of the class \( C \), it is essential to know not only \( \min C \) and \( \max C\), but also the poset structure of \( (C, \preceq) \). 
In the remainder of this subsection, we focus on lower and upper descent intervals that satisfy a specific condition. 
These intervals are particularly relevant to the $0$-Hecke modules discussed in the next section.

A diagram \( D \in \mathcal{D}_n \) contains a {\em strictly upper-right pair} if there are two boxes \((x_1, y_1)\) and \((x_2, y_2)\) in \( D \) satisfying:
\begin{enumerate}[label = {\rm (\roman*)}]
\item 
\( x_1 < x_2 \) and \( y_1 < y_2 \), and  

\item
no other boxes of \( D \) lie inside the smallest rectangle enclosing \((x_1, y_1)\) and \((x_2, y_2)\).  
\end{enumerate}
A diagram \( D \in \mathcal{D}_n \) is {\em free of a strictly upper-right configuration} if it has no strictly upper-right pairs.
Pictorially, \( D \) is a diagram that does not contain a subdiagram of the following form
\begin{equation}\label{Fig: these diagram are not allowed}
\begin{tikzpicture}[baseline=10mm]
\def \hp {10mm}
\def \vp {5mm}
\foreach \c in {0,1,2}{
    \foreach \d in {0,1,2,3}{
        \draw[dotted] (\hp*\c,\vp*\d) rectangle (\hp*\c+\hp,\vp*\d+\vp);
}}
\filldraw [color=black!5] (\hp*0,\vp*0) rectangle (\hp*1,\vp*1);
\draw (\hp*0,\vp*0) rectangle (\hp*1,\vp*1) node[xshift=-\hp*0.5,yshift=-\vp*0.5] {\tiny $(x_1,y_1)$};
\filldraw[color=black!5] (\hp*2,\vp*3) rectangle (\hp*3,\vp*4);
\draw (\hp*2,\vp*3) rectangle (\hp*3,\vp*4) node[xshift=-\hp*0.5,yshift=-\vp*0.5] {\tiny $(x_2,y_2)$};
\node at (\hp*3.5,\vp*0.1) {.};
\end{tikzpicture} 
\end{equation}
For the sake of simplicity in notation, let $D^x$ (or $(D)^x$) denote the diagram obtained by reflecting $D$ across the $x$-axis and shifting it as needed, ensuring that it remains an element of $\mathcal{D}_n$.
Similarly, for a filling $F$, let $F^x$ be the filling obtained by reflecting $F$ across the $x$-axis and shifting it upwards appropriately. 
To define a partial order on $\ST(D)$, we set $T \ble U$ if and only if $\readingTBLR(T) \preceq_L \readingTBLR(U)$.
For \(\zeta \in \SG_n\), let \(\zeta \cdot T\) denote the tableau obtained from \(T\) by replacing each entry \(i\) with \(\zeta(i)\) for \(1 \leq i \leq n\).
With these definitions, the following theorem can be stated.

\begin{lemma}\label{properties of free of a strictly upper-right configuration}
Let $D \in \mathfrak{D}_n$ and $T \in \ST(D^x)$. 

\begin{enumerate}[label = {\rm (\alph*)}]
\item  
$T= s_{i_r}\cdots s_{i_2}s_{i_1} \cdot T'_{D^x}$ for some nonnegative integer $r$,
where $i_k$ is strictly upper-left of $i+1$ in $s_{i_{k-1}}\cdots s_{i_2}s_{i_1} \cdot T'_{D_x}$ for $1\le k \le r$.
For the definition of $T'_{D^x}$, refer to \eqref{def of TD and TD prime}.

\item 
Suppose that $D$ is free of a strictly upper-right configuration.
Then, for \( 1 \leq i \leq n-1 \), \( i \) is strictly upper-left of \( i+1 \) in \( T \) if and only if \( (i, i+1) \) is a comparable pair but not a covering relation in \( P_{T^x} \).
\end{enumerate}
\end{lemma}
\begin{proof}
(a)  
We aim to prove that \( T = s_{i_r} \cdots s_{i_1} \cdot T_{D^x}' \), where each \(i_k\) satisfies the conditions stated in this lemma. 
To this end, we define an $H_n(0)$-action on \(\mathbb{C}\ST(D^x)\) as follows:
For $1 \leq i \leq n-1$ and $T \in \ST(\scrD^x)$,
\begin{align}\label{Hn0-action star}
\pi_i \star T:= 
\begin{cases}
s_i \cdot T & \text{if $i$ is strictly upper-left of $i+1$ in $T$}, \\ 
0 & \text{if $i$ is lower-left of $i+1$ in $T$}, \\
T & \text{otherwise.}
\end{cases}
\end{align}
Note that $s_i \cdot T \in \ST(\scrD^x)$ whenever $i$ is strictly upper-left of $i+1$ in $T$.
We need to check that the operators $\pi_i$ satisfy the defining relations of the $0$-Hecke algebra (see \eqref{Hn0-action star}).

We now verify that the operators \(\pi_i\) satisfy the defining relations of the \(0\)-Hecke algebra:
\begin{enumerate}
\item[(i)] \(\pi_i^2 = \pi_i\): \ 
If \(i\) is strictly upper-left of \(i+1\) in \(T\), then \(\pi_i \star T = s_i \cdot T\).
Since \(i\) is strictly lower-right of \(i+1\) in \(s_i \cdot T\), it follows that \(\pi_i \star (s_i \cdot T) = s_i \cdot T\). 
In all other cases, \(\pi_i \star T = 0\) or \(\pi_i \star T = T\), it is clear that \(\pi_i^2  \star T = \pi_i \star T\).

\item[(ii)] \(\pi_i \pi_j = \pi_j \pi_i\) for \(|i - j| > 1\): \ 
If \(|i - j| > 1\), the sets \(\{i, i+1\}\) and \(\{j, j+1\}\) are disjoint, so the actions of \(\pi_i\) and \(\pi_j\) commute.

\item[(iii)] \(\pi_i \pi_{i+1} \pi_i = \pi_{i+1} \pi_i \pi_{i+1}\): \
Consider the relative positions of \(i,i+1,i+2\) in \(T\):
\begin{itemize}
\item \(i\) is strictly upper-left of \(i+1\).
\begin{itemize}[leftmargin=3mm]
\item If $i+2$ is strictly lower-right of $i+1$, then $\pi_i \pi_{i+1} \pi_i\star T = s_i s_{i+1} s_i \cdot T$ and $\pi_{i+1} \pi_{i} \pi_{i+1} \star T = s_{i+1} s_{i} s_{i+1} \cdot T$.

\item If $i+2$ is strictly above and weakly right of $i+1$, then $\pi_i \pi_{i+1} \pi_i\star T  = \pi \pi_{i+1} \star (s_i \cdot T) = 0$ and $\pi_{i+1} \pi_i \pi_{i+1} \star T = 0$.

\item If $i+2$ is strictly upper-left of $i+1$ and strictly lower-right of $i$, then $\pi_i \pi_{i+1} \pi_i\star T  = \pi \pi_{i+1} \star (s_i \cdot T) = \pi \star (s_{i+1}s_i \cdot T) = s_{i+1}s_i \cdot T$ and $\pi_{i+1} \pi_i \pi_{i+1} \star T = \pi_{i+1} \pi_i \star T =  \pi_{i+1} \star (s_i \cdot T) = s_{i+1}s_i \cdot T$.

\item If $i+2$ is upper-right of $i$, then $\pi_i \pi_{i+1} \pi_i\star T  = \pi \pi_{i+1} \star (s_i \cdot T) = 0$ and $\pi_{i+1} \pi_i \pi_{i+1} \star T = \pi_{i+1} \pi_i \star T =  \pi_{i+1} \star (s_i \cdot T) = 0$.

\item If $i+2$ is strictly upper-left of $i$, then $\pi_i \pi_{i+1} \pi_i\star T  = \pi \pi_{i+1} \star (s_i \cdot T) = \pi_i \star (s_i \cdot T) = s_i \cdot T$ and $\pi_{i+1} \pi_i \pi_{i+1} \star T = \pi_{i+1} \pi_i \star T =  \pi_{i+1} \star (s_i \cdot T) = s_i \cdot T$.
\end{itemize}
\newpage

\item \(i\) is lower-left of \(i+1\). In this case $\pi_i \pi_{i+1} \pi_i \star T = 0$.
\begin{itemize}[leftmargin=3mm]
\item If $i+2$ is strictly lower-right of $i$, then
$i+2$ is strictly lower-right  of $i+1$.
So $\pi_{i+1} \pi_{i} \pi_{i+1} \star T = \pi_{i+1} \pi_{i} \star (s_{i+1} \cdot T) = \pi_{i+1} \star (s_i s_{i+1} \cdot T) = 0$.

\item If $i+2$ is upper-right of $i$ and strictly lower-right of $i+1$, then $
\pi_{i+1} \pi_i \pi_{i+1} \star T = \pi_{i+1} \pi_i \star (s_{i+1} \cdot T) = 0$.

\item If $i+2$ is upper-right of $i+1$, then $\pi_{i+1} \pi_i \pi_{i+1} \star T = 0$.

\item If $i+2$ is strictly upper-left of $i+1$, then $\pi_{i+1} \pi_i \pi_{i+1} \star T = \pi_{i+1} \pi_i \star T =  0$.
\end{itemize}

\item \(i\) is strictly lower-right of \(i+1\). \begin{itemize}[leftmargin=3mm]
\item If $i+2$ is strictly lower-right of $i$, then $\pi_{i} \pi_{i+1} \pi_i \star T = \pi_i \pi_{i+1} \star T = \pi_i \star (s_{i+1} \cdot T) = s_i s_{i+1} \cdot T$ and 
$\pi_{i+1} \pi_{i} \pi_{i+1} \star T = \pi_{i+1} \pi_{i} \star (s_{i+1} \cdot T) = \pi_{i+1} \star (s_i s_{i+1} \cdot T) = s_i s_{i+1} \cdot T$.

\item If $i+2$ is upper-right of $i$ and strictly lower-right of $i+1$, then
$\pi_i \pi_{i+1} \pi_i \star T = \pi_i \pi_{i+1} \star T = \pi_i \star (s_{i+1} \star T) = 0$ and 
$\pi_{i+1} \pi_i \pi_{i+1} \star T = \pi_{i+1} \pi_i \star (s_{i+1} \cdot T) = 0$.

\item If $i+2$ is strictly above and weakly left of $i$ and strictly lower-right of $i+1$, then 
$\pi_i \pi_{i+1} \pi_i \star T = \pi_i \pi_{i+1} \star T = \pi_i \star (s_{i+1} \star T) = s_{i+1} \cdot T$ and 
$\pi_{i+1} \pi_i \pi_{i+1} \star T = \pi_{i+1} \pi_i \star (s_{i+1} \cdot T) = \pi_{i+1} \star (s_{i+1} \cdot T) = s_{i+1} \cdot T$.

\item If $i+2$ is upper-right of $i+1$, then $\pi_i \pi_{i+1} \pi_i \star T = \pi_i \pi_{i+1} \star T = 0$ and 
$\pi_{i+1} \pi_i \pi_{i+1} \star T = 0$.

\item If $i+2$ is strictly upper-left of $i+1$, then $\pi_{i} \pi_{i+1} \pi_{i} \star T = T = \pi_{i+1} \pi_i \pi_{i+1} \star T$.
\end{itemize}
\end{itemize}
In all configurations, we verify that $\pi_i \pi_{i+1} \pi_i \star T = \pi_{i+1} \pi_i \pi_{i+1} \star T$ holds.
\end{enumerate}

These verifications confirm that the operators \(\pi_i\) satisfy the \(0\)-Hecke algebra relations. 
Let \(\sfM_{D^x}\) denote the resulting module.
In \(\ST(D^x)\), each tableau, except the tableau \(T'_{D^x}\), contains a pair $(i,i+1)$ such that $i$ is strictly lower-right of $i+1$.
This implies that $\sfM_D$ is generated by $T'_D$, which guarantees that every \(T \in \ST(D^x)\) can be expressed as \(T = s_{i_r} \cdots s_{i_2} s_{i_1} \cdot T_{D^x}'\), which each $i_k$ satisfies the conditions stated in (a).

(b)  
We first address the ``if'' direction.  
By the given assumption, it follows from \cref{def: diagram posets} that \(i\) is positioned to the lower-left of \(i+1\) in \(T^x\).  
Equivalently, \(i\) is positioned to the upper-left of \(i+1\) in \(T\).  
Since \(T\) is a standard tableau, it is impossible for \(i\) and \(i+1\) to occupy the same column in \(T\).  
Furthermore, as the shape of \(T^x\) is free of a strictly upper-right configuration, \(i\) and \(i+1\) cannot occupy the same row in \(T\).  
Thus, \(i\) is strictly above and strictly to the left of \(i+1\) in \(T\).

Next, we consider the ``only if'' direction.  
Under the assumption, \(i+1\) is strictly above and strictly to the right of \(i\) in \(T^x\), implying that \(i \preceq_{T^x} i+1\).  
Since the shape of \(T^x\) is free of a strictly upper-right configuration, there exists a box within the smallest rectangular subdiagram containing the boxes with \(i\) and \(i+1\).  
This ensures that \((i, i+1)\) forms a comparable pair that is not in the covering relation in \(P_{T^x}\).  
\end{proof}

\begin{theorem}\label{Thm: in case of Fsw = Fdown}
Under the same hypothesis as in \cref{Thm: Descriptions of minC_D and maxC_D}, we have the following. 

\begin{enumerate}[label = {\rm (\alph*)}]
\item If $\scrD:=D_{S;\rho}$ is free of a strictly upper-right configuration, then 
\[
\max C_{S;\rho} = \SGL{P_{F_{\scrD}^\rightarrow}}=[\readingLRTB(T'_{\scrD}),w_1(\set(\sfr(\scrD)))]_L.
\] 
Furthermore, $C_{S;\rho} = \{\SGL{P_{T^x}} \mid T \in \ST(\scrD^x)\}$, and $(C_{S;\rho},\preceq) \cong (\ST(\scrD^x),\ble)$ as posets.

\item If $\scrE:=D_{\sigma;S}$ is free of a strictly upper-right configuration, then 
\[
\min C_{\sigma;S} = \SGL{P_{F_{\scrE}^\downarrow}} = [w_0(\set(\sfc(D))^\rmc),\readingTBLR(T_\scrD)]_L.
\]
Furthermore, $C_{\sigma;S} = \{\SGL{P_{T^x}} \mid T \in \ST(\scrE^x)\}$ and $(C_{\sigma;S},\preceq) \cong (\ST(\scrE^x),\ble)$ as posets. 
\end{enumerate}
\end{theorem}

\begin{proof}
Since the proof for (b) is similar to that for (a), we provide the proof for (a) only.

Assume that the diagram \(\mathscr{D}\) is free of a strictly upper-right configuration.
This assumption guarantees that during the construction of \(F_{\mathscr{D}}^\smallnearrow\), the entry \(\bfx\) chosen in \textbf{Case 1} is always the leftmost entry in the uppermost row among all entries greater than \(i-1\) in \(F_{i-1}\). 
This ensures that
$F_{\mathscr{D}}^\smallnearrow = F_{\mathscr{D}}^\rightarrow$.
Based on \cref{Thm: Descriptions of minC_D and maxC_D}(a), this leads to
\[
\max C_{S;\rho} = \SGL{P_{F_{\mathscr{D}}^\rightarrow}}.
\]
Furthermore, by \cref{Prop: two kinds of P_D and intervals}(a) and (b), we have 
\[
\SGL{P_{F_{\mathscr{D}}^\rightarrow}} = [\readingLRTB(T'_{\scrD}), w_1(\set(\sfr(\scrD)))]_L.
\]
Combining these two equalities yields the first assertion.

Next, we prove the second assertion.
Let us consider the map
\[
\phi: \ST(\scrD^x) \to {\rm Int}(n), \quad T \mapsto \SGL{P_{T^x}}.
\]

\begin{enumerate}[label = {\it Claim \arabic*.}]
\item The diagram 
\[
\begin{tikzcd}
T \arrow{r}{\phi} \arrow[d,swap,"s_i \cdot"] & \SGL{P_{T^x}} \arrow[d,"\cdot s_i"] \\ 
s_i \cdot T \arrow[r,"\phi"]  &[2em] \SGL{P_{(s_i \cdot T)^x}} 
\end{tikzcd}
\]
commutes. 
Here, \(i\) is strictly upper-left of \(i+1\) in \(T\).

To verify this claim, consider the poset \(P_{(s_i \cdot T)^x}\), which is obtained from \(P_{T^x}\) by swapping the labels \(i\) and \(i+1\). 
Since \(\scrD\) is free of a strictly upper-right configuration, we see that $P_{(s_i \cdot T)^x} \neq P_{T^x}$.
By \cref{Lem: regular poset and indirectly comparable pair}, if \(P_{T^x}\) is regular, then \(P_{(s_i \cdot T)^x}\) is also regular.
Since $P_{(s_i \cdot T)^x} = s_i \cdot P_{T^x}$, it follows from \eqref{Eq: si P_I equivalent P_I s_i} that 
\[
\SGL{P_{(s_i \cdot T)^x}} = \SGL{s_i \cdot P_{T^x}} = \SGL{P_{T^x}} s_i.
\]

\item The map $\phi$ is injective.

Assume that $\phi(T)=\phi(U)$ for $T,U \in \ST(\scrD^x)$.
By \cref{properties of free of a strictly upper-right configuration}(a), 
$T=\zeta \cdot T_{\scrD^x}'$ and $U=\zeta' \cdot T_{\scrD^x}'$
for some $\zeta,\zeta' \in \SG_n$.
Hence, 
\[
\phi(T)=[w_0(S), \rho]_L \zeta^{-1}= [w_0(S), \rho]_L {\zeta'}^{-1}=\phi(U).
\]
This implies that $\zeta=\zeta'$, and therefore $T=U$.

\item ${\rm Im}(\phi)=C_{S;\rho}$, that is, $C_{S;\rho} = \{\SGL{P_{T^x}} \mid T \in \ST(\scrD^x)\}$.

By \cref{characterization of the equivalence relation},  
every interval in $C_{S;\rho}$ is of the form 
\[
[w_0(S),\rho]_L s_{i_1}s_{i_2}\cdots s_{i_r},
\]
where $(i_{k},i_{k}+1)$ is a comparable pair that is not in the covering relation in $P_{[w_0(S),\rho]_L s_{i_1}s_{i_2}\cdots s_{i_{k-1}}}$ for $1\le k \le r$. 
Observe that $\phi(T_{\scrD^x}') = [w_0(S),\rho]_L \in C_{S;\rho}$.
Let $T= s_{i_r}\cdots s_{i_2}s_{i_1} \cdot T_{\scrD^x}'$ for some nonnegative integer $r$,
where $i_k$ is strictly upper-left of $i_k+1$ in $s_{i_{k-1}}\cdots s_{i_2}s_{i_1} \cdot T_{\scrD^x}'$ for $1\le k \le r$. 
Using {\it Claim 1} repeatedly, we derive that  
\[
\phi(T)= [w_0(S),\rho]_L s_{i_1}s_{i_2}\cdots s_{i_r}.
\]
Now, the desired result follows from \cref{properties of free of a strictly upper-right configuration}(b).
\end{enumerate}

By {\it Claim 2} and {\it Claim 3}, the map $\phi$ induces a bijection $\phi:\ST(\scrD^x) \to C_{S;\rho}$.
Let \( T = s_{i_r} \cdots s_{i_2} s_{i_1} \cdot T_{\scrD^x}' \) for some nonnegative integer \( r \),  
where \( i_k \) is strictly upper-left of \( i_k+1 \) in \( s_{i_{k-1}} \cdots s_{i_2} s_{i_1} \cdot T_{\scrD^x}' \) for \( 1 \leq k \leq r \).  
Similarly, let \( U = s_{j_s} \cdots s_{j_2} s_{j_1} \cdot T_{\scrD^x}' \) for some nonnegative integer \( s \),  
where \( j_k \) is strictly upper-left of \( j_k+1 \) in \( s_{j_{k-1}} \cdots s_{j_2} s_{j_1} \cdot T_{\scrD^x}' \) for \( 1 \leq k \leq s \).  
Suppose that $T \ble U$. 
Since $\readingTBLR(T_{\scrD^x}') = w_0(S)$, it follows from the definition of $\ble$ that 
$$
s_{i_r} \cdots s_{i_2} s_{i_1} w_0(S) \preceq_L s_{j_s} \cdots s_{j_2} s_{j_1} w_0(S),
$$
which implies that
\begin{equation}\label{poset relations preserved}
w_0(S)s_{i_1}s_{i_2} \cdots  s_{i_r}\preceq_R w_0(S)s_{j_1}s_{j_2} \cdots  s_{j_s}.
\end{equation}
On the other hand, {\it Claim 1} says that 
\[
\phi(T)= [w_0(S),\rho]_L s_{i_1}s_{i_2}\cdots s_{i_r} \quad \text{and} \quad 
\phi(U)= [w_0(S),\rho]_L s_{j_1}s_{j_2}\cdots s_{j_s}.
\]
Thus, \eqref{poset relations preserved} ensures that $\phi(T)\preceq_R \phi(U)$.
Therefore, $\phi:(\ST(\scrD^x),\ble) \to (C_{S;\rho},\preceq)$
is a poset isomorphism.
\end{proof}

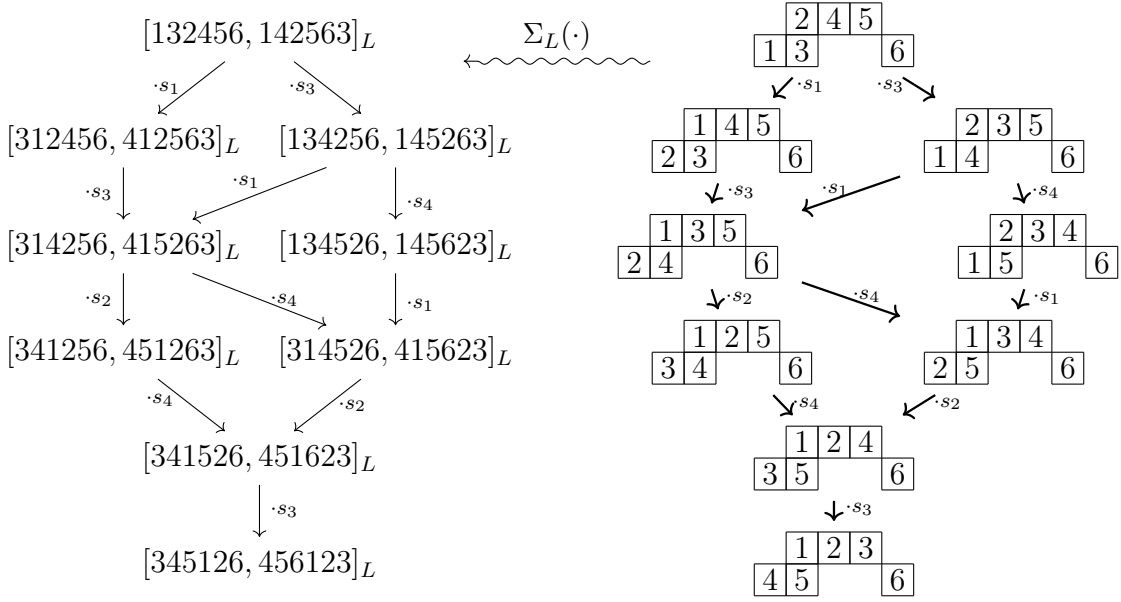
\begin{figure}[t]
\centering
\begin{tikzpicture}
\def \hp {9mm}
\def \vp {14mm}
\node at (\hp*0,\vp*0) (A11)  {$[132456,142563]_L$};
\node at (\hp*-2,\vp*-1) (A21)  {$[312456,412563]_L$};
\node at (\hp*2,\vp*-1) (A22)  {$[134256,145263]_L$};
\node at (\hp*-2,\vp*-2) (A31)  {$[314256,415263]_L$};
\node at (\hp*2,\vp*-2) (A32)  {$[134526,145623]_L$};
\node at (\hp*-2,\vp*-3) (A41)  {$[341256,451263]_L$};
\node at (\hp*2,\vp*-3) (A42)  {$[314526,415623]_L$};
\node at (\hp*0,\vp*-4) (A51)  {$[341526,451623]_L$};
\node at (\hp*0,\vp*-5) (A61)  {$[345126,456123]_L$};

\draw[->,color=black] (A11) -- (A21) node[left,midway] {\tiny $\cdot s_1$};
\draw[->,color=black] (A11) -- (A22) node[left,midway] {\tiny $\cdot s_3$};
\draw[->,color=black] (A21) -- (A31) node[left,midway] {\tiny $\cdot s_3$};
\draw[->,color=black] (A22) -- (A31) node[above,pos=0.6] {\tiny $\cdot s_1$};
\draw[->,color=black] (A22) -- (A32) node[right,pos=0.7] {\tiny $\cdot s_4$};
\draw[->,color=black] (A31) -- (A41) node[left,pos=0.5] {\tiny $\cdot s_2$};
\draw[->,color=black] (A31) -- (A42) node[right,pos=0.5] {\tiny $\cdot s_4$};
\draw[->,color=black] (A32) -- (A42) node[right,pos=0.6] {\tiny $\cdot s_1$};
\draw[->,color=black] (A41) -- (A51) node[left,pos=0.4] {\tiny $\cdot s_4$};
\draw[->,color=black] (A42) -- (A51) node[right,pos=0.5] {\tiny $\cdot s_2$};
\draw[->,color=black] (A51) -- (A61) node[right,pos=0.5] {\tiny $\cdot s_3$};
\def \vp {7mm}
\def \hpp {40mm}
\draw [->,decorate,decoration={snake,amplitude=.4mm,segment length=3mm,post length=1mm}] (\hpp+\hp*1.3,\vp*-0.5) -- (\hp*3,\vp*-0.5) node[above,pos=0.5] {\small $\SGL{\cdot}$};
\node at (\hp*4+\hpp ,\vp*0) (G1) {$\ctabT{
\empty & 2 & 4 & 5 \\ 
1 & 3 & \empty & \empty & 6}$};
\node at (\hp*2.5+\hpp,\vp*-2) (H1) {
$\ctabT{\empty & 1 & 4 & 5 \\ 
2 & 3 & \empty & \empty & 6}$};
\node at (\hp*6.5+\hpp,\vp*-2) (H2) {
$\ctabT{\empty & 2 & 3 & 5 \\
1 & 4 & \empty & \empty & 6}$};
\node at (\hp*2+\hpp,\vp*-4) (I1) {
$\ctabT{
\empty & 1 & 3 & 5 \\ 
2 & 4 & \empty & \empty & 6}$};
\node at (\hp*7+\hpp,\vp*-4) (I2) {$\ctabT{
\empty & 2 & 3 & 4 \\ 
1 & 5 & \empty & \empty & 6}$};
\node at (\hp*2.5+\hpp,\vp*-6) (J1) {
$\ctabT{
\empty & 1 & 2 & 5 \\ 
3 & 4 & \empty & \empty & 6}$};
\node at (\hp*6.5+\hpp,\vp*-6) (J2) {
$\ctabT{ 
\empty & 1 & 3 & 4 \\ 
2 & 5 & \empty & \empty & 6}$};
\node at (\hp*4+\hpp,\vp*-8) (K1) {
$\ctabT{
\empty & 1 & 2 & 4 \\ 
3 & 5 & \empty & \empty & 6}$};
\node at (\hp*4+\hpp,\vp*-10) (L1) {
$\ctabT{ 
\empty & 1 & 2 & 3 \\ 
4 & 5 & \empty & \empty & 6}$};
\draw[->,line width=0.3mm] (G1) -- (H1) node[right,pos=0.4] {\tiny $\cdot s_1$};
\draw[->,line width=0.3mm] (G1) -- (H2) node[left,pos=0.4] {\tiny $\cdot s_3$};
\draw[->,line width=0.3mm] (H1) -- (I1) node[right,pos=0.4] {\tiny $\cdot s_3$};
\draw[->,line width=0.3mm] (H2) -- (I1) node[left,pos=0.4] {\tiny $\cdot s_1$};
\draw[->,line width=0.3mm] (H2) -- (I2) node[right,pos=0.4] {\tiny $\cdot s_4$};
\draw[->,line width=0.3mm] (I1) -- (J1) node[right,pos=0.4] {\tiny $\cdot s_2$};
\draw[->,line width=0.3mm] (I1) -- (J2) node[right,pos=0.4] {\tiny $\cdot s_4$};
\draw[->,line width=0.3mm] (I2) -- (J2) node[right,pos=0.4] {\tiny $\cdot s_1$};
\draw[->,line width=0.3mm] (J1) -- (K1) node[right,pos=0.4] {\tiny $\cdot s_4$};
\draw[->,line width=0.3mm] (J2) -- (K1) node[right,pos=0.4] {\tiny $\cdot s_2$};
\draw[->,line width=0.3mm] (K1) -- (L1) node[right,pos=0.4] {\tiny $\cdot s_3$};
\end{tikzpicture}
\caption{The elements of the class $C_{\{2\};142563}$ (or $C_{345126;\{3\}}$) and the set of standard fillings on $\scrD$}
\label{Fig: the elements of the class and its relationship}
\end{figure}

\begin{example}
Consider two subsets $S = \{2\}$ and $S' = \{ 3\}$ in $[5]$, with permutations $\rho = 142563$ and $\sigma = 345126$ in $\SG_6$.
Observe that $w_0(S) = 132456 \preceq_L \rho$ and $\sigma \preceq_L w_1(S) = 456123$, indicating that both conditions hold. 
Following the construction outlined in \cref{Algo: Construction of D_alpha_rho} and \cite[Algorithm 5.4]{24CKO}, we obtain the diagrams $D_{S;\rho}$ and $D_{\sigma;S}$ as the same diagram
\[
\scrD:=\{(1,1), (2,1), (2,2), (3,2), (4,2), (5,1)\},
\]
which can be visualized as
$\ctabT{
\empty & ~ & ~ & ~ \\
~ & ~ & \empty & \empty & ~}$.
This implies that $\SGL{P_{F_\scrD^\downarrow}} = [w_0(S),142563]_L$ and $\SGL{P_{F_\scrD^\rightarrow}} = [345126,w_1(S)]_L$, where the standard fillings $F_\scrD^\downarrow$ and $F_\scrD^\rightarrow$ are
\[
\ctab{F_\scrD^\downarrow}{
\empty & 2 & 4 & 5 \\
1 & 3 & \empty & \empty & 6}
\quad \text{and} \quad 
\ctab{F_\scrD^\rightarrow}{
\empty & 1 & 2 & 3\\
4 & 5 & \empty & \empty & 6}.
\]
Note that $\scrD$ is free of a strictly upper-right configuration.
According to \cref{Thm: in case of Fsw = Fdown} we have that 
\[
\min C_{S;\rho} = \min C_{\sigma;S} = \SGL{P_{F_{\scrD}^\downarrow}} \quad \text{and} \quad \max C_{S;\rho} = \max C_{\sigma;S} = \SGL{P_{F_{\scrD}^\rightarrow}}.
\]
On the other hand, the elements of $\ST(\scrD^x)$ are
\begin{align*}
& \ctabT{1 & 3 & \empty & \empty & 6\\
\empty & 2 & 4 & 5}
\quad
\ctabT{1 & 4 & \empty & \empty & 6\\
\empty & 2 & 3 & 5}
\quad 
\ctabT{2 & 3 & \empty & \empty & 6\\
\empty & 1 & 4 & 5}
\quad 
\ctabT{2 & 4 & \empty & \empty & 6\\
\empty & 1 & 3 & 5}
\quad \cdots \quad 
\ctabT{4 & 5 & \empty & \empty & 6 \\ 
\empty & 1 & 2 & 3}.
\end{align*} 
By comparing $C_{S;\rho}$ and $\ST(\scrD^x)$ (see \cref{Fig: the elements of the class and its relationship}), we conclude that 
\[
C_{S;\rho} = C_{\sigma;S'} = \{\SGL{P_{T^x}} \mid T \in \ST(\scrD^x)\}. 
\]
\end{example}

\begin{remark}
Suppose that \( D_{S;\rho} \) is free of a strictly upper-right configuration. Then \( C_{S;\rho} \) possesses a unique lower descent interval. This follows by applying \cref{Prop: two kinds of P_D and intervals} and \cref{distinguished representatives} to \cref{Thm: in case of Fsw = Fdown}.
Furthermore, the following identity can be derived:
\[
(C_{S;\rho}, \preceq) \cong ([w_0(S), \readingLRTB(T'_{D_{S;\rho}})]_R,\preceq_R)
\]
by 
\cref{Thm: Descriptions of minC_D and maxC_D} and \cref{Thm: in case of Fsw = Fdown}. 
However, the tableau descriptions in \cref{Thm: in case of Fsw = Fdown} provide additional combinatorial insight and are more convenient from this perspective.
Similarly, if \( D_{\sigma;S} \) is free of a strictly upper-right configuration, then \( C_{\sigma;S} \) possesses a unique upper descent interval. 
Furthermore,
\[
(C_{\sigma;S}, \preceq) \cong ([\readingBTLR(T_{D_{\sigma;S}}'), w_1(S)]_R,\preceq_R).
\]
\end{remark}

It is quite interesting to observe that, in \cref{Thm: in case of Fsw = Fdown}(a), $\max C_{S;\rho}$ is an upper descent interval, and in \cref{Thm: in case of Fsw = Fdown}(b), $\min C_{\sigma;S}$ is a lower descent interval. In the following, we present a representation-theoretic interpretation of this observation. To do this, we introduce the notions of a projective cover and an injective hull of a finitely generated \( H_n(0) \)-module.

Let $A,B$ be finitely generated $H_n(0)$-modules.
A surjective $H_n(0)$-module homomorphism $f:A\to B$ is called an \emph{essential epimorphism} if an $H_n(0)$-module homomorphism $g: X\to A$ is surjective
whenever $f \circ g:X\to B$ is surjective.
A \emph{projective cover} of $A$ is an essential epimorphism $f:P \to A$ with $P$ projective, which always exists and is unique up to isomorphism.
Similarly, let $M,N$ be finitely generated $H_n(0)$-modules with $N \subsetneq M$.
We say that $M$ is an \emph{essential extension} of $N$ if $X\cap N \ne 0$ for all nonzero submodules $X$ of $M$.
An injective $H_n(0)$-module homomorphism $\iota: M \ra \bdI$ with $\bdI$ injective is called an \emph{injective hull} of $M$ if $\bdI$ is an essential extension of $\iota(M)$, which always exists and is unique up to isomorphism.

The $0$-Hecke algebra is a Frobenius algebra, which implies that it is self-injective.
Consequently, finitely generated projective and injective modules coincide
(see \cite[Proposition 4.1]{02DHT}, \cite[Proposition 4.1]{05Fayers}, and \cite[Proposition 1.6.2]{91Benson}).
Furthermore, as shown in \cite[Section 3.2]{22JKLO}, each projective indecomposable module corresponds to a weak Bruhat interval, more precisely, 
for $\alpha \models n$,
\begin{equation*}\label{Eq: bfP and sfB}
\bfP_\alpha \cong \sfB([w_0(\set(\alpha)^\rmc,w_1(\set(\alpha)^\rmc]_L) \quad \text{as $H_n(0)$-modules.}
\end{equation*}

With this preparation, we can derive the following result.

\begin{corollary}\label{how to find injective hull}
Under the same hypothesis as in \cref{Thm: Descriptions of minC_D and maxC_D}, we have the following. 

\begin{enumerate}[label = {\rm (\alph*)}]
\item If $\scrD:=D_{S;\rho}$ is free of a strictly upper-right configuration, then 
\[
\sfB([w_0(\DesLR{L}{\readingBTLR(F_\scrD^\rightarrow)}),w_1(\set(\sfr(\scrD)))]_L)
\]
is an injective hull of \( \sfB([w_0(S), \rho]_L) \).

\item If $\scrE:=D_{\sigma;S}$ is free of a strictly upper-right configuration, then 
\[
\sfB([w_0(\set(\sfc(D))^\rmc),w_1(\DesLR{L}{\readingLRBT(F_\scrE^\downarrow)})]_L).
\]
is a projective cover of $\sfB([\sigma, w_1(S)]_L)$.
\end{enumerate}   
\end{corollary}

\begin{proof}
(a) Since $\scrD$ is free of a strictly upper-right configuration, \cref{Thm: in case of Fsw = Fdown} implies that  
\[
\sfB([w_0(S), \rho]_L) \cong \sfB([\readingBTLR(F_\scrD^\rightarrow)^{-1},\readingLRBT(F_\scrD^\rightarrow)^{-1}]_L) \quad \text{as $H_n(0)$-modules}.
\]
Define $I := \DesLR{L}{\readingBTLR(F_\scrD^\rightarrow})$ and $J := \DesLR{L}{\readingLRBT(F_\scrD^\rightarrow)}$.
By the definition of $F_\scrD^\rightarrow$, we observe that $I \subseteq J$ and $J = \set(\sfr(\scrD))$. 
Since $\readingLRBT(F_\scrD^\rightarrow)^{-1} = w_1(J)$, we have from \cite[Theorem 4.6]{23BS} that  $\sfB([w_0(I),w_1(J)]_L)$ is an injective hull of $\sfB([\readingBTLR(F_\scrD^\rightarrow)^{-1},\readingLRBT(F_\scrD^\rightarrow)^{-1}]_L)$.
Finally, since the injective hull is preserved under $H_n(0)$-module isomorphisms, the result follows.

(b) The proof for (b) follows in a similar manner to that of (a).
\end{proof}

\section{Lower and upper descent intervals from quotient modules and submodules of projective indecomposable \texorpdfstring{$H_n(0)$}{Hn0}-modules}
\label{Lower and Upper Descent Weak Bruhat Intervals from Tableau-Cyclic}

Let \(S \subseteq S' \subseteq [n-1]\). In the context of the representation theory of \(0\)-Hecke algebras, lower and upper descent intervals arise from appropriately selected quotient modules and submodules of the projective module with the basis \([w_0(S), w_1(S')]_L\):
\[
\sfB([w_0(S), w_1(S')]_L) \cong \bigoplus_{\alpha} \bfP_\alpha \quad \text{as $H_n(0)$-modules,}
\]
where $\alpha$ ranges over compositions of $n$ with $S \subseteq {\rm set}(\alpha)^{\rm c} \subseteq S'$.
Specifically, for each \(\rho \in [w_0(S), w_1(S')]_L\), let \(M\) denote the submodule of \(\sfB([w_0(S), w_1(S')]_L)\) generated by the set \([w_0(S), w_1(S')]_L \setminus [w_0(S), \rho]_L\), and let \(N\) represent the submodule generated by \([\rho, w_1(S')]_L\).
Then, we obtain the following \(H_n(0)\)-module isomorphisms:
\[
\sfB([w_0(S), w_1(S')]_L)/M \cong \sfB([w_0(S), \rho]_L)
\]
and 
\[
N \cong \sfB([\rho, w_1(S')]_L).
\]
It should be noted that the latter modules were previously discussed in \cite[Section 4.2]{02DHT}.

In this section, the results presented in \cref{Sec: The equivalence classes of weak Bruhat Intervals} are illustrated through an examination of lower descent intervals that arise from significant quotient modules and submodules of projective indecomposable \(H_n(0)\)-modules. 

\subsection{Lower descent intervals from quotient modules of projective indecomposable \texorpdfstring{$H_n(0)$}{Hn0}-modules}
\label{Lower descent intervals from quotient modules}

We begin by introducing the quotient modules under consideration. 
Let \(\alpha\) be a composition of \(n\).
\begin{itemize}[leftmargin=6mm, itemsep = 0.5em]
\item
In~\cite{15BBSSZ}, Berg--Bergeron--Saliola--Serrano--Zabrocki construct an indecomposable \(H_n(0)\)-module \(\mDIF{\alpha}\) by defining a \(0\)-Hecke action on the set \(\SIT(\alpha)\) of \emph{standard immaculate tableaux of shape \(\alpha\)}. 
The image of this module under the quasisymmetric characteristic is the {\em dual immaculate quasisymmetric function} indexed by \(\alpha\).

\item
In~\cite{19Searles}, Searles constructs an indecomposable \(H_n(0)\)-module \(X_\alpha\) by defining a \(0\)-Hecke action on the set \(\SET(\alpha)\) of \emph{standard extended tableaux of shape \(\alpha\)}. 
The image of this module under the quasisymmetric characteristic is the {\em extended Schur function} indexed by \(\alpha\).

\item
In~\cite[Section 4.2]{22CKNO1}, Choi--Kim--Nam--Oh construct an \(H_n(0)\)-module \(\mYQS{\alpha}\), which is generally not indecomposable, by defining a \(0\)-Hecke action on the set \(\SYCT(\alpha)\) of \emph{standard Young composition tableaux of shape \(\alpha\)}.
\footnote{
There are two remarks. First, a Young composition tableau of shape \(\alpha\) is defined as a filling of \(\tcd(\alpha)\). 
Second, although permuted Young composition tableaux are discussed in~\cite[Section 4.2]{22CKNO1}, only the case \(\sigma = {\rm id}\) is considered in this work.}
The image of this module under the quasisymmetric characteristic is the {\em Young quasisymmetric Schur function} indexed by \(\alpha\). 
We here focus on the canonical submodule \(\mYQS{\alpha,\calC}\) 
of \(\mYQS{\alpha}\) defined in \cite[page. 7767]{22CKNO1}. 
\end{itemize}
All definitions related to \(\SIT(\alpha)\), \(\SET(\alpha)\), \(\SYCT(\alpha)\), the {\em fundamental quasisymmetric function} $F_\alpha$ and the {\em quasisymmetric characteristic} $\ch$ can be found in~\cite[Section 2 and Section 4]{22CKNO1}.

It was shown in \cite[Corollary 4.6]{22CKNO1} that there exists a sequence of surjective \(H_n(0)\)-homomorphisms among these modules, given by
\begin{equation}\label{Eq: sequence of surjective homomorphism}
\begin{tikzcd}
\bfP_\alpha \arrow[r, tail, twoheadrightarrow] &
\mDIF{\alpha} \arrow[r, tail, twoheadrightarrow] & \mESF{\alpha}  \arrow[r, tail, twoheadrightarrow] & \mYQS{\alpha,\calC} \arrow[r, tail, twoheadrightarrow] & \bfF_\alpha.
\end{tikzcd}
\end{equation}
Jung--Kim--Lee--Oh \cite{22JKLO} demonstrated that the intervals arising from these modules are lower descent intervals. In this paper, we show that these intervals satisfy the condition in \cref{Thm: in case of Fsw = Fdown}, which allows us to provide an explicit description of \((C,\preceq)\).
\vspace*{3mm}

\noindent
{\bf Convention.} 
For each \(\bfY_\alpha\) in \cref{Table: quotient modules} and \cref{Table: our consideration}, let \(\calB(\bfY_\alpha)\) denote the basis of tableaux associated with \(\bfY_\alpha\). 
We assume the existence of a reading function 
\[
\mathrm{w}: \calB(\bfY_\alpha) \to \SG_n, \quad T \mapsto \mathrm{w}(T).
\]
From this point onward, we will refer to the set \(\{\mathrm{w}(T) \mid T \in \calB(\bfY_\alpha)\}\) and its equivalence class as \(\RWrM{\mathrm{w}}{\bfY_\alpha}\) and \(C_{\bfY_\alpha}\), respectively.

\begin{table}[ht]	
\centering
\tabulinesep=1.2mm
\small
\begin{tabu}{c|c|c|c}
\tabucline[1.1pt]{-} 
\tiny{\textbf{quasisymmetric functions} $\{y_\alpha\mid \alpha \models n\}$} & {\tiny \textbf{$H_n(0)$-module $\mbalpha$}} & \tiny{$\ch(\mbalpha)$} & {\tiny \textbf{tableau-basis $\calB(\bfY_\alpha)$}}
\\ \hline \hline
{\tiny \text{ribbon Schur func.} $\{s_\alpha\}$ (\cite{99Stanley}) }& \tiny{$\bfP_{\alpha}$ (\cite{79Norton})} & $s_{\alpha^\rmc}$
& 
$\substack{\text{\tiny standard ribbon} \\ \text{\tiny tableaux of shape $\trd(\alpha)$}}$ 
\\ \hline
{\tiny \text{dual immaculate func.} $\{\DIF{\alpha}\}$ (\cite{14BBSSZ})} &\tiny{ $\mDIF{\alpha}$ (\cite{15BBSSZ})}& \tiny{$\DIF{\alpha}$}
& $\substack{\text{standard dual immaculate} \\ \text{ tableaux of shape $\tcd(\alpha)$}}$ 
\\ \hline
{\tiny  \text{extended Schur func.} $\{\ESF{\alpha}\}$ (\cite{19AS})}
&\tiny{ $\mESF{\alpha}$ (\cite{19Searles})} & \tiny{$\ESF{\alpha}$}
& $\substack{\text{standard extended } \\ \text{tableaux of shape $\tcd(\alpha)$}}$ 
\\ \hline
{\tiny \text{Young quasisymm. Schur func.}
$\{\YQS{\alpha}\}$ (\cite{13LMvW})}
& 
\tiny{$\mYQS{\alpha}$ (\cite{22CKNO1})}& \tiny{$\YQS{\alpha}$}
&
$\substack{\text{standard Young composition} \\ \text{ tableaux of shape $\tcd(\alpha)$}}$ 
\\ \hline
{\tiny \text{fundamental quasisymm. func.} $\{F_\alpha\}$ (\cite{99Stanley})}
&\tiny{
$\bfF_{\alpha}$ (\cite{79Norton})}& \tiny{$F_{\alpha^\rmc}$}
&
$\substack{\text{The standard ribbon} \\ \text{tableau $T_\alpha$}}$
\\ \tabucline[1.1pt]{-}
\end{tabu}
\caption{Quasisymmetric functions, associated \( H_n(0) \)-modules, quasisymmetric characteristics, and tableau-bases.}
\label{Table: quotient modules}
\end{table}

\subsubsection{$\bfP_\alpha$}
Recall that 
\begin{align*}
\bfP_\alpha\cong \sfB([w_0(\set(\alpha)^\rmc),w_1(\set(\alpha)^\rmc)]_L) \quad \text{as $H_n(0)$-modules}.
\end{align*}

A {\em standard ribbon tableau} (SRT) of shape $\alpha$ is a filling $T$ of the ribbon diagram $\trd(\alpha)$ with $\{1,2,\ldots,n\}$ such that the entries are all distinct, the entries in each row increase from left to right, and the entries in each column increase from bottom to top.
Let $\SRT(\alpha)$ be the set of SRTs of shape $\alpha$.
Let $T_\alpha \in \SRT(\alpha)$ be the standard ribbon tableau obtained by filling $\trd(\alpha)$ 
with the entries $1,2,\ldots,n$ from bottom to top and from left to right. 
From \cite[Lemma 5.2]{97KT} (or \cite[Section 3.2]{16Huang}) we have
\begin{align}\label{Eq: w(SRT) = descent class}
\RWrM{\readingTBLR}{\bfP_\alpha} = \{\readingTBLR(T) \mid T \in \SRT(\alpha)\} = [w_0(\set(\alpha^\rmc)),w_1(\set(\alpha^\rmc))]_L.
\end{align}

\begin{proposition}\label{Prop: Case of bP_alpha}
Let \(\alpha\) be a composition of \(n\). 
Then  
\[
C_{\bfP_\alpha} = \{[w_0(\set(\alpha)^\rmc),w_1(\set(\alpha)^\rmc)]_L\}. 
\]
\end{proposition}
\begin{proof}
Since $\RWrM{\readingLRTB}{\bfP_\alpha} =  [w_0(\set(\alpha)^\rmc),w_1(\set(\alpha)^\rmc)]_L$, 
the diagram $D_{\set(\alpha)^\rmc;w_1(\set(\alpha)^\rmc)}$ constructed by \cref{Algo: Construction of  D_alpha_rho} is the ribbon diagram $\trd(\alpha)$.
More precisely, 
\[
D_{\set(\alpha)^\rmc;w_1(\set(\alpha)^\rmc)} = \{(i,j+k_i) \mid 1 \leq i \leq \ell(\alpha), \ 1 \leq j \leq \alpha_i\},
\]
where $k_{\ell(\alpha)}:=0$ and  $k_i := \sum_{i < r \leq \ell(\alpha)} (\alpha_r - 1)$ for $1 \leq i \leq \ell(\alpha)-1$.
For example, if $\alpha = (1,2,1,2,1,3)$, then 
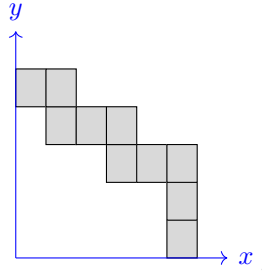
\begin{equation*}\label{Eq: diagrams of 121213}
    \begin{tikzpicture}[baseline = 20mm]
\foreach \c in {4}{
    \filldraw[color=black!15] (\hp*0,\vp*\c) rectangle (\hp*1,\vp*\c+\vp);
    \draw (\hp*0,\vp*\c) rectangle (\hp*1,\vp*\c+\vp);
}
\foreach \c in {3,4}{
    \filldraw [color=black!15] (\hp*1,\vp*\c) rectangle (\hp*2,\vp*\c+\vp);
    \draw (\hp*1,\vp*\c) rectangle (\hp*2,\vp*\c+\vp);
}
\foreach \c in {3}{
    \filldraw[color=black!15] (\hp*2,\vp*\c) rectangle (\hp*3,\vp*\c+\vp);
    \draw (\hp*2,\vp*\c) rectangle (\hp*3,\vp*\c+\vp);
}
\foreach \c in {2,3}{
    \filldraw[color=black!15] (\hp*3,\vp*\c) rectangle (\hp*4,\vp*\c+\vp);
    \draw (\hp*3,\vp*\c) rectangle (\hp*4,\vp*\c+\vp);
}
\foreach \c in {2}{
    \filldraw[color=black!15] (\hp*4,\vp*\c) rectangle (\hp*5,\vp*\c+\vp);
    \draw (\hp*4,\vp*\c) rectangle (\hp*5,\vp*\c+\vp);
}
\foreach \c in {0,1,2}{
    \filldraw[color=black!15] (\hp*5,\vp*\c) rectangle (\hp*6,\vp*\c+\vp);
    \draw (\hp*5,\vp*\c) rectangle (\hp*6,\vp*\c+\vp);
}
\draw[blue,->] (0,0) -- (\hp*7,0) node[anchor=west]{\textcolor{blue}{\footnotesize $x$}};
\draw[blue,->] (0,0) -- (0,\vp*6) node[anchor=south]{\textcolor{blue}{\footnotesize $y$}};
\node at (\hp*8.2,-1.5mm) {.};
\end{tikzpicture}
\end{equation*}
Now, the desired result follows from the property that
\[F^\downarrow_{D_{\set(\alpha)^\rmc;w_1(\set(\alpha)^\rmc)}} = F^\rightarrow_{D_{\set(\alpha)^\rmc;w_1(\set(\alpha)^\rmc)}}.
\]
\end{proof}

\subsubsection{$\bfF_\alpha$}
It is well known that 
\begin{align*}
\bfF_\alpha\cong \sfB([w_0(\set(\alpha)^\rmc),w_0(\set(\alpha)^\rmc)]_L) \quad \text{as $H_n(0)$-modules}
\end{align*}
(for instance, see \cite[Section 3.2]{22JKLO}).
By the definitions of $\bfF_\alpha$ and weak Bruhat interval modules 
in \cref{Subsec: modules of 0Hecke alg} and \cite[Section 3.2]{22JKLO}, 
we see that $\bfF_\alpha \cong \sfB([\sigma,\sigma]_L)$ for all $\sigma \in \SG_n$ with $\set(\alpha)^{\rm c}=\DesLR{L}{\sigma}$.

\begin{proposition}\label{class of singleton}
Let $\alpha$ be a composition of $n$. Then 
\begin{align*}
&C_{\bfF_\alpha}=\{[\sigma,\sigma]_L \mid \sigma \in [w_0(\set(\alpha)^{\rm c}), w_0(\set(\alpha))w_0]_R \},
\text{ and } \\
&
(C_{\bfF_\alpha},\preceq) \cong (\SRT(\alpha),\ble) \quad \text{(as posets).}
\end{align*}
\end{proposition}
\begin{proof}
The first assertions follow from \eqref{variation of descent class}.
For the second assertion, consider the map
\[
[w_0(\set(\alpha)^{\rm c}), w_0(\set(\alpha))w_0]_R \to \SRT(\alpha), \ \sigma \mapsto T_\sigma,
\]
where $T_\sigma$ denotes the SRT of shape $\alpha$ such that $\readingTBLR(T_\sigma) = \sigma^{-1}$ (see \cite[Lemma 5.2]{97KT}).
This map establishes an isomorphism between 
$(C_{\bfF_\alpha},\preceq)$ and $(\SRT(\alpha),\ble)$, as required.
\end{proof}

\subsubsection{\(\mDIF{\alpha}\)}
It was shown in \cite[Section 3.2]{22JKLO} that 
\begin{align*}
\mDIF{\alpha} \cong \sfB([w_0(\set(\alpha)^\rmc),\readingRLBT(\sinkSIT{\alpha})]_L) \quad \text{as $H_n(0)$-modules}.
\end{align*}
Here, \(\sinkSIT{\alpha}\) denotes the standard immaculate tableau of shape \(\alpha\), constructed by first filling the first column with entries \(1, 2, \dots, \ell(\alpha)\) from bottom to top. The remaining boxes are then filled with entries \(\ell(\alpha) + 1, \ell(\alpha) + 2, \dots, n\), moving left to right across each row, starting from the topmost row.

Let $D(\mDIF{\alpha},\readingRLBT)$ be the diagram consisting of 
\[
\{(i,1) \mid (1,i) \in \tcd(\alpha)\} \cup \{(i,j+k_i) \mid (j,i) \in \tcd(\alpha) \text{ and } j \geq 2\},
\]
where $k_{\ell(\alpha)} := 0$ and $k_i:= \sum_{i < r \leq \ell(\alpha)} (\alpha_r-1)$ for $1 \leq i < \ell(\alpha)$.
For example, if $\alpha = (3,2,4)$, then 
\begin{equation}\label{Eq: diagrams of 324}
    \begin{tikzpicture}[baseline = 20mm]
\foreach \c in {0,5,6}{
    \filldraw[color=black!15] (\hp*0,\vp*\c) rectangle (\hp*1,\vp*\c+\vp);
    \draw (\hp*0,\vp*\c) rectangle (\hp*1,\vp*\c+\vp);
}
\foreach \c in {0,4}{
    \filldraw [color=black!15] (\hp*1,\vp*\c) rectangle (\hp*2,\vp*\c+\vp);
    \draw (\hp*1,\vp*\c) rectangle (\hp*2,\vp*\c+\vp);
}
\foreach \c in {0,1,2,3}{
    \filldraw[color=black!15] (\hp*2,\vp*\c) rectangle (\hp*3,\vp*\c+\vp);
    \draw (\hp*2,\vp*\c) rectangle (\hp*3,\vp*\c+\vp);
}
\draw[blue,->] (0,0) -- (\hp*4,0) node[anchor=west]{\textcolor{blue}{\footnotesize $x$}};
\draw[blue,->] (0,0) -- (0,\vp*8) node[anchor=south]{\textcolor{blue}{\footnotesize $y$}};
\node at (\hp*5.2,-1.5mm) {.};
\end{tikzpicture}
\end{equation}

\begin{lemma}{\rm (cf. \cite[Theorem 5.6]{24CKO})}
\label{Lem: SGL PDV = interval from w0 and sinkSIT}
Let $\alpha$ be a composition. Then we have
\[
\SGL{P_{F^\downarrow_{D(\mDIF{\alpha},\readingRLBT)}}} = [w_0(\set(\alpha)^\rmc),\readingRLBT(\sinkSIT{\alpha})]_L.
\]
\end{lemma}
\begin{proof}
By considering the definition of $\sinkSIT{\alpha}$, we observe that for $1 \leq i \leq \ell(\alpha)$, 
\begin{align*}
&(\sinkSIT{\alpha})_{1,i} = i, \quad  \quad (\sinkSIT{\alpha})_{i,j} = j+k_i+(\ell(\alpha)-1) \quad \text{for } 2 \leq j \leq \alpha_i, \quad \text{and} \\ 
&\DesLR{L}{\readingRLBT(\sinkSIT{\alpha})} = \{\ell(\alpha),\ell(\alpha)+1,\ldots,n-1\}.
\end{align*}
Now, applying \cref{Algo: Construction of  D_alpha_rho} to $(\set(\alpha)^\rmc,\readingRLBT(\sinkSIT{\alpha}))$,
we obtain the following sets: 
\begin{align*}
&X_i(\set(\alpha)^\rmc;\readingRLBT(\sinkSIT{\alpha})) = \{i\} \cup \{j+k_i+(\ell(\alpha)-1) \mid 2 \leq j \leq \alpha_i\} \quad \text{for $1 \leq i \leq \ell(\alpha)$}, \\ 
&Y_1(\readingRLBT(\sinkSIT{\alpha})) = \{1,2,\ldots,\ell(\alpha)\}, \quad \text{and} \\  
&Y_j(\readingRLBT(\sinkSIT{\alpha})) = \{j+\ell(\alpha)-1\} \quad \text{for $2 \leq j \leq n-\ell(\alpha)+1$},
\end{align*}
From this data, we can construct the desired diagram as required.
\end{proof}

For an SIT \(\mathcal{T}\), let \(\readingRLBTBT(\mathcal{T})\) denote the reading word obtained from \(\mathcal{T}\) as follows:
\begin{enumerate}[label = {\rm (\roman*)}]
\item Begin by reading the entries in each row from right to left, starting with the bottommost row and moving upward, excluding the entries in the first column.

\item Then, read the entries in the first column from bottom to top.
\end{enumerate}
Let $\sourceSIT{\alpha}$ be the SIT 
shape $\alpha$ obtained by filling $\tcd(\alpha)$ with $1,2,\ldots,n$ from left to right and from bottom to top.
With these definitions, we have the following.

\begin{proposition}\label{Prop: The case of mDIF}
Let $\alpha$ be a composition of $n$, let
and $\scrD:=D(\mDIF{\alpha},\readingRLBT)$.
\begin{enumerate}[label = {\rm (\alph*)}]
\item
$C_{\mDIF{\alpha}} = \{\SGL{P_{T^x}} \mid T \in \ST(\scrD^x)\}$.
Furthermore, 
\[
\min C_{\mDIF{\alpha}} = [w_0(\set(\alpha)^\rmc),\readingRLBT(\sinkSIT{\alpha})]_L 
\ \text{and} \ 
\max C_{\mDIF{\alpha}} = [\readingRLBTBT(\calT_\alpha),w_1([n-\ell(\alpha)])]_L.
\]

\item 
$(C_{\mDIF{\alpha}},\preceq) \cong (\ST(\scrD^x),\boldsymbol{\le})$ as posets.
\end{enumerate}
\end{proposition}
\begin{proof}
(a) By the construction of \(\scrD\), it is evident that \(\scrD\) is free of a strictly upper-right configuration. 
Therefore, the first assertion follows directly from \cref{Thm: in case of Fsw = Fdown}. 
For the second assertion, it follows from \cref{Lem: SGL PDV = interval from w0 and sinkSIT} that
\[
\min C_{\mDIF{\alpha}} = [w_0(\set(\alpha)^\rmc), \readingRLBT(\sinkSIT{\alpha})]_L.
\]
To establish the final equality, consider
\[
\max C_{\mDIF{\alpha}} = \SGL{P_{F_{\scrD}^\rightarrow}} = [\readingLRTB(T_{\scrD}'), \readingLRTB(T_{\scrD})]_L.
\]
Consider the tableau \(T'_{\scrD}\). 
Its reading word, \(\readingLRTB(T'_{\scrD})\), can be expressed as
\[
\underbrace{\alpha_1 \ \alpha_1 - 1 \ \cdots \ 2}_{\text{row 1}} \ \cdots \ \underbrace{n \ n - 1 \ \cdots \ n - \alpha_{\ell(\alpha)} + 2}_{\text{row } \ell(\alpha)} \ \underbrace{1 \ \alpha_1 + 1 \ \cdots \ \alpha_1 + \cdots + \alpha_{\ell(\alpha) - 1} + 1}_{\text{col. 1}}.
\]
Here, `row \(i\)' and `col. \(1\)' refer to the \(i\)th row and the first column in \(\sourceSIT{\alpha}\), respectively. Comparing this word with $\readingRLBTBT(\sourceSIT{\alpha})$, we have that \( \readingLRTB(T'_{\scrD}) = \readingRLBTBT(\sourceSIT{\alpha})\).
On the other hand, we directly have that 
\[
\readingLRTB(T_{\scrD}) = n \ n-1 \ \cdots \ \ell(\alpha) + 1 \ 1 \ 2 \ \cdots \ \ell(\alpha) = w_1([n - \ell(\alpha)]).
\]
This completes the proof.

(b) Since \(\scrD\) is free of a strictly upper-right configuration, it follows from \cref{Thm: Descriptions of minC_D and maxC_D} and \cref{Thm: in case of Fsw = Fdown} that
\[
(C_{\mDIF{\alpha}}, \preceq) \cong (\ST(\scrD^x), \ble) \quad \text{as posets}.
\]
\end{proof}

\begin{example}  
We illustrate \cref{Prop: The case of mDIF} using the diagram \(\scrD\) given in \eqref{Eq: diagrams of 324E}, where \(\alpha = (3, 2, 4)\). 
Consider the tableaux  
\[
\ctab{\sourceSIT{\alpha}}{
6 & 7 & 8 & 9 \\ 
4 & 5 \\ 
1 & 2 & 3} 
\quad \text{and} \quad  
\ctab{\sinkSIT{\alpha}}{
3 & 4 & 5 & 6 \\ 
2 & 7 \\ 
1 & 8 & 9},
\]
then we have 
\begin{align*}
& \readingRLBT(\calT_\alpha) = 321549876, \quad \ \readingRLBT(\calT_\alpha) = 981726543 \quad \text{and} \\ 
& \readingRLBTBT(\calT'_\alpha) = 325987146, \quad \readingRLBTBT(\calT'_\alpha) = 987654123.
\end{align*}
According to it follows from \cref{Prop: The case of mDIF} that we compute
\begin{align*}
& \min C_{\mDIF{\alpha}} = [w_0(\{1,2,4,6,7,8\}),981726543]_L \quad \text{and} \\ 
& \max C_{\mDIF{\alpha}} = [325987146,w_1(\{1,2,3,4,5,6\})]_L.
\end{align*}

On the other hand, since \(\scrD\) is free of a strictly upper-right configuration, it follows from the standard fillings 
\[
F_{\scrD}^\downarrow  = 
\ctabT{
1 \\
2 \\ 
\empty & 4\\ 
\empty & \empty & 6 \\
\empty & \empty & 7 \\ 
\empty & \empty & 8 \\ 
3 & 5 & 9}
\quad \text{and} \quad 
F_{\scrD}^\rightarrow = 
\ctabT{
1 \\
2 \\ 
\empty & 3\\ 
\empty & \empty & 4 \\
\empty & \empty & 5 \\ 
\empty & \empty & 6 \\ 
7 & 8 & 9}
\]
that we have
\begin{align*}
&\min C_{\mDIF{\alpha}} = \SGL{P_{F_{\scrD}^\downarrow}} = [321549876, 981726543]_L \quad \text{and} \\ 
&\max C_{\mDIF{\alpha}} = \SGL{P_{F_{\scrD}^\rightarrow}} = [325987146, 987654123]_L.
\end{align*}

Finally, the elements of $\ST(\scrD^x)$ are
\begin{align*}
\ctabT{3 & 5 & 9 \\ 
\empty & \empty & 8\\
\empty & \empty & 7\\
\empty & \empty & 6\\
\empty & 4\\
2 \\
1}
\quad
\ctabT{4 & 5 & 9 \\ 
\empty & \empty & 8\\
\empty & \empty & 7\\
\empty & \empty & 6\\
\empty & 3\\
2 \\
1}
\quad
\ctabT{4 & 6 & 9 \\ 
\empty & \empty & 8\\
\empty & \empty & 7\\
\empty & \empty & 5\\
\empty & 3\\
2 \\
1}
\quad
\ctabT{4 & 6 & 9 \\ 
\empty & \empty & 8\\
\empty & \empty & 7\\
\empty & \empty & 5\\
\empty & 3\\
2 \\
1}
\quad
\ctabT{4 & 7 & 9 \\ 
\empty & \empty & 8\\
\empty & \empty & 6\\
\empty & \empty & 5\\
\empty & 3\\
2 \\
1}
\quad \cdots \quad 
\ctabT{7 & 8 & 9 \\ 
\empty & \empty & 6\\
\empty & \empty & 5\\
\empty & \empty & 4\\
\empty & 3\\
2 \\
1}.
\end{align*} 
By comparing \(C_{\mDIF{\alpha}}\) and \(\ST(\scrD^x)\), we observe that \((C_{\mDIF{\alpha}},\preceq) \cong (\ST(\scrD^x),\ble)\) as posets.
\end{example}

Combining \cref{how to find injective hull} with \cref{Prop: The case of mDIF} yields the following corollary.

\begin{corollary}{\rm (cf. \cite[Theorem 4.11]{22CKNO2})}
\label{Coro: injective hull of V}
Let $\alpha$ be a composition of $n$ and $\scrD := D(\mDIF{\alpha},\readingRLBT)$.
Then 
\[
\sfB([w_0(I),w_1(\set(\sfr(\scrD)))]_L)
\]
is an injective hull of $\mDIF{\alpha}$,
where 
$I = [n-1]\setminus\{\beta_1,\beta_2, \ldots,\beta_{\ell(\alpha)-1},n-\ell(\alpha), \ldots,n-1\}$ and $\beta_i = \alpha_1+\cdots+\alpha_{i}-i$.
\end{corollary}
\begin{proof}
By \cref{how to find injective hull} we know that 
$\sfB([w_0(\DesLR{L}{\readingBTLR(F_\scrD^\rightarrow)}),w_1(\set(\sfr(\scrD)))]_L)$ is an injective hull of $\mDIF{\alpha}$.
By considering the diagram $\scrD$, it is straightforward to verify that $\DesLR{L}{\readingBTLR(F_\scrD^\rightarrow)} = [n-1]\setminus\{\beta_1,\ldots,\beta_{\ell(\alpha)-1},n-\ell(\alpha),\ldots,n-1\}$, which completes the proof.
\end{proof}
It should be remarked that in \cite[Theorem 4.11]{22CKNO2}, 
an injective hull of $\mDIF{\alpha}$ was directly constructed via an injective map from standard immaculate tableaux
to standard ribbon tableaux.
Compared to this, the current method is more uniform.

\subsubsection{\(\mESF{\alpha}\)}
It was shown in \cite[Section 3.2]{22JKLO} that 
\begin{align*}
\mESF{\alpha} \cong \sfB([w_0(\set(\alpha)^\rmc),\readingRLBT(\sinkSET{\alpha})]_L) \quad \text{as $H_n(0)$-modules}.
\end{align*}
Here, \(\sinkSET{\alpha}\) denotes the standard extended tableau of shape \(\alpha\), constructed by filling with the entries \(1, 2, \dots, \ell(\alpha)\) in each column from bottom to top,  starting from the leftmost column.
Let $D(\mESF{\alpha},\readingRLBT)$ be the diagram $\tcd(\alpha)^\rmt$, that is,
\[
\{(j,i) \mid (i,j) \in \tcd(\alpha)\}.
\]
For example, if $\alpha = (3,2,4)$, then 
\begin{equation}\label{Eq: diagrams of 324E}
\begin{tikzpicture}[baseline = 6mm]
\foreach \c in {0,1,2}{
    \filldraw[color=black!15] (\hp*\c,\vp*0) rectangle (\hp*\c+\hp,\vp*1);
    \draw (\hp*\c,\vp*0) rectangle (\hp*\c+\hp,\vp*1);
}
\foreach \c in {0,1,2}{
    \filldraw [color=black!15] (\hp*\c,\vp*1) rectangle (\hp*\c+\hp,\vp*2);
    \draw (\hp*\c,\vp*1) rectangle (\hp*\c+\hp,\vp*2);
}
\foreach \c in {0,2}{
    \filldraw[color=black!15] (\hp*\c,\vp*2) rectangle (\hp*\c+\hp,\vp*3);
    \draw (\hp*\c,\vp*2) rectangle (\hp*\c+\hp,\vp*3);
}
\foreach \c in {2}{
    \filldraw[color=black!15] (\hp*\c,\vp*3) rectangle (\hp*\c+\hp,\vp*4);
    \draw (\hp*\c,\vp*3) rectangle (\hp*\c+\hp,\vp*4);
}
\draw[blue,->] (0,0) -- (\hp*5,0) node[anchor=west]{\textcolor{blue}{\footnotesize $x$}};
\draw[blue,->] (0,0) -- (0,\vp*4) node[anchor=south]{\textcolor{blue}{\footnotesize $y$}};
\end{tikzpicture} \ .
\end{equation}

\begin{lemma}{\rm (cf. \cite[Theorem 5.8]{24CKO})}
\label{Lem: P_FX_down = interval w0(S)}
Let $\alpha$ be a composition.
Then we have 
\[
\SGL{P_{F^\downarrow_{D(\mESF{\alpha},\readingRLBT)}}} = [w_0(\set(\alpha)^\rmc),\readingRLBT(\sinkSET{\alpha})]_L.
\]
\end{lemma}
\begin{proof}
The assertion follows in a similar manner as in \cref{Lem: SGL PDV = interval from w0 and sinkSIT}.\end{proof}

Let $\sourceSET{\alpha}$ be the standard extended tableau of shape $\alpha$ obtained by filling $\tcd(\alpha)$ with $1,2,\ldots,n$ in each row from left to right, starting from the bottommost row.
Then we have the following.

\begin{proposition}\label{Prop: Case where Y = X}
Let $\alpha$ be a composition of $n$, and let $\scrE:=D(\mESF{\alpha},\readingRLBT)$.
\begin{enumerate}[label = {\rm (\alph*)}]
\item
$C_{\mESF{\alpha}} = \{\SGL{P_{T^x}} \mid T \in \ST(\scrE^x)\}$.
Furthermore, 
\[
\min C_{\mESF{\alpha}} = [w_0(\set(\alpha)^\rmc),\readingRLBT(\sinkSET{\alpha})]_L 
\ \text{and} \ 
\max C_{\mESF{\alpha}} = [\readingBTRL(\sourceSET{\alpha}),w_1(\set((\tal')^\rmr))]_L.
\]

\item 
$(C_{\mESF{\alpha}},\preceq) \cong (\ST(\scrE^x),\ble)$ as posets.
\end{enumerate}
\end{proposition}
\begin{proof}
(a) By the construction of \( \scrE \), it is evident that the diagram \( \scrE \) is free of a strictly upper-right configuration. Therefore, by \cref{Thm: Descriptions of minC_D and maxC_D} and \cref{Thm: in case of Fsw = Fdown}, the first assertion follows.
The second assertion is a direct consequence of \cref{Lem: P_FX_down = interval w0(S)}.
To prove the final assertion, consider the equality
\[
\max C_{\mESF{\alpha}} = \SGL{P_{F_{\scrE}^\rightarrow}} = [\readingLRTB(T_{\scrE}'), \readingLRTB(T_{\scrE})]_L.
\]
By definition, we observe that 
\[
\readingLRTB(T_{\scrE}') = \readingBTRL((T_{\scrE}')^\rmt) \quad \text{and} \quad \readingLRTB(T_{\scrE}) = \readingBTRL((T_{\scrE})^\rmt).
\]
Since $(T'_{\scrE})^\rmt = \sourceSET{\alpha}$ and 
$(T_{\scrE})^\rmt = \sinkSET{\alpha}$, we have 
\[
\max C_{\mESF{\alpha}} = [\readingBTRL(\sourceSET{\alpha},\readingBTRL(\sinkSET{\alpha})]_L.
\]
By the definitions of $\readingLRTB$ and $T_{\scrE}$, since 
\[
\readingLRTB(T_{\scrE}) = \underbrace{\tilde{c}_{\tal'_1}+1 \ \cdots \ n}_{\text{row } 1} \ \underbrace{\tilde{c}_{\tal'_2}+1 \ \cdots \ \tilde{c}_{\tal'_1}}_{\text{row } 2} \ \cdots \ \underbrace{1 \ \cdots \ c_1}_{\text{row } \tal'_1},
\]
where \(c_j\) is the number of boxes in the \(j\)th row of \(\scrE\) (counted from the top) and \(\tilde{c}_j = \sum_{1 \leq i < j} c_i\),
we conclude that \(\readingLRTB(T_{\scrE}) = w_1(\set((\tal')^\rmr))\).
This establishes the desired result.

(b) Since the diagram \(\scrE\) is free of a strictly upper-right configuration, it follows from \cref{Thm: in case of Fsw = Fdown}(a) that 
\[
(C_{\mESF{\alpha}}, \preceq) \cong (\ST(\scrE^x), \ble) \quad \text{as posets}.
\]
\end{proof}

\begin{example}
Given $\alpha = (3,2,4)$, from the tableaux 
\[
\ctab{\sourceSET{\alpha}}{
6 & 7 & 8 & 9 \\ 
4 & 5 \\ 
1 & 2 & 3
}
\quad \text{and} \quad 
\ctab{\sinkSET{\alpha}}{
3 & 6 & 8 & 9 \\ 
2 & 5 \\
1 & 4 & 7
},
\]
we have that
\begin{align*}
& \readingRLBT(\sourceSET{\alpha}) = 321549876, \quad \readingRLBT(\sinkSET{\alpha}) = 741529863 \quad \text{and}\\ 
& \readingBTRL(\sourceSET{\alpha}) = 938257146, \quad \readingBTRL(\sinkSET{\alpha}) = 978456123.
\end{align*}
According to \cref{Prop: Case where Y = X}, the minimal and maximal elements of $C_{\mESF{\alpha}}$ are 
\[
\min C_{\mESF{\alpha}} = [\readingRLBT(\sourceSET{\alpha}),\readingRLBT(\sinkSET{\alpha})]_L 
\quad \text{and} \quad 
\max C_{\mESF{\alpha}} = [\readingBTRL(\sourceSET{\alpha}),\readingBTRL(\sinkSET{\alpha})]_L.
\]

On the other hand, we consider the right-hand diagram $\scrE$ given in \eqref{Eq: diagrams of 324}.
The diagram $\scrE$ is clearly free of a strictly upper-right configuration, so it follows from \cref{Thm: Descriptions of minC_D and maxC_D} that 
\[
F_{\scrE}^\downarrow  = \ctabT{
\empty & \empty & 6 \\
1 & \empty & 7 \\ 
2 & 4 & 8\\ 
3 & 5 & 9}
\quad \text{and} \quad 
F_{\scrE}^\smallnearrow = F_{\scrE}^\rightarrow = \ctabT{
\empty & \empty & 1 \\
2 & \empty & 3 \\ 
4 & 5 & 6 \\ 
7 & 8 & 9}.
\]
Then we compute 
\begin{align*}
&\min C_{\mESF{\alpha}} = \SGL{P_{F_{\scrE}^\downarrow}} = [321549876,74152983]_L = [w_0(\set(\alpha)^\rmc),74152983]_L \quad \text{and}\\ 
&\max C_{\mESF{\alpha}} = \SGL{P_{F_{\scrE}^\rightarrow}} = [938257146,978456123]_L = [938257146,w_1(\set((\tal')^\rmr))]_L.
\end{align*}

\end{example}

\begin{corollary}\label{Coro: injective hull of X}
Let $\alpha$ be a composition of $n$ and $\scrE := D(\mESF{\alpha},\readingRLBT)$.
Then 
\[
\sfB([w_0(\set(\sfr(\scrE))),w_1(\set(\sfr(\scrE)))]_L)
\]
is an injective hull of $\mESF{\alpha}$. 
In particular, it is a projective indecomposable module.
\end{corollary}
\begin{proof}
By \cref{how to find injective hull} we know that 
$\sfB([w_0(\DesLR{L}{\readingBTLR(F_\scrE^\rightarrow)}),w_1(\set(\sfr(\scrE)))]_L)$ is an injective hull of $\mESF{\alpha}$.
By the definition of the diagram $\scrE$, we observe that $\scrE$ contains no subdiagram of the form: 
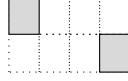
\begin{equation*}
\begin{tikzpicture}[baseline=8mm]
\foreach \c in {0,1,2,3}{
    \foreach \d in {0,1}{
        \draw[dotted] (\hp*\c,\vp*\d) rectangle (\hp*\c+\hp,\vp*\d+\vp);
}}
\filldraw [color=black!15] (\hp*0,\vp*1) rectangle (\hp*1,\vp*2);
\draw (\hp*0,\vp*1) rectangle (\hp*1,\vp*2);
\filldraw[color=black!15] (\hp*3,\vp*0) rectangle (\hp*4,\vp*1);
\draw (\hp*3,\vp*0) rectangle (\hp*4,\vp*1);
\end{tikzpicture} 
\end{equation*}
Using this observation and the definitions of $\readingBTLR$ and $F_\scrE^\rightarrow$, we see that the entry in the rightmost box of row $r$ ($r \geq 2$) in $F_\scrE^\rightarrow$ is contained in $\DesLR{L}{\readingBTLR(F_\scrE^\rightarrow)}$.
Thus we have $\DesLR{L}{\readingBTLR(F_\scrE^\rightarrow)} = \set(\tal^\rmr) = \set(\sfr(\scrE))$, completing the proof.
\end{proof}

\subsubsection{$\mYQS{\alpha,\calC}$}
The source tableau $\sourceSYCT{\alpha,\calC}$ in the canonical class $\calC$ is constructed by filling \(\tcd(\alpha)\) with entries \(1, 2, \dots, n\) in each row from left to right, starting from the bottommost row.
And, \(\sinkSYCT{\alpha,\calC}\) denotes the sink tableau in $\calC$, which can be obtained from the source tableau $\sourceSYCT{\alpha,\calC}$ in $C$ by the algorithm described in \cite[Section 4.1]{22BS} (or \cite[Algorithm 5.13]{24CKO}).
\footnote{In \cite[Section 4.1]{22BS}, the tableau \(\sinkSYCT{\alpha, \calC}\) is referred to as \(T_{\rm sup}\). 
Meanwhile, \cite[Algorithm 5.13]{24CKO} is based on standard reverse composition tableaux (SRCTs). 
To adapt the algorithm for our case, it is necessary to apply the natural bijection between SRCTs and SYCTs}
It was shown in \cite[Section 3.3]{22JKLO} that 
\begin{equation}\label{Eq: SYCT canonical is of the form (interval)}
\mYQS{\alpha,\calC} \cong \sfB([w_0(\set(\alpha)^\rmc),\readingRLBT(\sinkSYCT{\alpha,\calC})]_L) \quad \text{as $H_n(0)$-modules}.
\end{equation}
We here give an algorithm to the diagram from the pair $(w_0(\set(\alpha)^\rmc),\readingRLBT(\sinkSYCT{\alpha,\calC}))$.

\begin{algorithm}\label{Algo: produce the diagram for C}
\hfill 

\begin{enumerate}[label = {\it Step \arabic*.}]
\item 
For each \( 1 \leq i \leq \ell(\alpha) \), define
\[
X_i := \{(1, i), (2, i), \ldots, (\alpha_i, i)\}.
\]

\item Define \( {\sf BB}(\alpha) \) as the set of specific boxes in \( \tcd(\alpha) \) based on the following conditions:
\begin{itemize}
\item Include box \( (1, y) \) in column \( 1 \) if \( 1 \leq y < \ell(\alpha) \) and \( \alpha_y > 1 \), or if \( y = \ell(\alpha) \).

\item Include box \( (x, y) \) for \( x > 1 \) if there is no box strictly above it in the same column or the column immediately to the left.
\end{itemize}

Define the order \( \ble \) on \( {\sf BB}(\alpha) \) such that \( (x, y) \ble (u, v) \) if either \( x = u = 1 \) and \( y < v \), or if \( x < u \).

\item Arrange the elements of \( {\sf BB}(\alpha) \) in increasing order as 
\[
{\sf BB}(\alpha) = \{(b_1, d_1) \ble (b_2, d_2) \ble \cdots \ble (b_l, d_l)\}.
\]
Note that $b_1 = 1$.
Define \( Y_1 := \{(1, d_1), (1, d_1 - 1), \ldots, (1, 1)\} \) and set \( A_1 := \tcd(\alpha) \setminus Y_1 \).

\item For each \( j = 2, 3, \ldots, l \), construct the set \( Y_j \) based on the empty boxes in the diagram \( A_{j-1} \) as follows:
\begin{enumerate}[label = {\rm (\roman*)}]
\item If \( b_{j+1} = 1 \), then define
\[
Y_{j+1} := \{(1, d_{j+1}), (1, d_{j+1} - 1), \ldots, (1, d_j + 1)\}.
\]
Otherwise, set \( Y_{j+1} := \{(b_{j+1}, d_{j+1})\} \).
        
\item Let \( \kappa \) be the last box in \( Y_{j+1} \), and let \( c \) be its column index. 
Check if there is an empty box in column \( c+1 \) of \( A_j \) that is strictly below \( \kappa \). If such a box exists, add the lowermost of these boxes to \( Y_{j+1} \).
Otherwise, stop the process for this \( Y_{j+1} \).
\end{enumerate}
Define \( A_{j+1} := A \setminus Y_{j+1} \) using the updated \( Y_{j+1} \).

\item Using the sets \( X_i \) and \( Y_j \), construct the final diagram \( D_{\alpha,\calC} \) by 
\[
(i, j) \in D_{\alpha,\calC} \quad \text{if } X_i \cap Y_j \neq \emptyset.
\]
And, then return $D_{\alpha,\calC}$.
\end{enumerate}
\end{algorithm}

\begin{lemma} 
Given $\alpha \models n$, \cref{Algo: produce the diagram for C} returns a diagram in $\mathfrak{D}_n$.  
\end{lemma}  
\begin{proof}  
From \cref{Algo: produce the diagram for C}, every box in \(\tcd(\alpha)\) is uniquely assigned to both a set \(X_i\) and a set \(Y_j\). 
Since \(X_i\) forms the boxes in column \(i\) and \(Y_j\) forms the boxes in row \(j\), it follows that \(D_{\alpha, \calC} \in \mathfrak{D}_n\). 
\end{proof}

\begin{example}\label{Ex: class C of 25133}
Applying \cref{Algo: produce the diagram for C} to $\alpha = (2,5,1,3,3)$, one can see that 
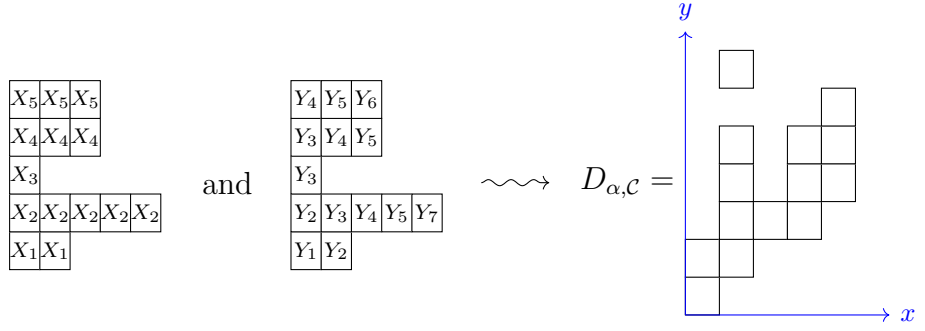
\[
\begin{tikzpicture}[baseline=10mm]
\foreach \c in {0,1}{
    \filldraw[color=black!0] (\hp*\c,\vp*0) rectangle (\hp*\c+\hp,\vp*1);
    \draw (\hp*\c,\vp*0) rectangle (\hp*\c+\hp,\vp*1);
    \node at (\hp*\c+\hp*0.5,\vp*0.5) {\tiny $X_1$};
}
\foreach \c in {0,1,2,3,4}{
    \filldraw[color=black!0] (\hp*\c,\vp*1) rectangle (\hp*\c+\hp,\vp*2);
    \draw (\hp*\c,\vp*1) rectangle (\hp*\c+\hp,\vp*2);
    \node at (\hp*\c+\hp*0.5,\vp*1.5) {\tiny $X_2$};
}
\foreach \c in {0}{
    \filldraw[color=black!0] (\hp*\c,\vp*2) rectangle (\hp*\c+\hp,\vp*3);
    \draw (\hp*\c,\vp*2) rectangle (\hp*\c+\hp,\vp*3);
    \node at (\hp*\c+\hp*0.5,\vp*2.5) {\tiny $X_3$};
}
\foreach \c in {0,1,2}{
    \filldraw[color=black!0] (\hp*\c,\vp*3) rectangle (\hp*\c+\hp,\vp*4);
    \draw (\hp*\c,\vp*3) rectangle (\hp*\c+\hp,\vp*4);
    \node at (\hp*\c+\hp*0.5,\vp*3.5) {\tiny $X_4$};
}
\foreach \c in {0,1,2}{
    \filldraw[color=black!0] (\hp*\c,\vp*4) rectangle (\hp*\c+\hp,\vp*5);
    \draw (\hp*\c,\vp*4) rectangle (\hp*\c+\hp,\vp*5);
    \node at (\hp*\c+\hp*0.5,\vp*4.5) {\tiny $X_5$};
}
\end{tikzpicture}
\quad \text{and} \quad  
\begin{tikzpicture}[baseline=10mm]
\foreach \c in {0,1}{
    \filldraw[color=black!0] (\hp*\c,\vp*0) rectangle (\hp*\c+\hp,\vp*1);
    \draw (\hp*\c,\vp*0) rectangle (\hp*\c+\hp,\vp*1);
}
\foreach \c in {0,1,2,3,4}{
    \filldraw[color=black!0] (\hp*\c,\vp*1) rectangle (\hp*\c+\hp,\vp*2);
    \draw (\hp*\c,\vp*1) rectangle (\hp*\c+\hp,\vp*2);
}
\foreach \c in {0}{
    \filldraw[color=black!0] (\hp*\c,\vp*2) rectangle (\hp*\c+\hp,\vp*3);
    \draw (\hp*\c,\vp*2) rectangle (\hp*\c+\hp,\vp*3);
}
\foreach \c in {0,1,2}{
    \filldraw[color=black!0] (\hp*\c,\vp*3) rectangle (\hp*\c+\hp,\vp*4);
    \draw (\hp*\c,\vp*3) rectangle (\hp*\c+\hp,\vp*4);
}
\foreach \c in {0,1,2}{
    \filldraw[color=black!0] (\hp*\c,\vp*4) rectangle (\hp*\c+\hp,\vp*5);
    \draw (\hp*\c,\vp*4) rectangle (\hp*\c+\hp,\vp*5);
}
\node at (\hp*0.5,\vp*0.5) {\tiny $Y_1$};
\node at (\hp*0.5,\vp*1.5) {\tiny $Y_2$};
\node at (\hp*1.5,\vp*0.5) {\tiny $Y_2$};
\node at (\hp*0.5,\vp*3.5) {\tiny $Y_3$};
\node at (\hp*0.5,\vp*2.5) {\tiny $Y_3$};
\node at (\hp*1.5,\vp*1.5) {\tiny $Y_3$};
\node at (\hp*0.5,\vp*4.5) {\tiny $Y_4$};
\node at (\hp*1.5,\vp*3.5) {\tiny $Y_4$};
\node at (\hp*2.5,\vp*1.5) {\tiny $Y_4$};
\node at (\hp*1.5,\vp*4.5) {\tiny $Y_5$};
\node at (\hp*2.5,\vp*3.5) {\tiny $Y_5$};
\node at (\hp*3.5,\vp*1.5) {\tiny $Y_5$};
\node at (\hp*2.5,\vp*4.5) {\tiny $Y_6$};
\node at (\hp*4.5,\vp*1.5) {\tiny $Y_7$};
\end{tikzpicture}
\quad
\begin{tikzpicture}[baseline=16mm]
\def \hp {4.5mm}
\def \vp {5mm}
\draw [->,decorate,decoration={snake,amplitude=.4mm,segment length=3mm,post length=1mm}] (\hp*-6,\vp*3.5) -- (\hp*-4,\vp*3.5);
\node[left] at (\hp*0,\vp*3.5) {$D_{\alpha,\calC}=$};
\foreach \c in {0}{
    \filldraw[color=black!0] (\hp*\c,\vp*0) rectangle (\hp*\c+\hp,\vp*1);
    \draw (\hp*\c,\vp*0) rectangle (\hp*\c+\hp,\vp*1);
}
\foreach \c in {0,1}{
    \filldraw[color=black!0] (\hp*\c,\vp*1) rectangle (\hp*\c+\hp,\vp*2);
    \draw (\hp*\c,\vp*1) rectangle (\hp*\c+\hp,\vp*2);
}
\foreach \c in {1,2,3}{
    \filldraw[color=black!0] (\hp*\c,\vp*2) rectangle (\hp*\c+\hp,\vp*3);
    \draw (\hp*\c,\vp*2) rectangle (\hp*\c+\hp,\vp*3);
}
\foreach \c in {1,3,4}{
    \filldraw[color=black!0] (\hp*\c,\vp*3) rectangle (\hp*\c+\hp,\vp*4);
    \draw (\hp*\c,\vp*3) rectangle (\hp*\c+\hp,\vp*4);
}
\foreach \c in {1,3,4}{
    \filldraw[color=black!0] (\hp*\c,\vp*4) rectangle (\hp*\c+\hp,\vp*5);
    \draw (\hp*\c,\vp*4) rectangle (\hp*\c+\hp,\vp*5);
}
\foreach \c in {4}{
    \filldraw[color=black!0] (\hp*\c,\vp*5) rectangle (\hp*\c+\hp,\vp*6);
    \draw (\hp*\c,\vp*5) rectangle (\hp*\c+\hp,\vp*6);
}
\foreach \c in {1}{
    \filldraw[color=black!0] (\hp*\c,\vp*6) rectangle (\hp*\c+\hp,\vp*7);
    \draw (\hp*\c,\vp*6) rectangle (\hp*\c+\hp,\vp*7);
}
\draw[blue,->] (0,0) -- (\hp*6,0) node[anchor=west]{\textcolor{blue}{\footnotesize $x$}};
\node[right] at (\hp*6.7,\vp*-0.3) {$.$};

\draw[blue,->] (0,0) -- (0,\vp*7.5) node[anchor=south]{\textcolor{blue}{\footnotesize $y$}};
\end{tikzpicture}
\]
\end{example}

\begin{lemma}\label{Lem: SGL PD_canonical is of the form w0S and sink}
Let $\alpha$ be a composition of $n$ and let $\scrC:=D(\mYQS{\alpha,\calC},\readingRLBT)$.
Then we have 
\[
\SGL{P_{F_{\scrC}^\downarrow}} = [w_0(\set(\alpha)^\rmc),\readingRLBT(\sinkSYCT{\alpha,\calC})]_L. 
\]
\end{lemma}
\begin{proof}
From \cref{Prop: two kinds of P_D and intervals}(a) and (c), we have  
\[
\SGL{P_{F_\scrC^\downarrow}} = [w_0(\set(\sfc(\scrC))^\rmc, \readingTBLR(T_\scrC)]_L.
\]
According to \cref{Algo: produce the diagram for C}, the boxes in \(X_i\) (\(1 \leq i \leq \ell(\alpha)\)) are arranged in the \(i\)th column of \(\scrC\), which implies that \(\sfc(\scrC) = \alpha\). 
Thus, it follows that \(w_0(\set(\sfc(\scrC))^\rmc) = w_0(\set(\alpha)^\rmc)\).

Next, we show that \(\readingTBLR(T_\scrC) = \readingRLBT(\sinkSYCT{\alpha})\). 
Recall the construction of the tableau \(\sinkSYCT{\alpha}\) as described in \cite[Section 4.1]{22BS}, as well as the definition of \(T_\scrC\). 
Comparing this method with \cref{Algo: produce the diagram for C}, one can easily observe that the entries in the \(i\)th row of \(\sinkSYCT{\alpha}\) from right to left correspond precisely to those in the \(i\)th column of \(T_\scrC\) from top to bottom
(refer to \cref{Ex: class C of 25133}).
This implies that \(\readingTBLR(T_\scrC) = \readingRLBT(\sinkSYCT{\alpha})\).
\end{proof}

Given a standard filling $T$ with $1,2,\ldots,n$, let $\overline{T}$ be the filling obtained from $T$ by allocating $\overline{T}_{i,j} = n+1-T_{i,j}$.
For an SYCT $\tau \in \calC$ and a standard filling $R$ of the same shape, let $\readingSYCTR{R}(\tau)$ denote the word obtained by reading the entries of $\tau$ according to the order specified by $R$.
Here, the order of $R$ refers to the sequence dictated by the boxes labeled $1,2,\ldots,n$ in $R$.

\begin{proposition}\label{Prop: Case where Y = hS}
Let $\alpha$ be a composition of $n$.
Let $\scrC:=D_{\alpha,\calC}$.
\begin{enumerate}[label = {\rm (\alph*)}]
\item
$C_{\mYQS{\alpha,\calC}} = \{\SGL{P_{T^x}} \mid T \in \ST(\scrC^x)\}$.
Furthermore, 
\[
\min C_{\mYQS{\alpha,\calC}} = [w_0(\set(\alpha)^\rmc),\readingRLBT(\sinkSYCT{\alpha,\calC})]_L \ \text{and} \ \max C_{\mYQS{\alpha,\calC}} =  [\readingSYCTR{\overline{\sinkSYCT{\alpha,\calC}}}(\sourceSYCT{\alpha,\calC}),w_1(\set(\sfr(\scrC)))]_L.
\]

\item 
$(C_{\mYQS{\alpha,\calC}},\preceq) \cong (\ST(\scrC^x),\ble)$ as posets.
\end{enumerate}
\end{proposition}
\begin{proof}
(a) First, we prove that \(\scrC\) is free of a strictly upper-right configuration. 
To do this, we recall the set ${\sf BB}(\alpha)$ in {\it Step 3} of \cref{Algo: produce the diagram for C}.
Letting $l_\alpha = |{\sf BB}(\alpha)|$, we have 
\[
{\sf BB}(\alpha)=\{(b_1, d_1) \ble (b_2, d_2) \ble \cdots \ble (b_{l_\alpha}, d_{l_\alpha})\}.
\]  
For \(1 \leq j \leq l_\alpha\), define the subdiagram \(\scrC|_{[j]}\) of $\scrC$ by 
\((u,v) \in \scrC|_{[j]}\) if \(X_u \cap Y_v \neq \emptyset\) for \(1 \leq u \leq \ell(\alpha)\) and \(1 \leq v \leq j\). 
For the definition of $X_i$ and $Y_k$, see \cref{Algo: produce the diagram for C}.
By definition, it follows that $\scrC=\scrC|_{[l_\alpha]}$.
The proof proceeds by induction on \( j \). 
If \( j = 1 \), then  
\(\scrC|_{[1]} = \{(1, 1), (1, 2), \ldots, (1, b_1)\}\),  
which is evidently free of a strictly upper-right configuration.
Now, assume that for all \(1 \leq j < l_\alpha\), \(\scrC|_{[j]}\) is free of a strictly upper-right configuration. 
Suppose, for contradiction, that \(\scrC|_{[j+1]}\) has a strictly upper-right pair. 
Then there are two pairs \((x_1, j') \in \scrC|_{[j]}\) and \((x_2, j+1) \in \scrC|_{[j+1]}\) such that
\begin{enumerate}[label = {\rm (\roman*)}]
\item \(x_1 < x_2\), \(j' \leq j < j+1\), and 

\item no box in \(\scrC|_{[j+1]}\) lies within the smallest rectangular subdiagram containing \((x_1, j')\) and \((x_2, j+1)\).
\end{enumerate} 
Note that the boxes \((x_1, j')\) and \((x_2, j+1)\) are constructed from $(X_{x_1},Y_{j'})$ and $(X_{x_2},Y_{j+1})$, respectively.
Let \(c_1\) and \(c_2\) be the column indices of boxes in \(\tcd(\alpha)\) corresponding to $(X_{x_1},Y_{j'})$ and $(X_{x_2},Y_{j+1})$, respectively. 
Since \(x_1 < x_2\), we have two cases:
\begin{itemize}
\item{\it \(c_1 < c_2\)}: \
Here, \((c_2-1, x_2)\) lies in \(\tcd(\alpha)\) and belongs to \(Y_q\) for some \(j' \leq q < j+1\).  
This implies that \((x_2, q)\) is a box of the rectangular diagram determined by \((x_1, j')\) and \((x_2, j+1)\), contradicting the assumption.
 
\item{\it \(c_1 > c_2\)}: \
Let \(x_2'\) be the largest row index such that \((c_2, x_2') \in \tcd(\alpha)\) and \(x_2' < x_2\).  
Here, \((c_2, x_2')\) belongs to \(Y_q\) for some \(j' \leq q < j+1\).  
This again implies that \((x_2, q)\) is a box of the rectangular diagram determined by \((x_1, j')\) and \((x_2, j+1)\), leading to a contradiction.
\end{itemize}
Thus, \(\scrC|_{[j+1]}\) is free of a strictly upper-right configuration.
Thus, by induction, it follows that \(\scrC|_l = \scrC\) is free of a strictly upper-right configuration. 
Consequently, we derive the first equality
\[
C_{\mYQS{\alpha,\calC}} = \{\SGL{P_{T^x}} \mid T \in \ST(\scrC^x)\}.
\]
The second equality follows directly from \cref{Thm: in case of Fsw = Fdown}(b) and \cref{Lem: SGL PD_canonical is of the form w0S and sink}.
Finally, recall that 
\[
\max C_{\mYQS{\alpha,\calC}} = \SGL{P_{F_\scrC^\rightarrow}} = [\readingLRTB(T'_\scrC), w_1(\set(\sfr(\scrC)))]_L,
\]
as shown in \cref{Thm: in case of Fsw = Fdown}(a). 
Therefore, the third equality follows from the equality
\[
\readingLRTB(T'_\scrC) = \readingSYCTR{\overline{\sinkSYCT{\alpha,\calC}}}(\sourceSYCT{\alpha,\calC}).
\]

(b) This assertion directly follows from \cref{Thm: in case of Fsw = Fdown}.
\end{proof}

\begin{example}
Let $\alpha = (2,5,1,3,3) \models 14$.
We validate \cref{Prop: Case where Y = hS}(a) using the tableaux \(\sourceSYCT{\alpha, \calC}\), \(\sinkSYCT{\alpha, \calC}\), and \(\overline{\sinkSYCT{\alpha, \calC}}\):
\[
\ctab{\sourceSYCT{\alpha,\calC}}{
12 & 13 & 14 \\ 
9 & 10 & 11 \\
8\\
3 & 4 & 5 & 6 & 7 \\ 
1 & 2 
}
\qquad 
\ctab{\sinkSYCT{\alpha,\calC}}{
9 & 12 & 13 \\ 
6 & 8 & 11 \\
5 \\
3 & 4 & 7 & 10 & 14 \\ 
1 & 2
}
\qquad 
\ctab{\overline{\sinkSYCT{\alpha,\calC}}}{
6 & 3 & 2 \\ 
9 & 7 & 4 \\
10 \\
12 & 11 & 8 & 5 & 1 \\ 
14 & 13
}
\]
From the definitions of \(\readingRLBT\) and \(\readingSYCTR{\overline{\sinkSYCT{\alpha, \calC}}}\), we compute that
\begin{align*}
&[\readingRLBT(\sourceSYCT{\alpha, \calC}), \readingRLBT(\sinkSYCT{\alpha, \calC})]_L = 
[w_0(\{1,3,4,5,6,9,10,12,13\}), \, 2 \, 1 \, 14 \, 10 \, 7 \, 4 \, 3 \, 5 \, 11 \, 8 \, 6 \, 13 \, 12 \, 9]_L, \\
&[\readingSYCTR{\overline{\sinkSYCT{\alpha, \calC}}}(\sourceSYCT{\alpha, \calC}), \readingSYCTR{\overline{\sinkSYCT{\alpha, \calC}}}(\sinkSYCT{\alpha, \calC})]_L = 
[7 \, 14 \, 13 \, 11 \, 6 \, 12 \, 10 \, 5 \, 9 \, 8 \, 4 \, 3 \, 2 \, 1, w_1(\{1,2,5,8,11,13\})]_L.
\end{align*}

On the other hand, since the diagram \(\scrC\) associated with the composition \(\alpha\) and \(\calC\), as described in \cref{Ex: class C of 25133} is free of a strictly upper-right configuration, it follows from \cref{Thm: in case of Fsw = Fdown} that the minimal and maximal elements of \(C_{\mYQS{\alpha, \calC}}\) are
\begin{align*}
& \min C_{\mYQS{\alpha,\calC}} = \SGL{P_{F_{\scrC}^\downarrow}} =  
[w_0(\{1,3,4,5,6,9,10,12,13\}), \, 2 \, 1 \, 14 \, 10 \, 7 \, 4 \, 3 \, 5 \, 11 \, 8 \, 6 \, 13 \, 12 \, 9]_L, \\ 
& \max C_{\mYQS{\alpha,\calC}} = \SGL{P_{F_{\scrC}^\rightarrow}} = 
[7 \, 14 \, 13 \, 11 \, 6 \, 12 \, 10 \, 5 \, 9 \, 8 \, 4 \, 3 \, 2 \, 1, w_1(\{1,2,5,8,11,13\})]_L.
\end{align*}
Here,  
\[
F_{\scrC}^\downarrow  = 
\ctabT{
\empty & 3 \\ 
\empty & \empty & \empty & \empty & 12 \\ 
\empty & 4 & \empty & 9 & 13\\ 
\empty & 5 & \empty & 10 & 14 \\ 
\empty & 6 & 8 & 11 \\ 
1 & 7 \\ 
2}
\quad \text{and} \quad 
F_{\scrC}^\smallnearrow = F_{\scrC}^\rightarrow = 
\ctabT{
\empty & 1 \\ 
\empty & \empty & \empty & \empty & 2 \\ 
\empty & 3 & \empty & 4 & 5 \\ 
\empty & 6 & \empty & 7 & 8 \\ 
\empty & 9 & 10 & 11 \\ 
12 & 13 \\ 
14}.
\]
Thus, by comparing these results, we confirm that
\[
\min C_{\mYQS{\alpha,\calC}} = [\readingRLBT(\sourceSYCT{\alpha,\calC}), \readingRLBT(\sinkSYCT{\alpha,\calC})]_L \quad \text{and} \quad \max C_{\mYQS{\alpha,\calC}} = [\readingSYCTR{\overline{\sinkSYCT{\alpha,\calC}}}(\sourceSYCT{\alpha,\calC}), \readingSYCTR{\overline{\sinkSYCT{\alpha,\calC}}}(\sinkSYCT{\alpha,\calC})]_L.
\]
\end{example}

We have the following from \cref{how to find injective hull}.

\begin{corollary}\label{Coro: injective hull of hS}
Let $\alpha$ be a composition of $n$ and $\scrC := D(\mYQS{\alpha,\calC},\readingRLBT)$.
Then 
\[
\sfB([w_0(\DesLR{L}{\readingBTLR(F_\scrC^\rightarrow)}),w_1(\set(\sfr(\scrC)))]_L)
\]
is an injective hull of $\mYQS{\alpha,\calC}$.
\end{corollary}

\subsection{Upper descent intervals from submodules of projective indecomposable \texorpdfstring{$H_n(0)$}{Hn0}-modules I}
\label{Upper descent intervals from submodules I}

We begin by introducing the submodules under consideration.
Let \(\alpha\) be a composition of \(n\). 
\begin{itemize}[leftmargin=6mm, itemsep = 0.5em]
\item
In~\cite[Section 6]{24NSVVW}, Niese--Sundaram--van Willigenburg--Vega--Wang constructs an indecomposable \(H_n(0)\)-module \(\mRDIF{\alpha}\)\footnotemark[5]
by defining a \(0\)-Hecke action on the set of \emph{standard dual immaculate tableaux of shape \(\alpha\)}. 
The image of this module under the quasisymmetric characteristic is the {\em row-strict dual immaculate quasisymmetric function} indexed by \(\alpha\).

\item
In~\cite[Section 7]{24NSVVW}, Niese--Sundaram--van Willigenburg--Vega--Wang construct an indecomposable \(H_n(0)\)-module \(\mRESF{\alpha} \)
\footnote{In~\cite{24NSVVW}, the modules \( \mDIF{\alpha} \), \( \mRDIF{\alpha} \) and \(\mRESF{\alpha}\) are denoted by \( \mathcal{W}_\alpha \), \( \mathcal{V}_\alpha \) and \( \mathcal{Z}_\alpha \), respectively.}
by defining a \(0\)-Hecke action on the set of {\em  standard extended tableaux of shape \(\alpha\)}. 
The image of this module under the quasisymmetric characteristic is the {\em row-strict extended Schur function} indexed by \(\alpha\).

\item
In~\cite[Section 3]{22BS}, Bardwell--Searles construct an \(H_n(0)\)-module \(\mRQS{\alpha}\)\footnote{In~\cite{22BS}, the module \( \mRQS{\alpha} \) is denoted by \( \mathbf{R}_\alpha \). }
by defining a \(0\)-Hecke action on the set of {\em standard Young row-strict  tableaux of shape \(\alpha\)}. 
The image of this module under the quasisymmetric characteristic is the {\em Young row-strict quasisymmetric Schur function} indexed by \(\alpha\).
\end{itemize}

\begin{table}[ht]	
\centering
\tabulinesep=1.2mm
\small
\begin{tabu}{c|c|c|c}
\tabucline[1.1pt]{-} 
{\tiny \textbf{quasisymmetric functions $\{y_\alpha \mid \alpha \models n\}$}}
& {\tiny \textbf{$H_n(0)$-module $\mbalpha$}}
& {\tiny $\ch(\mbalpha)$} 
& {\tiny \textbf{tableau-basis} {\tiny $\calB(\bfY_\alpha)$}}
\\ \hline
{\tiny \text{row-strict dual immaculate func.} $\{\RDIF{\alpha}\}$ (\cite{22NSvWVW})} 
& {\tiny $\mRDIF{\alpha}$ (\cite{24NSVVW})}
& {\tiny $\RDIF{\alpha}$}
& $\substack{\text{standard dual immaculate} \\ \text{tableaux of shape $\tcd(\alpha)$}}$ 
\\ \hline
{\tiny \text{row-strict extended Schur func.} $\{\RESF{\alpha}\}$ (\cite{22NSvWVW})}
& {\tiny $\mRESF{\alpha}$ (\cite{24NSVVW})}
& {\tiny $\RESF{\alpha}$}  
& $\substack{\text{standard extended} \\ \text{tableaux of shape $\tcd(\alpha)$}}$ 
\\ \hline
{\tiny \text{row-strict Young quasisymm. Schur func.
$\{\RQS{\alpha}\}$ (\cite{15MN})}}
& {\tiny $\mRQS{\alpha}$ (\cite{22BS})}
& {\tiny $\RQS{\alpha}$}
& $\substack{\text{standard Young composition} \\ \text{ tableaux of shape $\tcd(\alpha)$}}$ 
\\ \tabucline[1.1pt]{-}
\end{tabu}
\caption{Quasisymmetric functions, associated $H_n(0)$-modules, quasisymmetric characteristics, and tableau-bases}
\label{Table: our consideration}
\end{table}

To begin with, we observe that these modules can be obtained from the quotient modules in \cref{Lower descent intervals from quotient modules} by applying an anti-involution.
Let us review the definitions of (anti-) involution twists.
Given an automorphism $\mu$ of $H_n(0)$ and a left $H_n(0)$-module $M$, we define $\mu[M]$ by the left $H_n(0)$-module with the same underlying space as $M$ and with the action $\cdot_\mu$ defined by
\[
h \cdot_\mu v:= \mu(h) \cdot v \quad \text{for $h \in H_n(0)$ and $v \in M$}.
\]
We define $\bfT^+_\mu : \Rmod \to \Rmod$ to be the covariant functor, called the \emph{$\mu$-twist}, sending a left $H_n(0)$-module $M$ to $\mu[M]$ and an $H_n(0)$-module homomorphism $f:M \to N$ to $\mathbf{T}^+_\mu(f): \mu[M] \to \mu[N]$ defined by $\mathbf{T}^+_\mu(f)(v) = f(v)$ for $v \in M$.

Similarly, given an anti-automorphism $\nu$ of $H_n(0)$ and a left $H_n(0)$-module $M$, we define $\nu[M]$ by the left $H_n(0)$-module with $M^*$, the dual space of $M$, as the underlying space and with the action $\cdot^\nu$ defined by 
\begin{align}\label{Eq: T_nu^minus}
(h \cdot^\nu \delta)(v) := \delta(\nu(h) \cdot v)
\quad \text{for $h \in H_n(0)$, $\delta \in M^*$, and $v \in M$.}
\end{align}
We define $\bfT^-_\nu: \Rmod \to \Rmod$ to be the contravariant functor, called the \emph{$\nu$-twist}, sending an $H_n(0)$-module $M$ to $\nu[M]$ and an $H_n(0)$-module homomorphism $f:M \to N$ to $\bfT^-_\nu(f): \nu[N] \to \nu[M]$ defined by $\bfT^-_\nu(f)(\delta) = \delta \circ f$.

In this subsection, we consider two involutions $\upphi, \uptheta$ and an anti-involution $\upchi$ on $H_n(0)$ defined by 
\[
\upphi(\opi_i) = \opi_{n-i}, 
\quad 
\uptheta (\opi_i) = -\pi_i, 
\quad \text{and} \quad 
\upchi(\opi_i) = \opi_i \quad (1 \le i \le n-1).
\]
These (anti-)involutions were introduced by Fayers \cite[Proposition 3.2]{05Fayers} and commute with each other.

\begin{lemma}{\rm (cf. \cite[Table 2 and Section 4]{22JKLO})}
\label{involutive images of quotient modules}
Let $\alpha$ be a composition of $n$. 
Then we have the following $H_n(0)$-module isomorphisms:  
\begin{align}
 & \label{twist isom 1}\uptheta \circ \upchi[\bfP_\alpha] \cong \bfP_{\alpha^\rmt},\quad 
\uptheta \circ \upchi[\bfF_\alpha] \cong \bfF_{\alpha^\rmc}
,\text{ and }\\    
&\label{twist isom 2}\mRDIF{\alpha} \cong \uptheta \circ \upchi[\mDIF{\alpha}], \quad \mRESF{\alpha} \cong \uptheta \circ \upchi[ \mESF{\alpha}], \quad \mRQS{\alpha,\calC} \cong \uptheta \circ \upchi[\mYQS{\alpha,\calC}] .
\end{align}
Furthermore, we have a sequence of injective $H_n(0)$-module homomorphisms
\begin{equation}\label{Eq: Diagram of series I and II}  
\begin{tikzcd}
\bfF_{\alpha^\rmc} \arrow{r}{} &
\mRQS{\alpha,\calC} \arrow{r}{} & \mRESF{\alpha} \arrow{r}{} & \mRDIF{\alpha} \arrow{r}{} & \bfP_{\alpha^\rmt}
\end{tikzcd}. 
\end{equation}
\end{lemma}
\begin{proof}
The isomorphisms in \eqref{twist isom 1} and \eqref{twist isom 2} can be found in \cite[Table 2]{22JKLO} and \cite[Section 4]{22JKLO}, respectively.
And, the series in \eqref{Eq: Diagram of series I and II} is obtained from \eqref{Eq: sequence of surjective homomorphism} by applying the duality functor $\bfT^+_{\uptheta} \circ \bfT^-_{\upchi}$. 
\end{proof}

It was shown in \cite[Table 1]{22JKLO} that  
\( \uptheta \circ \upchi[\sfB([\sigma, \rho]_L)] \cong \sfB([\rho w_0, \sigma w_0]_L) \) as $H_n(0)$-modules.  
Since the modules in \eqref{Eq: sequence of surjective homomorphism} are, up to isomorphism, of the form \( \sfB([w_0(\set(\alpha)^\rmc), -]_L) \), this result implies that the modules in \eqref{Eq: Diagram of series I and II} are, up to isomorphism, of the form \( \sfB([-, w_1(\set(\alpha^\rmr))]_L) \).
We are now ready to state the main result of this subsection.

\begin{proposition}\label{twisted modules}
Let $\alpha$ be a composition of $n$.
Let $\scrD:=D(\mDIF{\alpha},\readingRLBT)$,  $\scrE:=D(\mESF{\alpha},\readingRLBT)$, and $\scrC:=D(\mYQS{\alpha,\calC},\readingRLBT)$.
Then we have the following.

\begin{enumerate}[label = {\rm (\alph*)}, itemsep = 0.3em]
\item 
$C_{\mRDIF{\alpha}}=\{
I w_0 \mid I \in C_{\mDIF{\alpha}}\}$, and $(C_{\mRDIF{\alpha}},\preceq) \cong (\ST((\scrD^\rmt)^x),\ble)$ as posets.

\item 
$C_{\mRESF{\alpha}} = \{I w_0 \mid I \in C_{\mESF{\alpha}}\}$, and $(C_{\mRESF{\alpha}},\preceq) \cong (\ST((\scrE^\rmt)^x),\ble)$ as posets.

\item 
$C_{\mRQS{\alpha,\calC}} = \{I w_0 \mid I \in C_{\mYQS{\alpha,\calC}}\}$, and 
$(C_{\mRQS{\alpha}},\preceq) \cong (\ST((\scrC^\rmt)^x),\ble)$ as posets.
\end{enumerate}
Here, $I w_0$ denotes the set $\{\zeta w_0 \mid \zeta \in I\}$.
\end{proposition}
\begin{proof}
By combining the isomorphism \( \uptheta \circ \upchi[\sfB([\sigma, \rho]_L)] \cong \sfB([\rho w_0, \sigma w_0]_L) \) with \eqref{twist isom 2}, we deduce the following:
\[
C_{\mRDIF{\alpha}} = \{I w_0 \mid I \in C_{\mDIF{\alpha}}\} \quad C_{\mRESF{\alpha}} = \{I w_0 \mid I \in C_{\mESF{\alpha}}\} \quad C_{\mRQS{\alpha, \calC}} = \{I w_0 \mid I \in C_{\mYQS{\alpha, \calC}}\}
\]

Additionally, from \cite[Theorem 3.6]{24CKO} it follows that \( \uptheta \circ \upchi (\bfM_P) \cong \bfM_{\overline{P}} \) for \( P \in \poset{n} \).  
Moreover, the proof of \cref{Prop: two kinds of P_D and intervals}(a) establishes that \( \overline{P}_{F_{D}^\downarrow} = P_{F_{D^\rmt}^\rightarrow} \) for all \( D \in \mathfrak{D}_n \).  
These results, combined with \cref{Prop: The case of mDIF}, \cref{Prop: Case where Y = X}, and \cref{Prop: Case where Y = hS}, confirm the desired poset isomorphisms.
\end{proof}

The following corollary follows from \cref{Coro: injective hull of V}, \cref{Coro: injective hull of X} and \cref{Coro: injective hull of hS}.
\begin{corollary}\label{covers of twisted modules}
Let $\alpha$ be a composition of $n$.
Let $\scrD:=D(\mDIF{\alpha},\readingRLBT)$,  $\scrE:=D(\mESF{\alpha},\readingRLBT)$, and $\scrC:=D(\mYQS{\alpha,\calC},\readingRLBT)$.
Then we have the following.
\begin{enumerate}[label = {\rm (\alph*)}, itemsep=0.3em]
\item 
$\sfB([w_1(\set(\sfr(\scrD))^\rmt), w_1(I^\rmt)]_L)$ is a projective cover of $\mRDIF{\alpha}$. 
Here, $I$ denotes the set in \cref{Coro: injective hull of V}(a).

\item $\sfB([w_0(\set(\sfr(\scrE))^\rmt), w_1(\set(\sfr(\scrE))^\rmt)]_L)$ is a projective cover of $\mRESF{\alpha}$.

\item 
$\sfB([w_0(\set(\sfr(\scrC))^\rmt), w_1(\DesLR{L}{\readingBTLR(F_\scrC^\rightarrow)}^\rmt)]_L)$ is a projective cover of $\mRQS{\alpha,\calC}$.
\end{enumerate}
\end{corollary}
\medskip

\subsection{Upper descent intervals from submodules of projective indecomposable \texorpdfstring{$H_n(0)$}{Hn0}-modules II}
\label{Upper descent intervals from submodules II} 

In this subsection, we discuss the submodules of a projective indecomposable \( H_n(0) \)-module that are related to the representation theory of \( 0 \)-Hecke--Clifford algebras.
 
Assume that \(\alpha\) is a peak composition of \(n\), meaning that \(\alpha_i \neq 1\) for \(1 \leq i < \ell(\alpha)\). 
A \emph{standard peak immaculate tableau} (SPIT) of shape \(\alpha\) is a standard immaculate tableau (SIT) \(\mathcal{T}\) of shape \(\alpha\) such that, for each \(1 \leq k \leq n\), the subdiagram of \(\tcd(\alpha)\) consisting of boxes filled with entries \(\leq k\) forms the diagram of a peak composition.
Let \(\SPIT(\alpha)\) denote the set of SPITs of shape \(\alpha\).
In \cite{22Searles}, Searles introduced a \(0\)-Hecke--Clifford module structure on \(\SPIT(\alpha)\).
Here, however, we consider only the \(0\)-Hecke module structure on \(\SPIT(\alpha)\).

Recall the module $\HnZmodR_{\alpha}$, where the \( H_n(0) \)-action on the \(\mathbb{C}\)-span of \(\SPIT(\alpha)\) is defined as follows: for each \( i = 1,2,\ldots,n-1 \) and \(\calT \in \SPIT(\alpha)\),
\begin{align}\label{eq: action on SIT}
\opi_i \cdot \calT = \begin{cases}
-\calT & \text{if \( i \) is weakly above \( i+1 \) in \(\calT\)},\\
0 & \text{if \( i \) and \( i+1 \) are in the first column of \(\calT\)},\\
s_i \cdot \calT & \text{otherwise}.
\end{cases}
\end{align}
We refer to \cite[Section 3.2]{24CNO} for the details.
\footnote{
It was shown in \cite{22Searles} that the image of the \(0\)-Hecke--Clifford module $\HnZmodR_{\alpha}\uparrow_{H_n(0)}^{HCl_n(0)}$ under the peak quasisymmetric characteristic is given by  
$
\QSQ{\alpha} = \sum_{\mathcal{T} \in \SPIT(\alpha)} K_{\Peak(\DesLR{L}{\readingLRTB(\mathcal{T})})}$,
where $K_{\Peak(\DesLR{L}{\readingLRTB(\mathcal{T})})}$ is the peak quasisymmetric function associated to $\Peak(\DesLR{L}{\readingLRTB(\mathcal{T})}$. 
This function is referred to as the \emph{quasisymmetric Schur \(Q\)-function} indexed by \(\alpha\).
We refer to \cite[Section 2]{24CNO} for the undefined notations.}

Let $\sinkSPIT{\alpha}$ be the SPIT of shape $\alpha$ whose entries in column $1$ are the first $\ell(\alpha)$ odd numbers, whose entries in column $2$ are the first $\ell(\alpha)-1$ or $\ell(\alpha)$ even numbers (depending on whether or not the last part of $\alpha$ is equal to 1), and whose entries in subsequent rows from top to bottom are the remaining numbers, increasing consecutively up each row.
Then the following lemma shows that $\RWrM{\readingLRTB}{\HnZmodR_{\alpha}}$ forms an upper descent interval.

\begin{lemma}{\rm (\cite[Lemma 3.18]{24CNO})}
\label{Lem: row reading word is interval}
Let $\alpha$ be a peak composition of $n$.
Then we have
\[
\HnZmodR_{\alpha} \cong \osfB([\readingLRTB(\sinkSPIT{\alpha}),w_1(\set(\alpha)^\rmr)]_L) \quad \text{as $H_n(0)$-modules}.
\]
\end{lemma}

Let $D(\HnZmodR_{\alpha},\readingLRTB)$ be the diagram whose $j$th row is
\[
\begin{cases}
\{(j+x-1,j) \mid 1 \leq x \leq \alpha_{\ell(\alpha)}\} & \text{if } j = \ell(\alpha),\\ 
\{(j,j), (j+1,j)\} \cup \{(k_j+x,j) \mid 3 \leq x \leq \alpha_j\} & \text{if $j < \ell(\alpha)$.}
\end{cases}
\]
Here, $k_j:=(\ell(\alpha)-1)+\sum_{j < r \leq \ell(\alpha)}(\alpha_{r}-2)$
for $1 \leq j \leq \ell(\alpha)-2$.
For example, if $\alpha = (3,2,4,2)$ and $\beta = (3,2,3,1)$, then the corresponding diagrams $D(\HnZmodR_{\alpha},\readingLRTB)$ and $D(\HnZmodR_{\beta},\readingLRTB)$ are
\begin{align*}
&\{(1,1),(2,1),(2,2),(3,2),(3,3),(4,3),(4,4),(5,4),(6,3),(7,3),(8,1)\}  \quad \text{and} \\ 
&\{(1,1),(2,1),(2,2),(3,2),(3,3),(4,3),(4,4),(5,3),(6,1)\},
\end{align*}
respectively. 
Graphically, these diagrams are represented as
\begin{equation}\label{Diagram: mHQSQ of 423}
\begin{tikzpicture}[baseline=10mm]
\foreach \c in {0,1,7}{
    \filldraw[color=black!15] (\hp*\c,\vp*0) rectangle (\hp*\c+\hp,\vp*1);
    \draw (\hp*\c,\vp*0) rectangle (\hp*\c+\hp,\vp*1);
}
\foreach \c in {1,2}{
    \filldraw[color=black!15] (\hp*\c,\vp*1) rectangle (\hp*\c+\hp,\vp*2);
    \draw (\hp*\c,\vp*1) rectangle (\hp*\c+\hp,\vp*2);
}
\foreach \c in {2,3,5,6}{
    \filldraw[color=black!15] (\hp*\c,\vp*2) rectangle (\hp*\c+\hp,\vp*3);
    \draw (\hp*\c,\vp*2) rectangle (\hp*\c+\hp,\vp*3);
}
\foreach \c in {3,4}{
    \filldraw[color=black!15] (\hp*\c,\vp*3) rectangle (\hp*\c+\hp,\vp*4);
    \draw (\hp*\c,\vp*3) rectangle (\hp*\c+\hp,\vp*4);
}
\draw[blue,->] (0,0) -- (\hp*9,0) node[anchor=west]{\textcolor{blue}{\footnotesize $x$}};
\draw[blue,->] (0,0) -- (0,\vp*5) node[anchor=south]{\textcolor{blue}{\footnotesize $y$}};
\end{tikzpicture}
\quad \text{and } \qquad 
\begin{tikzpicture}[baseline=10mm]
\foreach \c in {0,1,5}{
    \filldraw[color=black!15] (\hp*\c,\vp*0) rectangle (\hp*\c+\hp,\vp*1);
    \draw (\hp*\c,\vp*0) rectangle (\hp*\c+\hp,\vp*1);
}
\foreach \c in {1,2}{
    \filldraw[color=black!15] (\hp*\c,\vp*1) rectangle (\hp*\c+\hp,\vp*2);
    \draw (\hp*\c,\vp*1) rectangle (\hp*\c+\hp,\vp*2);
}
\foreach \c in {2,3,4}{
    \filldraw[color=black!15] (\hp*\c,\vp*2) rectangle (\hp*\c+\hp,\vp*3);
    \draw (\hp*\c,\vp*2) rectangle (\hp*\c+\hp,\vp*3);
}
\foreach \c in {3}{
    \filldraw[color=black!15] (\hp*\c,\vp*3) rectangle (\hp*\c+\hp,\vp*4);
    \draw (\hp*\c,\vp*3) rectangle (\hp*\c+\hp,\vp*4);
}
\draw[blue,->] (0,0) -- (\hp*7,0) node[anchor=west]{\textcolor{blue}{\footnotesize $x$}};
\draw[blue,->] (0,0) -- (0,\vp*5) node[anchor=south]{\textcolor{blue}{\footnotesize $y$}};
\node[right] at (\hp*9,\vp*-0.5) {.};
\end{tikzpicture}
\end{equation}
From \cref{Prop: two kinds of P_D and intervals} and \cref{Lem: row reading word is interval} it follows that 
\begin{equation}\label{Eq: RWR = PFD in case of Q}
\RWrM{\readingLRTB}{\HnZmodR_{\alpha}}=\SGL{P_{F_{D(\HnZmodR_{\alpha},\readingLRTB)}^\rightarrow}}.    
\end{equation}

For an SPIT $\calT$, we let $\readingdr(\calT)$ the reading word obtained from $\calT$ by reading the entries at $(1,j), (2,j-1)$ for $1 \leq j \leq \ell(\alpha)$ in the first two columns, and then continuing to read the entries top to bottom in each subsequent row, with the rows read in order from left to right.
For a positive integer $k$, let $2[k]:=\{2,4,\ldots,2k\}$.

\begin{theorem}\label{Thm: Class of Qr}
Let $\alpha$ be a peak composition of $n$.
Let $C_{\HnZmodR_{\alpha}}$ be the equivalence class of $[\readingLRTB(\sinkSPIT{\alpha}) 
,w_1(\set(\alpha^\rmr))]_L$ and $\scrG:=D(\HnZmodR_{\alpha},\readingLRTB)$.
\begin{enumerate}[label = {\rm (\alph*)}]
\item
$C_{\HnZmodR_{\alpha}} = \{\SGL{P_{T^x}} \mid T \in \ST(\scrG^x)\}$.
Furthermore, 
\begin{align*}
\min C_{\HnZmodR_{\alpha}} = [w_0(2[\ell(\alpha)-1]),\readingdr(\sourceSPIT{\alpha})]_L\quad \text{and}  \quad 
\max C_{\HnZmodR_{\alpha}} = [\readingLRTB(\sinkSPIT{\alpha}),w_1(\set(\alpha^\rmr))]_L.
\end{align*}

\item 
$C_{\HnZmodR_{\alpha}} \cong (\ST(\scrG^x),\ble)$ as posets.
\end{enumerate}
\end{theorem}
\begin{proof}
(a) By the construction of the diagram \( \scrG \), it is clear that \( \scrG \) is free of a strictly upper-right configuration. 
Therefore, by \cref{Thm: in case of Fsw = Fdown}, we have the first assertion that
\[
C_{\HnZmodR_{\alpha}} = \{\SGL{P_{T^x}} \mid T \in \ST(\scrG^x)\}.
\]

For the second assertion, let \(\min C_{\HnZmodR_{\alpha}} = [\sigma_0, \rho_0]_L\). 
Since \(\set(\sfc(\scrG)) = \{1, 3, \ldots, 2\ell(\alpha) - 1, 2\ell(\alpha), \ldots, n - 1\}\), it follows from \cref{Thm: in case of Fsw = Fdown} that
\[
\sigma_0 = w_0(\set(\sfc(\scrG))^\rmc) = w_0(2[\ell(\alpha) - 1]).
\]
For \(\rho_0\), \cref{Thm: Descriptions of minC_D and maxC_D} and \cref{Prop: two kinds of P_D and intervals} imply that $\rho_0 = \readingTBLR(T_{\scrG})$,
where \(T_{\scrG}\) is the standard tableau on \(\scrG\) (see \eqref{def of TD and TD prime}). 
Writing \(\readingTBLR(T_{\scrG})\) explicitly, we have
\begin{align*}
\readingTBLR(T_{\scrG}) = &
1 \ 
\underbrace{\alpha_1 + 1 \ 2}_{\text{col. 2}} \ 
\cdots \ 
\underbrace{\bfs_{\ell(\alpha)-1}(\alpha) + 1 \ \bfs_{\ell(\alpha)-2}(\alpha) + 2}_{\text{col. $\ell(\alpha)$}} \ 
\underbrace{\bfs_{\ell(\alpha)-1}(\alpha)+2 \ \bfs_{\ell(\alpha)-1}(\alpha)+3 \ \cdots \  n}_{\text{row $\ell(\alpha)$}} \\ 
& \underbrace{\bfs_{\ell(\alpha)-2}(\alpha) + 2 \ \cdots \ \bfs_{\ell(\alpha)-1}(\alpha)}_{\text{row $\ell(\alpha) - 1$}} \ 
\cdots \ 
\underbrace{3 \ 4 \ \cdots \ \alpha_1}_{\text{row 1}}.
\end{align*}
Here, $\bfs_j(\alpha):=\sum_{1 \leq r \leq j} \alpha_r$ for $1 \leq j < \ell(\alpha)$. 
From the definitions of \(\readingdr\) and \(\sourceSPIT{\alpha}\), we observe that the entries in column \(i\) (\(2 \leq i \leq \ell(\alpha)\)) correspond to \(\sourceSPIT{\alpha}(1, i)\) and \(\sourceSPIT{\alpha}(2, i - 1)\), while the entries in row $i$ ($1 \leq i \leq \ell(\alpha)$) correspond to the entries in row $i$ of $\sourceSPIT{\alpha}$. 
This gives that $\readingTBLR(T_{\scrG}) = \readingdr(\sourceSPIT{\alpha})$.
Thus we have
\[
\min C_{\HnZmodR_{\alpha}} = [w_0(2[\ell(\alpha) - 1]), \readingdr(\sourceSPIT{\alpha})]_L.
\]

For the third assertion, let \(\max C_{\HnZmodR_{\alpha}} = [\sigma_1, \rho_1]_L\). 
By \cref{Thm: Descriptions of minC_D and maxC_D} and \cref{Prop: two kinds of P_D and intervals}, we have $\sigma_1 = \readingBTLR(T'_{\scrG})$,
where \(T'_{\scrG}\) is the standard tableau on $\scrG$. 
Explicitly, \(\readingBTLR(T'_{\scrG})\) is written as
\[
\underbrace{k \ k + 1 \cdots k + \alpha_{\ell(\alpha)} - 1}_{\text{row $\ell(\alpha)$}} \ 
\underbrace{k - 2 \ k - 1 \ k + \alpha_{\ell(\alpha)} \ k + \alpha_{\ell(\alpha)} + 1 \ \cdots}_{\text{row $\ell(\alpha) - 1$}} 
\cdots 
\underbrace{1 \ 2 \ \cdots n}_{\text{row 1}},
\]
where \(k := 2\ell(\alpha) - 1\). 
Observing the correspondence with \(\sinkSPIT{\alpha}\), we have that $
\readingBTLR(T'_{\scrG}) = \readingLRTB(\sinkSPIT{\alpha})$.
For \(\rho_1\), since \(\set(\sfr(\scrG)) = \{\alpha_{\ell(\alpha)}, \alpha_{\ell(\alpha)} + \alpha_{\ell(\alpha) - 1}, \ldots, \alpha_{\ell(\alpha)} + \cdots + \alpha_1\}\), it follows from \cref{Thm: in case of Fsw = Fdown} that $\rho_1 = w_1(\set(\sfr(\scrG))) = w_1(\set(\alpha^\rmr))$.
Thus, we have
\[
\max C_{\HnZmodR_{\alpha}} = [\readingLRTB(\sinkSPIT{\alpha}), w_1(\set(\alpha^\rmr))]_L.
\]

(b) By the construction of \( \scrG \), it is clear that \( \scrG \) is free of a strictly upper-right configuration. Therefore, by \cref{Thm: Descriptions of minC_D and maxC_D} and \cref{Thm: in case of Fsw = Fdown}, we have
\[
(C_{\HnZmodR_{\alpha}}, \preceq) \cong (\ST(\scrG^x), \ble) \quad \text{as posets}.
\]
\end{proof}

\begin{example}
Let \(\alpha = (3,2,3,1) \models 9\). 
Consider the tableaux
\[
\ctab{\sinkSPIT{\alpha}}{
7 \\ 
5 & 6 & 8 \\ 
3 & 4 \\ 
1 & 2 & 9
}
\quad \text{and} \quad
\ctab{\sourceSPIT{\alpha}}{
9 \\ 
6 & 7 & 8 \\ 
4 & 5 \\ 
1 & 2 & 3
}.
\]
Using \cref{Thm: Class of Qr}, we compute the following:  
\begin{align*}
&\min C_{\HnZmodR_{\alpha}} = [w_0(\{2,4,6\}), 142659783]_L \quad \text{and} \\
&\max C_{\HnZmodR_{\alpha}} = [756834129, w_1(\{1,4,6\})]_L.
\end{align*}
Thus, \(C_{\HnZmodR_{\alpha}}\) is the equivalence class of \([756834129, w_1(\{1,4,6\})]_L\).

On the other hand, the diagram \(D_{756834129; \{1,4,6\}}\) corresponds to the diagram \(\scrG\) shown in \eqref{Diagram: mHQSQ of 423}. 
Notably, \(\scrG\) is free of a strictly upper-right configuration. 
Given the standard fillings 
\[
F_\scrG^\downarrow = 
\ctabT{
\empty & \empty & \empty & 6 \\ 
\empty & \empty & 4 & 7 & 8  \\ 
\empty & 2 & 5 \\ 
1 & 3 & \empty & \empty & \empty & 9
}
\quad \text{and} \quad 
F_\scrG^\rightarrow = 
\ctabT{
\empty & \empty & \empty & 1 \\ 
\empty & \empty & 2 & 3 & 4  \\ 
\empty & 5 & 6 \\ 
7 & 8 & \empty & \empty & \empty & 9
},
\]
it follows from \cref{Thm: in case of Fsw = Fdown} and \cref{Rem: reading permutations from F} that 
we confirm 
\begin{align*}
&\min C_{\HnZmodR_{\alpha}} = [132547689, 142659783]_L \quad \text{and} \\ 
&\max C_{\HnZmodR_{\alpha}} = [756834129, 967845123]_L.
\end{align*}

To further explore \(C_{\HnZmodR_{\alpha}}\), we note that \(P_{F^x}\) corresponds to all the posets obtained from \(F \in \ST(\scrG^x)\) by \cref{Thm: Class of Qr}(b). 
By computational enumeration, it is found that there are precisely 594 elements in \(\ST(\scrG^x)\).
\end{example}

We close this section by providing the projective cover and the injective hull of $\HnZmodR_{\alpha}$.

\begin{theorem}\label{Thm: Proj. and Inj. of Qr}
Let $\alpha$ be a peak composition of $n$ and let $\scrG:=D(\HnZmodR_{\alpha},\readingLRTB)$.

\begin{enumerate}[label = {\rm (\alph*)}]
\item $
\osfB([w_0(2[\ell(\alpha)-1]),w_1(\DesLR{L}{\readingLRBT(F_\scrG^\downarrow)})]_L)$ is a projective cover of $\HnZmodR_{\alpha}$.

\item $
\osfB([w_0(\set(\alpha^\rmr)),w_1(\set(\alpha^\rmr))]_L)$ is an injective hull of $\HnZmodR_{\alpha}$.
In particular, it is a projective indecomposable module.
\end{enumerate}
\end{theorem}
\begin{proof}
Since \(\scrG\) is free of a strictly upper-right configuration, it follows from \cref{how to find injective hull} that  
\[
\osfB\big([w_0(\set(\sfc(\scrG))^\rmc), w_1(\DesLR{L}{\readingLRBT(F_\scrG^\downarrow)})]_L\big)
\]
is a projective cover of \(\HnZmodR_{\alpha}\). 
Observing that  
\[
\set(\sfc(\scrG)) = \{1, 3, \ldots, 2\ell(\alpha) - 1, 2\ell(\alpha), \ldots, n - 1\},
\]  
we deduce the assertion (a).

Similarly, using \cref{how to find injective hull}, we also find that  
\[
\osfB\big([w_0(\DesLR{L}{\readingBTLR(F_\scrG^\rightarrow)}), w_1(\set(\sfr(\scrG)))]_L\big)
\]
is an injective hull of \(\HnZmodR_{\alpha}\). 
Since  
\[
\DesLR{L}{\readingBTLR(F_\scrG^\rightarrow)} = \{\alpha_{\ell(\alpha)}, \alpha_{\ell(\alpha)} + \alpha_{\ell(\alpha)-1}, \ldots, \alpha_{\ell(\alpha)} + \cdots + \alpha_2\} = \set(\sfr(\scrG)),
\]  
the assertion (b) follows. 
Moreover, this injective hull coincides with a projective indecomposable module.
\end{proof}

As an important consequence of \cref{Thm: Proj. and Inj. of Qr}(b), we derive the following corollary.

\begin{corollary}\label{Coro: indecomposability of Qr}
Let $\alpha$ be a peak composition of $n$.
Then the $H_n(0)$-module $\HnZmodR_{\alpha}$ is indecomposable.
\end{corollary}
\vspace*{5mm}

\section{Further avenues}

\begin{enumerate}
\item
In \cref{Lem: regular poset and indirectly comparable pair} (c), it is shown that $s_i \cdot P \Ceq P $. This naturally leads to the question of whether 
$s_i \cdot P \overset{M}{\not\simeq} P$ holds. Addressing this question would be of significant interest.

\item
A natural direction for future work is to extend the results of \cref{Sec: The equivalence classes of weak Bruhat Intervals} to Coxeter groups of types $B$ and $D$. For example, the following questions arise:

\begin{enumerate}
\item Given an equivalence class of weak Bruhat intervals defined by lower and upper descent sets, is there a family of diagrammatic objects or fillings that are in bijection with the elements of the class? 

\item Can one formulate a result on injective hulls and projective covers analogous to \cref{how to find injective hull}? 
\end{enumerate}
\end{enumerate}

\vspace*{10mm}

\bibliographystyle{abbrv}

\end{document}